\RequirePackage{snapshot}
\documentclass[onefignum,onetabnum]{siamart171218}

\usepackage[english]{babel}

\usepackage{mathtools}

\usepackage{amssymb}  

\usepackage{cancel}

\usepackage{commath}
\usepackage{graphicx}
\usepackage{booktabs}
\usepackage{placeins} 
\usepackage{framed}
\usepackage[normalem]{ulem}
\usepackage{enumerate}
\usepackage{cleveref}
\usepackage{algorithmic}

\ifpdf
  \DeclareGraphicsExtensions{.eps,.pdf,.png,.jpg}
\else
  \DeclareGraphicsExtensions{.eps}
\fi

\newsiamremark{remark}{Remark}
\newsiamremark{hypothesis}{Hypothesis}
\crefname{hypothesis}{Hypothesis}{Hypotheses}
\newsiamthm{claim}{Claim}

\headers{HiPhom\texorpdfstring{$\varepsilon$}- for advection diffusion reaction problems}{G. Conni, S. Piccardo, S. Perotto, G.M. Porta, M. Icardi}

\title{HiPhom\texorpdfstring{$\pmb{\varepsilon}$}-: HIgh order Projection-based HOMogenisation for advection diffusion\\ reaction problems \thanks{Submitted to the editors September 25, 2023.
}
}

\author{Giovanni Conni\thanks{Department of Computer Science, KU Leuven, Celestijnenlaan 200A, 3001, Leuven, Belgium}; {Department of Mechanical Engineering, KU Leuven, Celestijnenlaan 300b,
3001 Leuven, Belgium} (\email{giovanni.conni@kuleuven.be})
	\and Stefano Piccardo\thanks{CERMICS, Ecole des Ponts, 77455 Marne-la-Vallée, France}; {Inria, 2 rue Simone Iff, 75589 Paris, France};
	 {LaC\`aN, Universitat Polit\`ecnica de Catalunya, 08034 Barcelona, Spain} (\email{stefano.piccardo@enpc.fr})
	\and Simona Perotto$^*$\thanks{MOX-Department of Mathematics, Politecnico di Milano, Piazza Leonardo da Vinci 32, I-20133, Milano, Italy} (\email{simona.perotto@polimi.it})
	\and Giovanni Michele Porta\thanks{Department of Civil and Environmental Engineering, Politecnico di Milano, Piazza Leonardo da Vinci 32, I-20133, Milano, Italy} (\email{giovanni.porta@polimi.it})
	\and Matteo Icardi\thanks{School of Mathematical Sciences, University of Nottingham, NG7 2RD, UK (\email{matteo.icardi@nottingham.ac.uk})}
	}

\usepackage{amsopn}

   \newcommand{\Rey}{\operatorname{\mathit{R\kern-.04em e}}}
    \newcommand{\Damk}{\operatorname{\mathit{D\kern-.06em a}}}
    \newcommand{\Pecl}{\operatorname{\mathbb P e}} 

    \newcommand*\dd{\mathop{}\!\mathrm{d}}
    
\ifpdf
\hypersetup{
  pdftitle={HiPhom\texorpdfstring{${\varepsilon}$}-},
  pdfauthor={G. Conni, S. Piccardo, S. Perotto, G.M. Porta, M. Icardi}
}
\fi

\begin{document}

\maketitle

\begin{abstract}
We propose a new model reduction technique for multiscale scalar transport problems that exhibit dominant axial dynamics. 
To this aim, we rely on the separation of variables to combine a Hierarchical Model (HiMod) reduction with a two-scale asymptotic expansion. We extend the two-scale asymptotic expansion to an arbitrary order and exploit the high-order correctors to define the HiMod modal basis which approximates the transverse dynamics of the flow, while we adopt a finite element discretisation to model the leading stream. 
The resulting method, which is named HiPhom$\varepsilon$ (HIgh-order Projection-based HOMogEnisation), is successfully assessed both in steady and unsteady advection-diffusion-reaction settings. The numerical results confirm the very good performance of HiPhom$\varepsilon$ which improves the accuracy and the convergence rate of HiMod and extends the reliability of the standard homogenised solution to transient and pre-asymptotic regimes.
\end{abstract}



\section{Introduction}
In the last decades, computational science has progressed in becoming a reliable prediction tool for the analysis of real-world phenomena.
To solve complex problems accurately yet efficiently, numerical methods often need to resort to strategies that can reduce model complexity, and thus the computational burden, while carefully balancing the computational power and fidelity to a high-complexity reference model. Model Order Reduction (MOR) methods include a wide class of numerical and analytical methodologies \cite{chinesta2011short,HRS16,qu2013model,QMN16} which achieve this goal, for instance, by reducing the dimension of the state space of the problem, leading to an improvement in the computational performance with respect to classical numerical methods, such as finite elements or finite volumes.
This approach has proven successful in reducing the computational burden associated with complex dynamical systems, e.g., in fluid dynamics applications, high-dimensional inverse problems and uncertainty quantification \cite{Benner,Hinze,Lassila14}.

Efficient numerical simulations and MOR techniques turn out to be particularly relevant in the presence of dynamic processes acting on space and time scales much smaller than those tied to the final quantities of interest. In such cases, a Direct Numerical Simulation of such fine-scale dynamics would constrain the resolution of the space-time discretisation to a level at which the computational cost becomes unaffordable. This is the case, for example, of porous media~\cite{bear1968physical}, where the fine-scale features have characteristic spatial dimensions much smaller than the ones associated with the problem of interest (e.g., single pore versus aquifer/reservoir scales). These problems are typically tackled by resorting to averaged models, that implicitly represent pore-scale features through effective parameters. This unavoidably leads to a loss of information which, sometimes, may limit the modelling properties of averaged models.

In this work, we propose a new numerical approach that is aimed at capturing efficiently small- and large-scale dynamics, combining MOR with a proper multiscale formulation. We consider, as a first application, standard advection-diffusion-reaction equations. Our approach is based on the consideration that often the advective field is characterized by a preferential direction and yet exhibits smaller but non-negligible local fluctuations and variations in intensity along the other directions. This is the case, for example, of geophysical flows, such as in rivers and aquifers, or fluid dynamics in pipes and blood vessels. In these cases, transport exhibits large scale dynamics along the flow stream direction and fine scale transverse dynamics in the directions perpendicular to the flow. Thanks to the multiscale structure, the problem can be tackled by a two-scale expansion or a volume-averaging technique to derive macroscopic transport equations where the effect of the velocity spatial variations are taken into account through an effective diffusion term \cite{auriault1995taylor,tesiPop, municchi2020macroscopic}, known as hydrodynamic (or Taylor) dispersion \cite{taylor1953dispersion}. This upscaling is based on retaining the lowest order terms in a two-scale asymptotic expansion.
Extensions of upscaling approaches, that can accurately approximate pre-asymptotic and high-order features such as mixing and reactions (see, e.g., \cite{porta2015continuum,sund2017upscaling}), are of great relevance in many applications, such as chemically or biologically activated reactive processes, and can be coupled with solute transport \cite{bolster2011mixing}. In this context, the key challenge is to be able to capture not only the average concentration in the system but also local time-  and space-dependent fluctuations induced by small-scale features of the problem, e.g., the advective field spatial arrangement or the local geometry.

The method presented in this work stems from the Hierarchical Model (HiMod) reduction introduced in \cite{perotto2010hierarchical} and combines it with two-scale asymptotic expansions and homogenisation theory.
This combination allows us to take advantage of the mathematical structure of the problem at hand to build an optimal reduced setting.
In line with our objectives, HiMod focuses on problems where a dominant direction can be recognized. Model reduction is achieved by splitting the solution into longitudinal and transverse components, by employing a separation of variables. In the original formulation in~\cite{perotto2010hierarchical}, HiMod employs finite elements along the longitudinal axis and a spectral approximation along the transverse direction. For the spectral approximation, we propose here a problem-specific set of modal basis functions derived from the homogenisation procedure, tailored to approximating the transverse dynamics. The use of homogenisation and perturbation expansions to derive multiscale numerical methods has been extensively investigated in the literature across fluid and solid mechanics applications \cite{Abdulle2014,Efendiev2008,Geers2010,municchi2020macroscopic}. Here, we consider a formulation tailored for directional transport settings. Taking advantage of this particular setup, we analytically derive a set of recursive equations that can be used to obtain effective transverse modes for the spectral expansion employed by HiMod. 
This combined approach is named here HIgh order Projection-based HOMogEnisation, here introduced as HiPhom${\varepsilon}$ model reduction, where ${\varepsilon}$ is
the standard symbol used in perturbation-based expansions.
Through this approach, an enriched one-dimensional model is obtained, which is significantly less computationally intensive than the full two- or three-dimensional approximations, while being richer than a simple one-dimensional approximation, such as the one yielded by a standard two-scale homogenisation approach. Notably, HiPhom$\varepsilon$ approach includes the effects of the transverse dynamics through high-order modes allowing a very accurate description of the transitions towards the self-similar asymptotic long-time long-space concentration profiles.

The work is organized as follows. In section~\ref{sec:stArt} we introduce the two reduction approaches HiPhom$\varepsilon$ relies on, namely 
the HiMod formulation in section~\ref{sec:himod} and 
the homogenisation approach
in section~\ref{sec:tae}. The reference  setting is a standard scalar advection-diffusion-reaction problem.
In section~\ref{sec:hoasy} we 
generalize the two-scale expansion carried out in section~\ref{sec:tae} to an arbitrary order, and we explicitly prove   
a recursive formula to compute the physic-based modal basis characterizing a HiPhom$\varepsilon$ expansion.
We formalize
the HiPhom$\varepsilon$ approach in section~\ref{sec:hip}, while we numerically assess the associated performance in section~\ref{num_Sec}, both in a steady and in an unsteady regime. 
The outperformance of HiPhom$\varepsilon$ with respect to a standard HiMod reduction is verified for all the considered case studies. 

\section{Problem statement and background}\label{sec:stArt}
We identify the computational domain with a long and thin rectangle defined by 
\begin{equation} \label{eq:2}
{\Omega }_{\varepsilon} := \Big\{ {\left( {x,z} \right) \in \mathbb{R}^2\ :\ 0 < x < L,\; - \frac{l}{2}  < z < \frac{l}{2} } \Big\},
\end{equation}
where the length $L$ and the width $l$ of the domain are such that $\varepsilon:=l/L\ll1$. In particular, we distinguish the following three portions of the boundary, $\partial \Omega_\varepsilon$, of $\Omega_\varepsilon$
\begin{equation} \label{eq:3bis}
\begin{aligned}
&{{\Gamma}^i_{\varepsilon} }:=\Big\{ {\left( {x,z} \right) \in \mathbb{R}^2\ : \  x=0,\; - \frac{l}{2}  \leq z \leq \frac{l}{2} } \Big\},\\[2mm]
&{\Gamma}^w_{\varepsilon}:= \Big\{ \left( {x,z} \right) \in \mathbb{R}^2\ : \ 
0 < x < L,\;  z=-\frac{l}{2},\frac{l}{2}\Big\},\\[2mm]
&{\Gamma}^o_{\varepsilon}:= \Big\{ \left( {x,z} \right) \in \mathbb{R}^2\ : \
x=L,\; - \frac{l}{2}\leq z \leq \frac{l}{2} \Big\},
\end{aligned}
\end{equation}
coinciding with the inflow, the wall and the outflow boundary, respectively.\\ 
As a reference differential problem, we consider a standard unsteady scalar advection-diffusion-reaction (ADR) problem, namely
\begin{equation}\label{eq:6}
\left\{
\begin{array}{ll}
\partial_t c+\mathcal{L}_z c=f \quad &\mbox{in } {{\Omega}_{\varepsilon}}\times \left(0,T\right] \\[2mm]
D_\varepsilon \nabla_{\nu_{\varepsilon}} c=0  &\mbox{on } {{\Gamma}^w_{\varepsilon}}\times \left(0,T\right]\\[2mm]
c=c_B   &\mbox{on } {{\Gamma}^i_{\varepsilon}}\times \left(0,T\right]  \\[2mm]
 D_\varepsilon \nabla_{\nu_{\varepsilon}} c=0   &\mbox{on } {{\Gamma}^o_{\varepsilon}}\times \left(0,T\right]  \\[2mm]
c=\tilde c &\mbox{in } {\Omega}_{\varepsilon}\times \{t=0\},
\end{array}
\right.
\end{equation}
where $c=c(x,z,t)$ is the (unknown) solute concentration; $\mathcal{L}_z$ is the spatial differential operator defined by
\begin{equation}\label{eq:Lfunc1}
    \mathcal L_z := -\nabla\cdot(D_\varepsilon\nabla) + \pmb{u}\cdot \nabla + \sigma,
\end{equation}
with $D_\varepsilon\in \mathbb R $ the diffusion coefficient, 
$\pmb{u} = [u(z),0]^T$ the advective field, here assumed to coincide with an axial flow, 
$\sigma \in \mathbb R$ the reaction coefficient; $f \in \mathbb R$ denotes the forcing term; $c_B\in \mathbb R$ is the non-homogeneous Dirichlet data; $\tilde c\in \mathbb R$ identifies the initial value of the concentration, here taken identically null for simplicity; 
$\nabla_{\nu_{\varepsilon}}$ denotes the gradient along the normal direction to the boundary of ${\Omega}_{\varepsilon}$,
$\nu_{\varepsilon}$ being the unit outward normal vector to  $\partial {\Omega}_{\varepsilon}$.
%
In particular, we consider advection-dominated problems, so that the global P\'eclet number for mass transfer
\begin{equation}
    \Pecl:=\frac{\bar{u}L}{2 D_\varepsilon}  
\end{equation}
with $\overline{ {u}}$ the average axial velocity,
becomes larger and larger for $\varepsilon$ going to zero, where we 
have assumed both 
$\overline{{u}}$ and $L$ to be $\mathcal{O}\left( 1 \right)$ and 
$D_\varepsilon := \varepsilon D$, with $D=\mathcal{O}\left( 1 \right)$ (namely,
$D_\varepsilon = \mathcal{O}\left( {{\varepsilon}}  \right)$), so that $\Pecl=\mathcal{O}\left( {{\varepsilon ^{ - 1}}} \right)$.

\subsection{HiMod reduction}\label{sec:himod}
Hierarchical Model (HiMod) reduction applies to the weak form of the full problem of interest, i.e., with reference to problem \eqref{eq:6}, to the formulation
\begin{equation}\label{eq:weakform}
\forall t>0 \  \ \mbox{find } c(t)\in V \ : \ \big(\partial_t c(t), v \big)+a(c(t),v)=\mathcal{F}(v) \ \ \forall v\in V,
\end{equation}
with $c(0)=0$, and where $V$ is the Sobolev space, $H^1_{{\Gamma}^i_{\varepsilon}}({\Omega}_{\varepsilon})$, of the $H^1({\Omega}_{\varepsilon})$-functions with null trace on the inflow boundary~\cite{ErnGuermond}; $a(\cdot,\cdot):V\times V\rightarrow\mathbb{R}$ is the bilinear form given by
$$
a(c(t), v) = \displaystyle \int_{{\Omega}_{\varepsilon}} \Big[ D_{\varepsilon} \nabla c(t) \cdot \nabla v + \pmb{u} \cdot \nabla c(t)v+ \sigma c(t) v \Big]\, d{\Omega}_{\varepsilon};
$$
$\mathcal{F}(\cdot):V\rightarrow\mathbb{R}$ is the linear form defined by
$$
\mathcal{F}(v)=\displaystyle \int_{{\Omega}_{\varepsilon}} fv\, \mathrm d{\Omega}_{\varepsilon} - a(R_{c_B}, v)\, d{\Omega}_{\varepsilon},
$$
with $R_{c_B}$ the lifting of the boundary value $c_B$ in \eqref{eq:6}.
Suitable regularity hypotheses are advanced on the problem data in order to ensure the 
existence and the uniqueness of the solution to the weak form \eqref{eq:weakform} by the Lax Milgram lemma~\cite{ErnGuermond}.

\vspace*{0.2cm}

\paragraph{The domain geometry} HiMod reduction moves from the assumption that the computational domain exhibits a fibre-bundle layout. 
We observe that domain $\Omega_{\varepsilon}$ in \eqref{eq:2} is compliant with this requirement.
This geometric hypothesis leads us to distinguish between a one-dimensional (1D) supporting fibre, aligned with the main dynamics of the problem, and a set of transverse fibres, parallel to the secondary dynamics. Thus,
domain $\Omega_{\varepsilon}$ in \eqref{eq:6} can be characterized by the relation 
\begin{equation}\label{eq:domainHiMod_1}
    \Omega_{\varepsilon}=\bigcup_{x\in  \Omega^{1D}_{\varepsilon}}\{x\}\times\gamma_x,
\end{equation}
with $\Omega^{1D}_{\varepsilon}$ the supporting fibre, and $\gamma_x\subset\mathbb{R}$ the transverse fibre associated with the generic point $x\in \Omega^{1D}_{\varepsilon}$. 
For simplicity, fibre $\Omega^{1D}_{\varepsilon}$ is here taken rectilinear, i.e.,  $\Omega^{1D}_{\varepsilon}=[x_i,x_o]\subset\mathbb{R}$, where subscript $i$ and $o$ stands for inflow and  outflow, respectively. 
We refer to~\cite{perotto2014hierarchical,JCFD20} for the case of a bent supporting fibre. 

Successively, HiMod reduction simplifies computations by mapping the physical domain $\Omega_{\varepsilon}$ into a reference environment, $\widehat{\Omega}$. 
This is performed first by introducing, for any
$x\in \Omega^{1D}_{\varepsilon}$, the map
$\psi_x:\gamma_x\rightarrow\widehat{\gamma}$ between the transverse fibre $\gamma_x$ and the reference one $\widehat{\gamma}=[0, 1]$. 
Then, maps $\psi_x$ are exploited to define the global map $\Psi:\Omega_{\varepsilon}\rightarrow\widehat{\Omega}$, with $\widehat{\Omega}=\cup_{x\in\widehat{\Omega}^{1D}}\{x\}\times\widehat{\gamma}$, $\widehat{\Omega}^{1D}$ being the reference supporting fibre. In particular, we assume that map $\Psi$ does not modify the such a fibre, so that $\widehat{\Omega}^{1D}\equiv {\Omega}^{1D}_{\varepsilon}$ and the generic point $(x, z)$ in $\Omega_{\varepsilon}$ is associated with the point $(\widehat{x},\widehat{z})$ in $\widehat{\Omega}$, with $\widehat{x}=x$ and $\widehat{z}=\psi_x(z)$.
In addition, since we limit to the 2D case, map $\psi_x$ coincides with the linear transformation
\begin{equation} \label{eq:linpsi}
    \widehat{z}=\psi_x(z)=\frac{1}{|\gamma_x|}\, z,
\end{equation}
with $|\gamma_x|$ the length of the transverse fibre. 
Some regularity assumptions are advanced on maps $\psi_x$ and $\Psi$ (i.e., $\psi_x$ is assumed to be a $C^1$-diffeomorphism, while $\Psi$ is taken differentiable) in order to avoid the presence of irregularities (such as kinks) along $\partial \Omega_\varepsilon$.
Finally, the fibre decomposition characterising $\Omega_\varepsilon$ is inherited by the domain boundary so that the three boundary portions in \eqref{eq:3bis} can be defined as
\begin{equation}\label{eq:boundaries}
    {\Gamma}^i_{\varepsilon}=\{x_i\}\times\gamma_{x_i}, \quad {\Gamma}^w_{\varepsilon}=\bigcup_{x\in\Omega^{1D}_{\varepsilon}}\{x\}\times\partial\gamma_x,\quad{\Gamma}^o_{\varepsilon}=\{x_o\}\times\gamma_{x_o}, 
\end{equation}
and, analogously for the corresponding boundaries, $ \widehat{\Gamma}^i$, $ \widehat{\Gamma}^w$, $ \widehat{\Gamma}^o$, in $\widehat \Omega$.

\vspace*{0.2cm}

\paragraph{The HiMod space} HiMod reduction approximates the dynamics aligned with the main and with the transverse fibres by means of a different discretisation, according to a separation of variables paradigm. \\
Following the seminal papers~\cite{enumath,perotto2010hierarchical}, we employ a finite element discretisation along the $x$-axis, while resorting to a modal expansion to describe the dynamics parallel to the $z$-direction.
With this aim, we associate the discrete space $V_{1D}$ of the continuous piecewise linear polynomials  with a partition, ${\mathcal T}_h$, of the supporting fibre $\Omega^{1D}_{\varepsilon}$ into subintervals, and the basis $\{\varphi_k\}_{k\in\mathbb{N}}$ of modal functions which satisfy an $L^2(\widehat{\gamma})$-orthonormality condition with the reference transverse fibre.
Approximations alternative to finite elements can be adopted to discretise the mainstream dynamics (we refer, for instance, to~\cite{JCFD20} where the authors employ an isogeometric approximation of the leading dynamics for haemodynamic modelling in 3D patient-specific geometries). \\
The boundary conditions on ${\Gamma}^i_{\varepsilon}$ and ${\Gamma}^o_{\varepsilon}$ are properly included in the space $V_{1D}$, while boundary data on ${\Gamma}^w_{\varepsilon}$ have to be taken into account by the modal functions. In particular, the imposition of the conditions on the lateral boundary via functions $\varphi_k$ deserves particular care. To this aim, we adopt the \emph{educated basis} approach proposed in~\cite{aletti2018himod} and successfully validated in several practical contexts~\cite{guzzetti2018hierarchical,JCFD20}. The idea is to identify the modal functions with the solution to an auxiliary (1D) Sturm-Liouville eigenvalue problem
\begin{equation}\label{eq:SLE1}
    \mathcal{L}_s\varphi_k({\widehat{z}})=\lambda_k\varphi_k({\widehat{z}}),
\end{equation}
solved on the reference fibre $\widehat{\gamma}$ and endowed with the same boundary conditions imposed on ${\Gamma}^w_{\varepsilon}$. 
The modal basis $\varphi_k$ imposes in an essential way the data to be assigned on the wall boundary.
In general, the differential operator $\mathcal{L}_s$  in \eqref{eq:SLE1} coincides with the symmetric part of the operator $\mathcal L_z$ in \eqref{eq:Lfunc1}. In the sequel, for simplicity, we choose $\mathcal{L}_s=\partial_{zz}$.
A rigorous convergence analysis for the educated HiMod reduction is available in~\cite{aletti2018himod}.

Thus, the HiMod reduced space can be formalized as 
\begin{equation} \label{eq:vmspace}
    V_m =\bigg\{v_m(x,z)=\sum_{k=0}^{m-1}\widetilde{v}_k(x)\varphi_k\left(\psi_x(z)\right) :  \widetilde{v}_k\in V_{1D}, x\in\Omega^{1D}_{\varepsilon},  z\in\gamma_x\bigg\}.
\end{equation}
The modal index $m$ categorizes the space $V_m$ in the hierarchy of reduced spaces. In particular, the value for $m$ is set by the user according to an a priori knowledge of the phenomenon at hand or, rather, by means of an automatic selection driven, for instance, by an a posteriori modelling error analysis~\cite{perotto2014coupled,perotto2015space,perotto2014survey}. 
Moreover, the user can adopt the same number of modal functions along the whole domain (uniform HiMod reduction), or locally tune such a number in order to describe spatial heterogeneities without overloading the whole computational effort (we refer the interested reader to the piece-wise and to the point-wise HiMod formulations in~\cite{zilioproc,perotto2014survey,perotto2015space}).

Finally, conformity and spectral approximability hypotheses are introduced on the HiMod space $V_m$, together with a standard $H^1(\Omega^{1D}_{\varepsilon})$-density requirement on the discrete space $V_{1D}$ in order to make the HiMod formulation well-defined.

\vspace*{0.2cm}

\paragraph{The HiMod reduced problem} The HiMod formulation of problem \eqref{eq:weakform} can be stated as
\begin{equation}\label{eq:HiMod}
\forall t>0 \  \  \mbox{find } c_m(t)\in V_m \ : \ \big(\partial_t c_m(t), v_m \big)+a(c_m(t),v_m)=\mathcal{F}(v_m) \ \ \forall v_m\in V_m,
\end{equation}
where $c_m(0)=0$.
 The assumptions advanced on spaces $V_m$ and $V_{1D}$ ensure the well-posedness of problem \eqref{eq:HiMod}, together with the convergence of the HiMod solution $c_m(t)$ to the weak solution $c(t)$ in \eqref{eq:weakform}, for any $t>0$
(for more details, we refer to~\cite{perotto2010hierarchical}).

The computational advantage of the HiMod reduction becomes evident when considering the algebraic counterpart of formulation \eqref{eq:HiMod}. To this aim, we exploit the modal expansion 
\begin{equation}\label{eq:hieq}
    c_m(t)(x,z) =\sum_{k=0}^{m-1}\widetilde{c}_k(t)(x)\varphi_k(\psi_x(z))=\sum_{k=0}^{m-1}\sum_{s=1}^{N_h}\widetilde{c}_{ks}(t)\theta_s(x)\varphi_k(\psi_x(z))
\end{equation}
for the HiMod concentration in \eqref{eq:HiMod} in terms of the modal functions and of the finite element basis, $\{\theta_s\}_{s=1}^{N_h}$ with $N_h={\rm dim}(V_{1D})$, and we pick the test function $v_m$ as $\theta_i \varphi_j$, with $i=1,\ldots, N_h$ and $j=0,\ldots, m-1$. Concerning the time-dependence, we resort to a standard $\theta$-method, after discretising the time window $[0, T]$ with the times $\{ t^n\}$ (see~\cite{perotto2015space} for an alternative time discretisation based on finite elements). \\
Thus, at each $t^n$, we are led to solve a linear system of order $mN_h$ whose unknowns coincide with the coefficients $\widetilde{c}_{ks}^{\, n}=\widetilde{c}_{ks}(t^n)$, for $k=0, \ldots, m-1$ and $s=1, \ldots, N_h$. This is equivalent to state that, for each time $t^n$, the full order problem \eqref{eq:weakform} reduces to a system of $m$ coupled 1D problems, which are solved along the supporting fibre $\Omega^{1D}_{\varepsilon}$, and whose coefficients collect information along the transverse fibres $\gamma_x$.\\ 
In particular, the HiMod stiffness matrix exhibits a block-sparsity pattern to be advantageously exploited when solving the reduced system~\cite{perotto2010hierarchical}. 

We remark that the inclusion of the transverse dynamics in the HiMod model automatically promotes the HiMod approximation in terms of reliability when compared with standard 1D reduced models. In addition, the possibility to deal with 1D problems instead of higher-dimensional settings relieves the user from certain computationally demanding steps, such as the meshing of the computational domain (now limited to a 1D partition of $\Omega^{1D}_{\varepsilon}$), with a consequent meaningful gain in terms of efficiency, in particular when compared with full 3D models~\cite{aletti2018himod,guzzetti2018hierarchical,blanco2015hybrid,JCFD20}. 

\subsection{Homogenisation via a two-scale asymptotic expansion}\label{sec:tae}
In this section, we highlight how 
a standard homogenisation technique leads to separate variables, analogously to a HiMod reduction. This compliance will provide us the motivation to formalize the HiPhom$\varepsilon$ method.

According to a standard homogenisation process for axial flows~\cite{Hornu:97}, problem \eqref{eq:6} is decomposed into a \textit{macro-scale} contribution, associated with the main flow dynamics, i.e., parallel to the $x$-direction, and a \textit{micro-scale} term, which takes into account the local fluctuations along the perpendicular $z$-direction. By means of the multiscale parameter $\varepsilon$, the dependence of the solution $c$ in  \eqref{eq:6} on the macro- and on the micro-scale is explicitly taken into account by introducing the rescaled  transverse (or fast) variable $y:=\frac{z}{\varepsilon}$ so that now $c = c\left(x,y,t;\varepsilon\right)$.\\
Consistently, domain $\Omega_\varepsilon$ in \eqref{eq:2} and the boundaries in \eqref{eq:3bis} are redefined as
%
$$
\begin{array}{c} 
\Omega:=\left(0,L\right)\times\left(-\displaystyle\frac{L}{2},\frac{L}{2}\right), \\[5mm] \Gamma^i:=\left\{0\right\}\times\left(-\displaystyle\frac{L}{2},\frac{L}{2}\right), \quad \Gamma^w:=\left(0,L\right)\times\left\{-\displaystyle\frac{L}{2},\frac{L}{2}\right\}, \quad 
\Gamma^o:=\left\{L\right\}\times \left(-\displaystyle\frac{L}{2},\frac{L}{2}\right).
\end{array}    
$$
In this new coordinate system, 
the first order transverse derivative in \eqref{eq:6} can be expressed as $\partial_z\, \cdot =\partial_z y\, \partial_y\, \cdot =\varepsilon^{-1}\partial_y\, \cdot$, so that the differential operator $\mathcal{L}_z$ in \eqref{eq:Lfunc1} changes into
\begin{equation}\label{eq:Lfunc2}
    \mathcal{L} := -\varepsilon D(\partial_{xx} + \varepsilon^{-2}\partial_{yy}) + u\partial_x  + \sigma,
\end{equation}
where the specific choice made in \eqref{eq:Lfunc1} for the advective field ($\pmb{u}=[u(z), 0]^T$) together with the assumption $D=\mathcal{O}\left( 1 \right)$ are exploited. Thus,
problem \eqref{eq:6} can be rewritten as
\begin{equation} \label{eq:8}
\left\{
\begin{array}{ll}
\partial_tc + \mathcal{L} c = f \quad &\mbox{in } \Omega\times \left(0,T\right] \\[2mm] 
\varepsilon D{\nabla}_\nu c=0 \quad &\mbox{on }\Gamma^w\times \left(0,T\right] \\[2mm]
c=c_B  \quad &\mbox{on } \Gamma^i\times \left(0,T\right]  \\[2mm]
\varepsilon D{\nabla}_\nu c =0  \quad &\mbox{on }\Gamma^o\times \left(0,T\right]  \\[2mm]
c=0 \quad &\mbox{in }\Omega\times \{ t=0\}, 
\end{array}
\right.
\end{equation}
where $\nu$ denotes the unit outward normal vector to the boundary $\partial \Omega$ of $\Omega$.\\
Now, by assuming $c$ to be regular enough with respect to $\varepsilon$,
 we asymptotically expand the solution to problem \eqref{eq:8} as
\begin{equation} \label{eq:10}
c(x,y,t;\varepsilon) = c_0(x, y, t) + \varepsilon c_1(x,y,t) + \varepsilon^2 c_2(x,y,t) + \mathcal{O}(\varepsilon^{3}).
\end{equation}
A similar expansion could be adopted also with respect to the time variable. However, following~\cite{auriault1995taylor},
we neglect such an issue, being beyond the goal of the paper.\\
It is standard to prove that the first term in the expansion depends only on the macro-scale, namely $c_0=c_0(x, t)$ (see \eqref{eq:solc0}). 

Goal of the homogenisation process is to exploit expansion \eqref{eq:10} to commute problem \eqref{eq:8} into an averaged equation which models the effective concentration $c_e(x,t):=c_0(x,t)+\varepsilon \overline{c}_1(x,t)$, with
\begin{equation}\label{eq:avgvelo}
	    \overline{c}_1(x,t):=\frac{1}{L}\int_{-\frac{L}{2}}^{\frac{L}{2}}c_1(x,y,t)\dd{y} \mbox{, } 
\end{equation}
the $y$-average of $c_1(x,y,t)$ in \eqref{eq:10}.
The dependence of $c_e$ on the fast variable $y$ is removed, the micro-scale contribution being approximated by the averaged term $\overline{c}_1(x,t)$. We will show that the effective concentration
$c_e$ turns out to be the solution to a 1D problem which is solved along the $x$-direction (see equation \eqref{eq:24b}), being
associated with the macro-scale only. 

Now, by replacing expansion \eqref{eq:10} into the differential equation in \eqref{eq:8}, we have
\begin{equation}\label{eq:11}
\begin{array}{lll}
\partial_tc_0+\varepsilon\partial_tc_1&=& \varepsilon D\partial_{xx}c_0+\varepsilon^{-1}D\partial_{yy}c_0+D\partial_{yy}c_1+\\[3mm] 
&+&\varepsilon D\partial_{yy}c_2-u\partial_x\left(c_0+\varepsilon c_1\right)- \sigma(c_0+\varepsilon c_1)+f+\mathcal{O}(\varepsilon^2),
\end{array}
\end{equation}
while the boundary conditions become
\begin{equation}\label{eq:12}
\begin{alignedat}{2}
	& D\big(\partial_yc_0+\varepsilon\partial_yc_1+\varepsilon^2\partial_yc_2\big)=\mathcal{O}(\varepsilon^3) \qquad && \text{on } \Gamma^w ,\\
	& c_0 = c_B \qquad && \text{on } \Gamma^i ,\\
	& D\partial_x c_0=0 \qquad && \text{on } \Gamma^o,
	\end{alignedat}
\end{equation}
%
with $c_0=0$ at the initial time. In particular, as shown below, the boundary data on $\Gamma^i$ and $\Gamma^o$ are not instrumental to the homogenisation process, so that the expansion \eqref{eq:10} is only partially involved in \eqref{eq:12}$_{2,3}$. %

The homogenisation process works as follows. We collect the terms in \eqref{eq:11}-\eqref{eq:12}  according to the successive powers of $\varepsilon$ (notice that 
we pair the terms of order $i$ on $\Omega$ with the ones of order $i+1$ on the boundary to comply with the same order of magnitude). Thus, we obtain a differential problem for each order $i$ of $\varepsilon$, whose solution is denoted by $c_{i+1}$. We proceed order by order, by focusing on low orders in this section: 

\vspace*{0.3cm}

\begin{itemize}
\item {\underline{order $-1$}} ($\varepsilon^{-1}$ in $\Omega$; $\varepsilon^0$ on $\partial \Omega$):
	\begin{equation}\label{eq:varepsmeno1}
	\begin{aligned}
		D\partial_{yy}c_0 &= 0 \quad \mbox{in } \Omega, \\
		D\partial_yc_0 &= 0 \quad \mbox{on } \Gamma^w.\\
	\end{aligned}
	\end{equation}
Equation \eqref{eq:varepsmeno1} identifies the solution $c_0$ up to a constant, and it confirms that $c_0$ does not depend on $y$, namely
	\begin{equation}\label{eq:solc0}
		c_0\left(x,y,t\right) = c_0\left(x,t\right).
	\end{equation}

\vspace*{0.3cm}
\item {\underline{order $0$}} ($\varepsilon^{0}$ in $\Omega$; $\varepsilon^1$ on $\partial \Omega$):
	\begin{equation}\label{eq:14}
	\begin{aligned}
	\partial_tc_0-D\partial_{yy}c_1+u\partial_xc_0 + \sigma c_0 &=f \quad &&\mbox{in } \Omega, \\
		D\partial_{y}c_1&=0 \quad &&\mbox{on } \Gamma^w.\\
	\end{aligned}
	\end{equation}
	Integrating the differential equation divided by $L$ with respect to $y$ on the interval $[-\textstyle\frac{L}{2},\textstyle\frac{L}{2}]$, and using the boundary condition, we obtain the so-called \textit{leading order equation}
	\begin{equation}\label{eq:15}
 \partial_t c_0 + \overline{u} \partial_xc_0+ \sigma c_0=f \quad \mbox{in } (0,L)\times\left(0,T\right],
	\end{equation}
being $\sigma$ and $f$ constant, and where the average velocity $\overline{u}$ along the transverse section is defined according to \eqref{eq:avgvelo} (notice that $\overline{u}$ is a constant, the physical velocity field being independent of $x$ and $t$).
We remark that, due to \eqref{eq:solc0}, equation \eqref{eq:15} depends only on the macro-scale independent space variable $x$, and characterizes the leading behaviour of the concentration. As expected, the flow is driven exclusively by the advection, while the effect of the diffusion appears at higher orders, due to the assumption $D_\varepsilon=O(\varepsilon)$. \\
To obtain an expression for $c_1$, we subtract equations \eqref{eq:14} and \eqref{eq:15}, thus recovering the boundary value problem
\begin{equation}\label{eq:16}
	\begin{aligned}
	-D\partial_{yy} c_1 + \left(u-\overline{u}\right)\partial_x c_0 = -D\partial_{yy} c_1 +  u'\partial_x c_0 &= 0
	 \quad &&\mbox{in } \Omega, \\
		D\partial_{y}c_1&=0 \quad &&\mbox{on } \Gamma^w,\\
\end{aligned}
\end{equation}
whose solution is still defined up to a constant. 
We take the average along the $y$-direction while imposing the boundary data on $\Gamma^w$ such that
 \begin{equation}\label{eq:16avg}
    0= \overline{u'}\partial_x c_0.
\end{equation}
Then, we subtract \eqref{eq:16} and \eqref{eq:16avg} and we get
\begin{equation}\label{eq:16new}
	-D\partial_{yy} c_1 + \left(u'-\overline{u'}\right)\partial_x c_0 = 0
	 \quad \mbox{in } \Omega. 
\end{equation}
 
Integrating once with respect to $y$, we have
\begin{equation}\label{eq:neumannugen}
    D\partial_yc_1=\left( \Phi(y)-\overline{u'}y\right)\partial_x c_0 + \alpha, 
\end{equation}
with $\Phi(y)$ the primitive of $u'$
with respect to $y$,
and
$\alpha=\alpha(x,t)$ an integration constant.
%
%
Now, exploiting the 
boundary condition in \eqref{eq:16}, the constant $\alpha$ has to satisfy the relations
\begin{subequations} \label{eq:eqa}
\begin{alignat}{1}
\alpha&=-\left(\Phi\left(-\frac{L}{2}\right)+\overline{u'}\frac{L}{2}\right)\partial_xc_0,\label{eq:eqa1}\\[2mm]
    \alpha&=-\left(\Phi\left(\frac{L}{2}\right)-\overline{u'}\frac{L}{2}\right)\partial_xc_0.\label{eq:eqa2}
\end{alignat}
\end{subequations}
%
Notice that the two relations do coincide, being  $\overline{u'}=L^{-1}\big[{\Phi(\frac{L}{2})-\Phi(-\frac{L}{2})}\big]$.
In particular, plugging \eqref{eq:eqa1} into \eqref{eq:neumannugen}, we obtain
%
\begin{equation}\label{eq:neumannugen2}
    D\partial_yc_1=\left[ \Phi(y)-\Phi\left(-\frac{L}{2}\right)-\left(\frac{L}{2}+y\right)\overline{u'}\right]\partial_x c_0
\end{equation}
that, after integrating again along $y$, becomes
\begin{equation}\label{eq:1_3_4}
\hspace*{-0.1cm}c_1\left(x,y,t\right)=\frac{1}{D} \left[\int \Phi({y})\dd{{y}}-\Phi\left(-\frac{L}{2}\right)y-\left(\frac{L}{2}y+\frac{y^2}{2}\right)\overline{u'}\right]\partial_x c_0 +b_1,
\end{equation}
with $b_1=b_1(x,t)$ an integration constant.

Thus, in the first contribution on the right-hand side, the dependence of $c_1$ on $x$ and $y$ can be separated into the derivative, $\partial_xc_0\left(x,t\right)$, with respect to the macro-scale, which depends on $x$ only, and a term, that we denote by ${\tilde \chi}_1(y)$, 
depending on the fast variable $y$ and on the velocity profile, so that \eqref{eq:1_3_4} becomes
\begin{equation} \label{eq:1_3_6}
		c_1(x,y,t) = {\tilde \chi}_1(y)\partial_x c_0(x,t) +b_1(x,t).
\end{equation}
Some hints about the computation of the integration constant $b_1$ will be provided in the next section when dealing with higher-order expansions.

\vspace*{0.3cm}

\item {\underline{order $1$}} ($\varepsilon$ in $\Omega$; $\varepsilon^2$ on $\partial \Omega$):
	\begin{equation}\label{eq:19}
	\begin{aligned}
	\partial_tc_1-D\left(\partial_{xx}c_0+\partial_{yy}c_2\right)+u\partial_xc_1 + \sigma c_1&=0 \quad &&\mbox{in } \Omega, \\
	D\partial_yc_2&=0 \quad &&\mbox{on } \Gamma^w.\\
	\end{aligned}
\end{equation}
Integrating the differential equation in \eqref{eq:19} with respect to $y$ on the interval $[-\frac{L}{2}, \frac{L}{2}]$, and using the boundary condition, we have
\begin{equation}\label{eq:20}
\begin{array}{ll}
&\partial_t\displaystyle\int_{-\frac{L}{2}}^{\frac{L}{2}}c_1(x,y,t)\dd{y}-LD\partial_{xx}c_0(x,t)\\[5mm]
&+\displaystyle\int_{-\frac{L}{2}}^{\frac{L}{2}}\big(u\partial_x c_1(x,y,t) + \sigma c_1(x,y,t)\big)\dd{y} =0,
\end{array}
\end{equation}
%
which can be rewritten in the form
\begin{equation}\label{eq:22}
\partial_t{\overline c}_1-D\partial_{xx}c_0+\overline u\partial_x{\overline c}_1 + \sigma \overline{c}_1=-\frac{1}{L}\int_{-\frac{L}{2}}^{\frac{L}{2}}u'(y) \partial_x c_1(x,y,t)\dd{y},
\end{equation}
where we have applied the generic 
average-fluctuation splitting,
$f(x,y)=\overline{f}(x)+f'(y)$, to function $u$, with $\overline{f}=\overline{f}(x)$ the $y$-average according to \eqref{eq:avgvelo} and $f'=f'(y)$ the fluctuation contribution associated with the $y$-direction (notice that function $u$ depends on $y$ only, so that $\overline u$ is constant while $u'=u'(y)$).
\end{itemize}

\vspace*{.3cm}

The asymptotic expansions derived above are instrumental to formulate a transport equation for the effective concentration $c_e=c_0+\varepsilon{\overline c}_1$. To this end, we sum the leading order equation \eqref{eq:15} with equation \eqref{eq:22} multiplied by $\varepsilon$, to obtain
\begin{equation}\label{eq:23}
\partial_t\left(c_0+\varepsilon{\overline c}_1\right)-\varepsilon D\partial_{xx}c_0+\overline u\partial_x\left(c_0+{\varepsilon\overline{ c}}_1\right) + \sigma(c_0+\varepsilon \overline{c}_1)=f- \displaystyle \frac{\varepsilon}{L}\int_{-\frac{L}{2}}^{\frac{L}{2}}u'\partial_x c_1\dd{y}.
\end{equation}
Then, after 
exploiting relation \eqref{eq:1_3_6}, the average-fluctuation splitting for $u'$ in the integral on the right-hand side, and the fact that $\overline{u}$ is constant,
we can rewrite equation \eqref{eq:23} as
\begin{equation}\label{eq:24}
	\partial_t c_e = \partial_x \Big[-\overline{u}c_e + \varepsilon D \Big( 1+ \frac{\overline{u}\overline{{\tilde \chi}_1}-\overline{u {\tilde \chi}_1}}{D}\Big)\partial_{x} c_e \Big] -\sigma c_e + f,
\end{equation}
the terms of order strictly greater than one being neglected. 
Thus, we are led to solve 
the $1$D advection-diffusion-reaction problem associated with the macro-scale 
%
%
\begin{equation}\label{eq:24b}
\left\{
\begin{array}{ll}
	\partial_t c_e = \partial_x \big[-\overline{u}c_e + D_{\rm eff}\, \partial_{x} c_e \big] -\sigma c_e + f \quad  &\mbox{in }(0, L)\times\left(0,T\right],\\[2mm]
	 c_e = c_B  \quad&\mbox{on } \{x = 0\}\times\left(0,T\right], \\[2mm]
	 D_{\rm eff}\partial_x c_e =0  \quad &\mbox{on } \{x = L\}\times\left(0,T\right], \\[2mm]
  c_e = 0 \quad &  \mbox{in }(0, L)\times\{t = 0\},
	\end{array}
 \right.
\end{equation}
where the so-called \textit{Taylor dispersion coefficient},
$D_{\rm eff}:=\varepsilon D \big[ 1+ (\overline{u}\overline{{\tilde \chi}_1}-\overline{u {\tilde \chi}_1})D^{-1}\big]\geq \varepsilon D$, collects the model contribution along the transverse direction.

Equation \eqref{eq:24b} represents the outcome of the homogenisation process. The effective concentration can be adopted as a 1D surrogate to the concentration $c$ in \eqref{eq:8}. However, $c_e$ is expected to offer a poorer approximation when compared with the HiMod solution in \eqref{eq:hieq}, since the HiMod approximation preserves the dependence on the transverse direction whereas $c_e$ definitely loses this information through the average  \eqref{eq:avgvelo}.\\
For this reason, to settle the HiPhom$\varepsilon$ method, we will refer to an asymptotic expansion such as the one in \eqref{eq:1_3_6} rather than to the actual solution of the homogenisation process. Indeed, $c_1$ depends on both $x$ and $y$ according to the separation of variable paradigm driving the HiMod expansion \eqref{eq:hieq}, in contrast to the effective concentration which is a function of the $x$ variable only.



\section{High-order asymptotic expansions for axial flows}
\label{sec:hoasy}

In this section we extend the two-scale expansion in the previous section to an arbitrary order. To this aim, we follow what done in~\cite{bakhvalov_homogenisation_1989} where such an extension is limited to the context of purely diffusive problems. Also in the high-order case, the goal we pursue is to derive a separable structure for the solution to problem \eqref{eq:8}, with the aim of being compliant with a HiMod discretization. 

Starting from the results in Section~\ref{sec:tae}, it follows that the first two terms of the expansion in \cref{eq:10} can be expressed as
\begin{equation}
    \label{eq:3.13}
    \begin{aligned}
    c_0(x,y,t) &= \chi_{0}^*(y) c_0(x,t) = c_0(x,t), \\[1mm]
    c_1(x,y,t) &= \chi_1^*(y) \partial_x c_0(x,t),
\end{aligned}
\end{equation}
with $\chi_{0}^*=\chi_{0}^*(y)$ constant and here set, for simplicity, to $1$  and 
$$
 \chi_{1}^*(y)=\frac{1}{D} \Big[\int \Phi({y})\dd{{y}}-\Phi\Big(-\frac{L}{2}\Big)y-\Big(\frac{L}{2}y+\frac{y^2}{2}\Big)\overline{u'}  + \beta_1 \Big],
$$
having assumed $b_1(x,t) = \beta_1(x,t) \partial_x c_0(x,t) $ in \eqref{eq:1_3_4}. 
To guarantee the uniqueness, it is necessary to add an extra condition on $\chi_1^*$. For simplicity, we impose 
\begin{equation}\label{eq:cond_chi1}
    \chi_{1}^*(0)=  F(0) \quad \text{with }F (y):= \frac{1}{D} \int\left[\int u'({y})\dd{{y}}\right]\dd{y}
\end{equation}
so that we have $\beta_1 = 0$, being $\Phi(y)=\int u'(y)\dd{y}$.
Of course, other conditions (e.g., a null average hypothesis or the assignment of a value at a specific point of the domain) are allowed since not affecting the final result.
In particular, moving from \eqref{eq:varepsmeno1} and \eqref{eq:16}, functions $\chi_0^*$ and $\chi_1^*$ are the solutions to problems
\begin{equation}\label{eq:homo2}
    \begin{array}{ll}
    D \partial_{yy} \chi_0^*(y) = 0 \quad&\mbox{in } \left(-\displaystyle\frac{L}{2}, \frac{L}{2}\right),\ \  \mbox{with} \ \  D\partial_y \chi_0^*\left( \pm \displaystyle\frac{L}{2}\right) = 0, \\[4mm]
    D \partial_{yy} \chi_1^*(y) = u'(y) \quad&\mbox{in } \left(-\displaystyle\frac{L}{2}, \frac{L}{2}\right),\ \  \mbox{with} \  \ D\partial_y \chi_1^* \left( \pm \displaystyle \frac{L}{2}\right)= 0,
    \end{array}
\end{equation}
where $\chi_0^*$ has been chosen equal to $1$ among all the possible constant solutions, while  $\chi_1^*$ additionally satisfies condition \eqref{eq:cond_chi1}. 

\smallskip

The current goal is to derive 
an explicit expression for $c_2(x,y,t)$ 
analogous to the ones in \eqref{eq:3.13} for $c_0$ and $c_1$, starting from the differential problem in \cref{eq:19}. 
To this aim, we adopt the following procedure:
\begin{enumerate}
    \item We employ 
      the definitions in \cref{eq:3.13} for the previous orders, $c_0$ and $c_1$, in problem \cref{eq:19}, so that
   we have
\begin{equation}\label{eq:19_bis}
	\begin{aligned}
	\chi_1^* \partial_{xt} c_0 -D\left(\partial_{xx}c_0+\partial_{yy}c_2\right)+u\chi_1^*\partial_{xx}c_0 + \sigma \chi_1^* \partial_{x}c_0&=0 \quad &&\mbox{in } \Omega, \\
	D\partial_yc_2&=0 \quad &&\mbox{on } \Gamma^w.\\
	\end{aligned}
\end{equation}
\item We exploit the leading order equation \cref{eq:15} to replace the time derivative of $c_0$ in \eqref{eq:19_bis}, to obtain
\begin{equation} \label{eq:c2_homo_1}
    D\partial_{yy} c_2 = \chi_1^* \partial_{x} (f - \overline{u} \partial_x c_0 - \sigma c_0) -D\partial_{xx}c_0 + u\chi_1^*\partial_{xx}c_0 + \sigma \chi_1^* \partial_{x}c_0,
\end{equation}
namely, 
\begin{equation} \label{eq:c2_homo}
    D\partial_{yy} c_2 = \left( u'\chi_1^* - D \right) \partial_{xx}c_0,
\end{equation}
being  $\partial_{x} f = 0$ since $f$ is assumed constant in \eqref{eq:6}, and after applying the average-fluctuation splitting to $u$.
\item We subtract from \cref{eq:c2_homo} the associated average along the $y$-direction while imposing the boundary data on $\Gamma^w$, to get
\begin{equation} \label{eq:3.15}
    D \partial_{yy} c_2 = \left(u'\chi_1^* - \overline{u'\chi_1^*} \right) \partial_{xx}c_0 =\left(u'\chi_1^*\right)'\partial_{xx}c_0 ,
\end{equation}
the average-fluctuation splitting being now applied to $u'\chi_1^*$.
\item By mimicking the computations leading to \eqref{eq:1_3_4}, we can express $c_2$ as
\begin{equation}\label{eq:homo7}
    c_2 (x, y, t) = \chi_2^*(y)\partial_{xx}c_0(x, t) ,
\end{equation}
with $\chi_2^*$ solution, up to a constant $b_2=b_2(x, t)$, to problem
\begin{equation}\label{eq:homo6}
     D\partial_{yy} \chi_2^* = \left(u'\chi_1^*\right)' \quad \mbox{in } \bigg(-\frac{L}{2}, \frac{L}{2}\bigg),\ \  \mbox{with} \  D\partial_y \chi_2^* \Big( \pm \displaystyle\frac{L}{2}\Big)= 0,
\end{equation}
where constant $b_2$ is fixed by the imposition of an extra condition analogous to the one assigned to $\chi_1^*$ in \eqref{eq:cond_chi1}, after assuming $b_2(x,t)=\beta_2(x, t)\partial_{xx}c_0(x, t)$ 
 (see Proposition~\ref{lem:theo1}).
\end{enumerate}

\smallskip

As a consequence, the asymptotic expansion \cref{eq:10} truncated at the second order can be reformulated as
\begin{equation}\label{eq:soasy}
\begin{array}{rcl}
 c(x, y, t; \varepsilon)&=&\displaystyle \sum_{j=0}^{2}\varepsilon^j c_j(x, y, t)+\mathcal{O}(\varepsilon^3)\\[5mm]
 &=&  \chi^*_0(y) \, c_0(x, t)+  \varepsilon\, \chi^*_1(y)\, \partial_x c_0(x, t)  + \varepsilon^2 \,\chi_2^*(y)\,\partial_{xx}c_0(x, t) + \mathcal{O}(\varepsilon^3)\\[2mm]
 &=& \displaystyle \sum_{j=0}^{2}\varepsilon^j 
 \, \chi^*_j(y)\, \partial_{x^j} c_0(x, t)+\mathcal{O}(\varepsilon^3),
\end{array} 
\end{equation}
where $\partial_{x^j}$ denotes the derivative of order $j$ with respect to the independent variable $x$, with $\partial_{x^0}$ the identity operator. This means that the concentration $c$ does coincide with a linear combination of terms depending on $(x, t)$ with coefficients depending on the $y$ variable and on the homogenisation parameter $\varepsilon$. \\
Relation \eqref{eq:soasy} 
can be generalized to 
a higher order, as stated in the result below:

\begin{proposition}
\label{lem:theo1}
For any $\varepsilon$, the solution $c = c(x,y,t;\varepsilon)$ to the homogenised problem \cref{eq:8} does coincide with the asymptotic expansion
\begin{equation}\label{prima_espans}
c(x,y,t;\varepsilon) := \displaystyle \sum_{i=0}^{+\infty}\varepsilon^i c_i(x,y,t),
\end{equation}
with $c_0$ the solution to the leading order equation \eqref{eq:15} and with 
\begin{equation}\label{eq:lemma11}
    c_{i}(x,y,t) :=  
    \chi_{i}^{*}(y)\, \partial_{x^{i}}c_0(x,t) \quad  {i> 0},
\end{equation}
where
function $\chi_i^*$ coincides, up to a constant, with the solution to the boundary value problem
\begin{equation}\label{eq:chi_i_plus_cond_a}
 D \partial_{yy} \chi_{i}^* = \left(u'\chi_{i-1}^*-D\chi_{i-2}^*\right)' 
 \quad \mbox{in } \left(-\frac{L}{2}, \frac{L}{2}\right),\ \  \mbox{with} \  \ D\partial_y \chi_i^* \left( \pm \frac{L}{2}\right)= 0
\end{equation} 
being $\chi^*_{-2}=\chi^*_{-1}=0$, $\chi_0^*=1$. In particular, to guarantee the uniqueness of the solution to problem \eqref{eq:chi_i_plus_cond_a}, we pick $\chi_i^*$ such that
\begin{equation}\label{eq:chi_i_plus_cond_b}
      \chi_i^*(0) = F_i(0),
\end{equation}
with  
\begin{equation}\label{Fdef}
F_i(y):=\frac{1}{D} \int\left[\int \phi_i({y})\dd{{y}}\right]\dd{y}
\end{equation}
and $\phi_i(y) = u'(y)\chi_{i-1}^*(y)-D\chi_{i-2}^*(y)$.
\end{proposition}

\begin{proof}\label{pr:proof1}
We proceed by induction. Statements 
\cref{eq:lemma11}-\eqref{eq:chi_i_plus_cond_b} are true for $i=0$ (actually, in Section~\ref{sec:tae} and in \eqref{eq:homo7}-\eqref{eq:homo6}, we have explicitly derived these relations also for $i=1$ and $i=2$, respectively where condition \eqref{eq:chi_i_plus_cond_b} is used to get rid of the definition up to a constant for $\chi_1^*$ and $\chi^*_2$).

Now, we assume that relations 
\cref{eq:lemma11}-\eqref{eq:chi_i_plus_cond_b} are true for a generic index $i$ and we prove they still hold for the next index, $i+1$.\\
As a first goal, we derive the differential problem satisfied by $c_{i+1}$, by generalizing computations done to derive the separated representation \eqref{eq:1_3_6} for $c_1$.
In \eqref{eq:11}, we collect terms 
associated with $\varepsilon^i$ in $\Omega$ and with $\varepsilon^{i+1}$
on $\partial \Omega$, so that we have
\begin{equation}\label{eq:homo8}
    \partial_t c_i = D\partial_{xx}c_{i-1}+D\partial_{yy}c_{i+1}-u\partial_x c_i - \sigma c_i \quad \mbox{in } \Omega,
\end{equation}
completed by the boundary condition $D\partial_y c_{i+1}=0$ on $\Gamma^w$. Now, we rewrite terms depending on both $c_{i-1}$ and $c_i$ by exploiting \eqref{eq:lemma11} and by resorting to the leading order equation \eqref{eq:15} in order to get rid of the time derivative on the left-hand side.
This yields
\begin{equation}\label{eq:terms_ci}
\begin{array}{rcl} 
 \partial_t c_i &=& -\chi_i^*\overline{u}\partial_{x^{i+1}}c_0-\chi_i^*\sigma\partial_{x^{i}} c_0,\\[3mm]
    D\partial_{xx}c_{i-1} &=& D \chi_{i-1}^{*}\partial_{x^{i+1}}c_0,\\[3mm]
    u\partial_x c_i &=& u\chi_{i}^{*} \partial_{x^{i+1}}c_0,\\[3mm]
    \sigma c_i &=& \sigma\chi_{i}^{*} \partial_{x^{i}}c_0.
\end{array}
\end{equation}
By replacing \eqref{eq:terms_ci} in \eqref{eq:homo8}, we obtain
\begin{equation}\label{eq:proof_thm1_num0}
     D\partial_{yy}c_{i+1}=\left(u'\chi_i^*-D\chi_{i-1}^*\right)\partial_{x^{i+1}}c_0.
\end{equation}
Consistently with the computations leading to \cref{eq:3.15}, we take the average of problem \eqref{eq:proof_thm1_num0} along the $y$-direction, while imposing the
boundary data for $c_{i+1}$ on $\Gamma^w$, which yields
\begin{equation}\label{eq:average_proof}
     0=\overline{\left(u'\chi_i^*-D\chi_{i-1}^*\right)}\partial_{x^{i+1}}c_0.
\end{equation}
Then, we subtract  \cref{eq:proof_thm1_num0} and \cref{eq:average_proof} and we apply the average-fluctuation splitting  to ${u'\chi_{i}^* -D\chi_{i-1}^*}$. This leads to the differential problem
\begin{equation}\label{eq:3_24}
    D\partial_{yy}c_{i+1} =
    \left(u' \chi_i^*-D\chi_{i-1}^*\right)'\partial_{x^{i+1}}c_0.
\end{equation}
Now, on the right-hand side, we adopt the notation $\phi_{i+1}=\left(u'\chi_i^*-D\chi_{i-1}^*\right)$, thus obtaining 
\begin{equation}\label{eq:cpi_yy_phi}
    D\partial_{yy}c_{i+1}  
    =(\phi_{i+1} - \overline{\phi_{i+1}})\partial_{x^{i+1}}c_0.
\end{equation}
After a first integration
with respect to $y$, we have
\begin{equation}\label{eq:cpi_y_phi}
    D\partial_{y}c_{i+1}   
    =(\Phi_{i+1} - \overline{\phi_{i+1}}y)\partial_{x^{i+1}}c_0 + \alpha_{i+1}(x,t),
\end{equation}
with $\Phi_{i+1}=\Phi_{i+1}(y)=\int \phi_{i+1}(y) \dd{y}$ the primitive of $\phi_{i+1}$ with respect to $y$, 
and $\alpha_{i+1}=\alpha_{i+1}(x,t)$ the integration constant. To compute 
$\alpha_{i+1}$, we exploit the boundary conditions 
satisfied by $\chi_{j}^*$, for $1\le j\le i$, at $y=\pm \frac{L}{2}$ (i.e., $D \partial_{y} \chi_{j}^*(\pm \frac{L}{2})=0$), 
and by $c_{i+1}$ on $\Gamma^w$ (i.e., $D\partial_y c_{i+1}= 0$ on $\Gamma^w$). It turns out that 
\begin{equation}\label{eq:bc_alpha}
     \alpha_{i+1}(x,t) =  - \left( \Phi_{i+1}\left(\pm \frac{L}{2}\right) - \left(\pm \frac{L}{2}\right)\overline{\phi_{i+1}} \right)\partial_{x^{i+1}}c_0.
\end{equation}
The two expressions for $\alpha_{i+1}$ do coincide, thanks to the definition of $\overline{\phi_{i+1}}=\big[\phi_{i+1}(\frac{L}{2})-\phi_{i+1}(-\frac{L}{2})\big]L^{-1}$. Thus, 
by substituting \cref{eq:bc_alpha} for the choice $y=-\frac{L}{2}$ into \cref{eq:cpi_y_phi}, it follows
\begin{equation}\label{eq:cpi_y_al}
    D\partial_{y}c_{i+1} 
    =\left[\Phi_{i+1}(y)- \Phi_{i+1}\left(- \frac{L}{2}\right) -\left( \frac{L}{2}+ y \right) \overline{\phi_{i+1}}  \right]\partial_{x^{i+1}}c_0. 
\end{equation}
Now, by integrating again with respect to $y$, we obtain 
\begin{equation}\label{ci+1_def}
     c_{i+1}  = \frac{1}{D}\bigg[\displaystyle \int \Phi_{i+1}(y) \dd{y} - \Phi_{i+1}\Big(-\frac{L}{2}\Big) y -  \Big( \frac{L}{2} y + \frac{y^2}{2} \Big) \overline{\phi_{i+1}} \bigg] \partial_{x^{i+1}}c_0 + b_{i+1},
\end{equation}
with $b_{i+1}=b_{i+1}(x,t)$ an integration constant. In particular, 
we assume that $b_{i+1}$ takes the form $b_{i+1}(x,t) = \beta_{i+1}(x,t) \partial_{x^{i+1}}c_0(x,t)$, so that
\eqref{ci+1_def} can be rewritten 
according to \eqref{eq:lemma11}, being
\begin{equation}\label{separa-ci+1}
     c_{i+1}(x, y, t)  ={\chi}_{i+1}^*(y) \, \partial_{x^{i+1}}c_0(x, t)
\end{equation}
with
\begin{equation}\label{practiceCHI}
{\chi}_{i+1}^* =  \frac{1}{D}\left[\displaystyle \int \Phi_{i+1}(y) \dd{y} - \Phi_{i+1}\left(-\frac{L}{2}\right) y -  \left( \frac{L}{2} y + \frac{y^2}{2} \right) \overline{\phi_{i+1}} + \beta_{i+1} \right]. 
\end{equation}
To get rid of constant $\beta_{i+1}$, we impose the extra condition 
\begin{equation}\label{chirel}
{\chi}_{i+1}^*(0) =  \displaystyle F_{i+1}(0),
\end{equation}
with $F_{i+1}$ defined according to  \eqref{Fdef},
which leads to pick $ \beta_{i+1} = 0$. By replacing the separated representation \eqref{separa-ci+1} for $c_{i+1}$ into the differential problem \eqref{eq:3_24} and into the boundary condition $D\partial_y c_{i+1}$ on $\Gamma^w$, it turns out that function $\chi^*_{i+1}$ solves the problem 
\begin{equation}\label{eq:chi_iplus1_plus_cond}
     D \partial_{yy} \chi_{i+1}^* = \big(u'\chi_{i}^*-D\chi_{i-1}^*\big)' \quad \mbox{with } D \partial_y \chi^*_{i+1}\Big( \pm \displaystyle\frac{L}{2}\Big) = 0,
\end{equation}
jointly with condition \eqref{chirel}, namely $\chi^*_{i+1}$
satisfies \eqref{eq:chi_i_plus_cond_a}-\eqref{eq:chi_i_plus_cond_b}. This concludes the proof.
\end{proof}
\vspace{2mm}

\begin{corollary}[separable representation in an axial flow regime]\label{proposition1}
For any $\varepsilon$, the solution $c = c(x,y,t;\varepsilon)$ to the homogenised problem \cref{eq:8} can be represented in a separated way as
\begin{equation}\label{eq:expansion_cm}
    c(x,y,t;\varepsilon) 
    =\sum_{i=0}^{+\infty}\chi_i^*(y)c_i^*(x,t; \varepsilon),
\end{equation}
with 
\begin{equation}\label{eq:cm}
    {c}^*_i(x,t; \varepsilon) = \varepsilon^i\partial_{x^i}c_0(x,t).
\end{equation}
In particular, we can consider the truncated of order $k$ of the asymptotic expansion in \eqref{eq:expansion_cm}, given by 
\begin{equation}\label{exp1}
{\mathcal C}_k(x,y,t;\varepsilon) =\sum_{i=0}^{k}\chi_i^*(y){ c}_i^*(x,t; \varepsilon)+\mathcal{O}\big(\varepsilon^{k+1}\big).
\end{equation}
\end{corollary}

\begin{proof}
The asymptotic expansion \eqref{eq:expansion_cm}-\eqref{eq:cm} follows simply by exploiting the explicit definition \eqref{eq:lemma11} of coefficients $c_i$ in \eqref{prima_espans}.
\end{proof}




\section{The HiPhom\texorpdfstring{$\pmb{\varepsilon}$}-   method}
\label{sec:hip}




In the two previous sections, we have shown that both a HiMod discretisation and a high-order asymptotic expansion lead to a representation for the solution to problem \eqref{eq:6} characterized by a separation of variables, which distinguishes between an axial and a transverse contribution. This is evident by comparing 
equation \eqref{eq:hieq} with the expansion in
\eqref{exp1}.\\
The main feature which discriminates the two modelling approaches consists in the selection of the functions associated with the transverse dynamics. 
A HiMod discretisation selects a-priori a set of modal basis functions (e.g., sinusoidal functions, Legendre polynomials or, more in general, the solutions to the auxiliary Sturm-Liouville problem \eqref{eq:SLE1}) to model the transverse behaviour. On the contrary, a high-order asymptotic expansion explicitly derives the set of functions $\{\chi^*_j\}$, to approximate the lateral fluctuations with respect to the longitudinal dynamics. 

Starting from these considerations, we propose a new model reduction procedure that we name HiPhom$\varepsilon$ (aka, HIgh order Projection-based HOMogenisation), which combines HiMod reduction with a high-order asymptotic expansion. In more detail, we adopt functions $\chi^*_i$ 
defined in Proposition~\ref{lem:theo1} as modal basis functions to perform a HiMod discretisation. The rationale behind this is to involve in the HiMod approximation a modal basis of functions which are strictly related to the problem at hand. Indeed, by construction, functions $\chi^*_i$  depend on 
the transverse distribution of the flow field as well as on the diffusion coefficient.\\
To replace modes 
$\varphi_k$ in \eqref{eq:hieq} with the new basis of functions, we have to ensure that functions $\chi^*_i$ are orthonormal with respect to the $L^2(\widehat{\gamma})$-scalar product. To this aim, we resort to the Gram-Schmidt algorithm.
For consistency with the reference domain $\widehat{\Omega}$ employed in a HiMod reduction, we map the vertical dimension $[-L/2, L/2]$ of the 
homogenised domain $\Omega$ into the HiMod reference fibre $\widehat \gamma=[0, 1]$ via the map
$\Theta: [-L/2, L/2] \to [0, 1]$ such that $\Theta (y)=\widehat{z}=L^{-1}\left(y+L/2\right)$, with $\widehat{z}$ and $y$ the independent variable associated with the transverse direction in $\widehat \Omega$ and $\Omega$, respectively. Thus, the new basis of modal functions orthonormal on $\widehat{\gamma}$ turns out to be defined by
\begin{equation}\label{eq:3_27}
    \chi_0(y(\widehat{z})):=1, \quad  \chi_i(y(\widehat{z})) := \frac{1}{a_i}\left(\chi^*_i(y(\widehat{z}))-\displaystyle \sum_{j=0}^{i-1}p_{i,j}\chi_j(y(\widehat{z}))\right),
\end{equation}
with $y(\widehat{z})=L \widehat{z} - L/2$, 
\begin{equation}
\label{eq:3_26}
    p_{i,j} := \displaystyle \int_{\widehat \gamma} \chi^*_i(y(\widehat{z}))\, \chi_j(y(\widehat{z}))\, d \widehat{z}, \quad
    a_i := \Big \|\chi^*_i(y(\widehat{z}))-\sum_{j=0}^{i-1}p_{i,j}\chi_j(y(\widehat{z}))\Big\|_{L^2{(\widehat \gamma)}}.
\end{equation}
Accordingly, we can introduce the HiPhom$\varepsilon$ reduced space
\begin{equation} \label{eq:vmspacenew}
    V_{m, \varepsilon} =\bigg\{v_{m, \varepsilon}(x,z)=\sum_{k=0}^{m-1}\widetilde{v}_{k, \varepsilon}(x)\chi_k\left(y(\psi_x(z))\right) :   \widetilde{v}_{k, \varepsilon}\in V_{1D},  x\in\Omega^{1D}_{\varepsilon},  z\in\gamma_x\bigg\},
\end{equation}
so that the HiMod expansion in \eqref{eq:hieq} is replaced by the HiPhom$\varepsilon$ approximation
\begin{equation}\label{HP_exp}
c_{m, \varepsilon}(t)(x,z) =\sum_{k=0}^{m-1}\widetilde{c}_{k, \varepsilon}(t)(x)\chi_k\left(y(\psi_x(z))\right)=\sum_{k=0}^{m-1}\sum_{s=1}^{N_h}\widetilde{c}_{ks, \varepsilon}(t)\theta_s(x)\chi_k\left(y(\psi_x(z))\right).
\end{equation}
At each time $t^n$ of the temporal discretisation,
coefficients $\{ \widetilde{c}_{ks, \varepsilon}(t^n)\}$ represent the actual unknowns of the HiPhom$\varepsilon$ discretisation. According to an offline paradigm,  modal functions $\chi_k$ can be precomputed once and for all, since independent of time. Then, the associated HiPhom$\varepsilon$ stiffness matrix is assembled so that,
at each $t^n$, the reduced system of order $m N_h$ defined in Section~\ref{sec:himod} is solved to provide coefficients $\widetilde{c}_{ks, \varepsilon}(t^n)$.

\begin{remark}
A cross-comparison between \eqref{HP_exp} and \eqref{exp1} 
hints an explicit relation between coefficients $\widetilde{c}_{k, \varepsilon}$ and $c^*_i$. However, HiPhom$\varepsilon$ just computes functions 
$\chi_k$ to yield the modal basis for the HiMod approximation,
thus skipping the explicit computation of $c_0$ and of the associated derivatives.
\end{remark}

\begin{remark}
In practice, to build the HiPhom$\varepsilon$ basis of functions $\{ \chi_i\}$, we compute functions $\chi^*_i$ directly by exploiting formula \eqref{practiceCHI}, instead of solving the differential problem \eqref{eq:chi_i_plus_cond_a}-\eqref{eq:chi_i_plus_cond_b}, and then we proceed by applying the Gram-Schmidt orthonormalisation in \eqref{eq:3_27}.
\end{remark}

Because of the problem-driven derivation,
we do expect that the modal basis, $\{\chi_k\}$, supporting a HiPhom$\varepsilon$ approximation  
turns out to be more informative when compared with a completely general set $\{ \varphi_k\}$ of functions.
The HiMod setting which is closest to a HiPhom$\varepsilon$ discretisation is the educated basis approach in \eqref{eq:SLE1}. 
In the next section, we will compare the performance of these two methods in order to identify which is the best choice. 



\section{Numerical assessment}\label{num_Sec}
In this section, we apply the HiPhom$\varepsilon$ model reduction to two benchmark case studies, 
characterized by a different velocity profile, by considering both steady and unsteady configurations.
In particular, we will investigate the reliability of the HiPhom$\varepsilon$ reduction, both from a qualitative and a quantitative viewpoint. 
\subsection{The Poiseuille flow} 
\label{sec:experiments}
The first configuration deals with 
the Poiseuille flow, namely 
a standard parabolic profile that models a steady, axisymmetric, laminar flow induced by a (constant) pressure gradient along the axial direction of an ideally infinite pipe.\\
In practice, we solve problem \eqref{eq:6} in 
a finite computational domain $\Omega_{\varepsilon} = (0,2)\times (-\varepsilon/2,\varepsilon/2)$ for $\varepsilon=0.2$, with $\Gamma_{\varepsilon}^i = \{x=0\}\times (-\varepsilon/2,\varepsilon/2)$, $\Gamma_{\varepsilon}^o = \{x=2\}\times (-\varepsilon/2,\varepsilon/2)$, and  $\Gamma_{\varepsilon}^w = [0,2]\times \{\varepsilon/2,\varepsilon/2\}$, and we assume to be in the presence of a stationary regime (i.e., $c=c(x, z)$).
We set the velocity field as
\begin{equation}\label{eq:poiseuille_case1}
    \mathbf{u}=[u(z), 0]^T =  \left[2\bar{v}\Big(1-\Big(\frac{z}{\varepsilon}\Big)^2\Big), 0\right]^T,
\end{equation}
with $\bar{v}=10$; 
$\sigma=1$; $f = 0$; $c_B = 1$. Moreover, we set $D=1$
consistently with the assumption $D=\mathcal{O}(1)$ in Section~\ref{sec:tae}, the diffusion coefficient being thus equal to $D_\varepsilon = \varepsilon D=0.2$.

\begin{figure}[!htb]
\hspace*{-.8cm}
\includegraphics[width=0.55\textwidth]{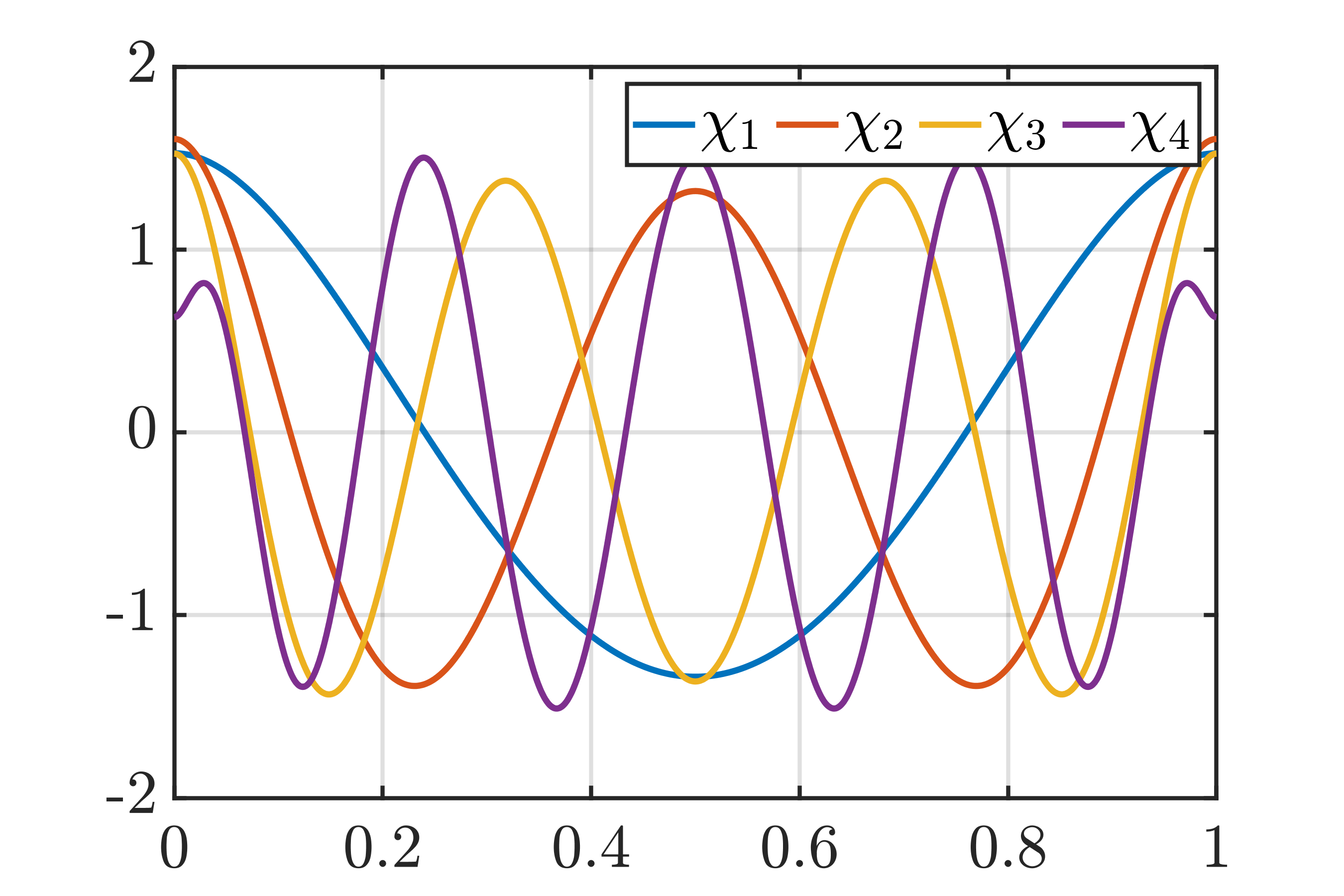}
\hspace*{-.8cm}
\includegraphics[width=0.55\textwidth]{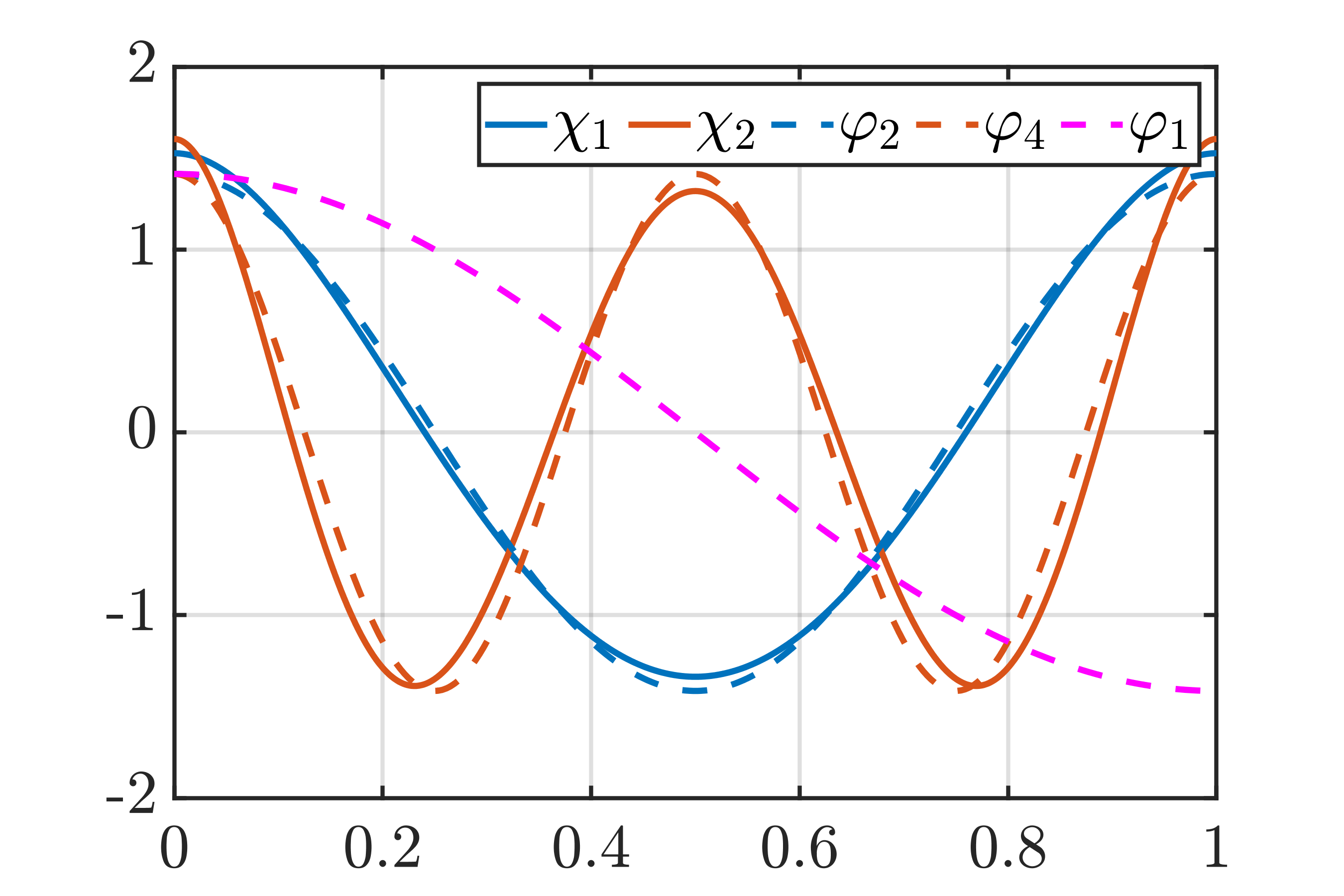}
\caption{Poiseuille flow test case: 
HiPhom$\varepsilon$ modal basis functions $\chi_i$, for $i=1, \ldots, 4$ (left); comparison 
among HiPhom$\varepsilon$ and HiMod modal basis functions (right).}\label{basis_fs}
\end{figure}
\begin{figure}[!htbp]
    \centering
    \hspace*{-.3cm}
    \includegraphics[height=1.55cm]{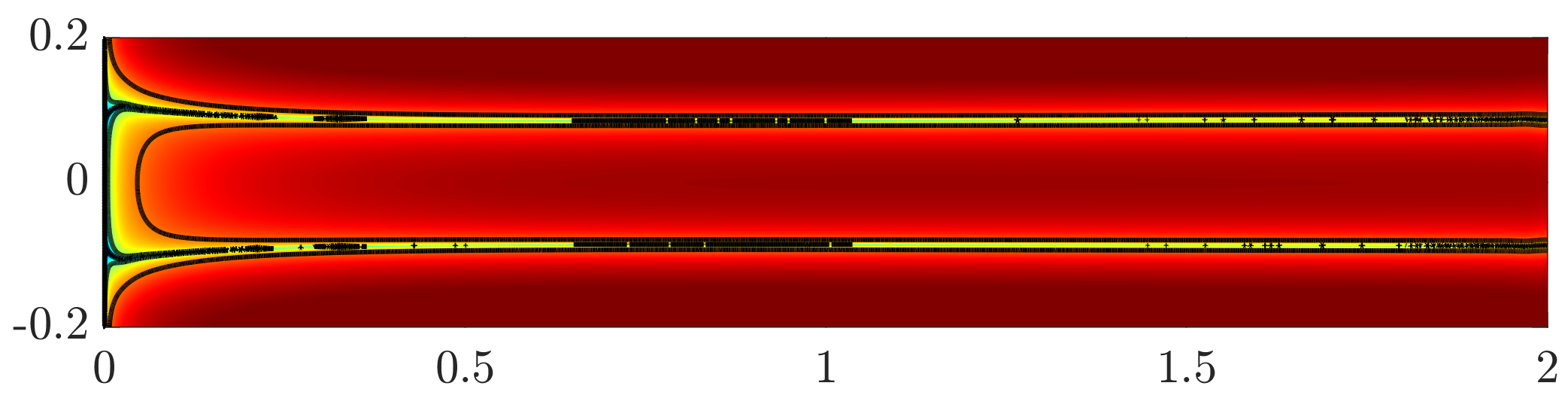}
    \hspace*{-.1cm}
    \includegraphics[height=1.55cm]{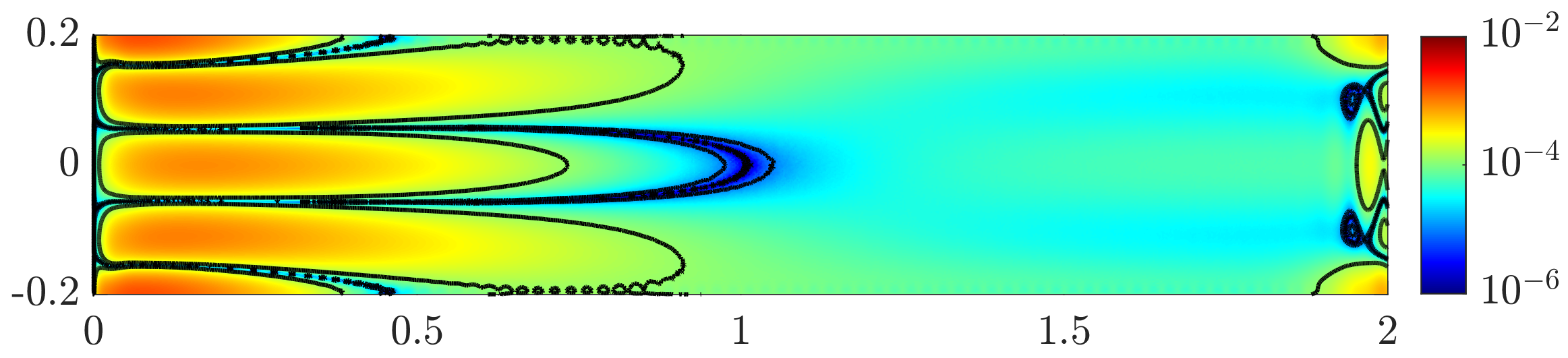}
    \hspace*{-.3cm}
    \includegraphics[height=1.55cm]{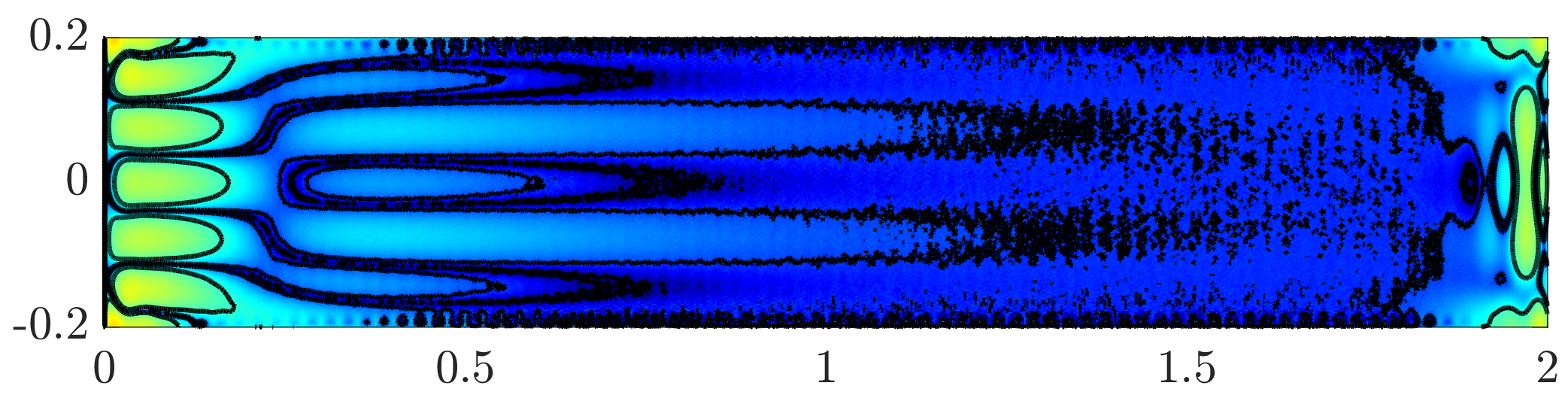}
    \hspace*{-.1cm}
    \includegraphics[height=1.55cm]{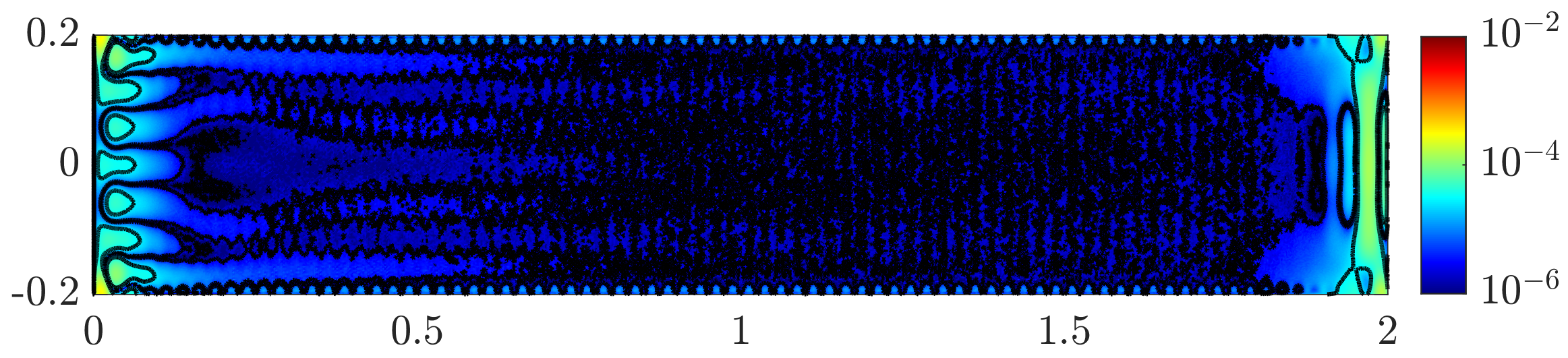}
     \caption{Poiseuille flow test case: spatial distribution of the absolute modelling error associated with the HiPhom$\varepsilon$ approximation for $m= 1$, $2$ (top, left and right) and $m=3$, $4$ (bottom, left and right).}
     \label{fig:errpoinewnew}
   \end{figure}
\begin{figure}[!htbp]
    \centering
    \hspace*{-.3cm}
    \includegraphics[height=1.55cm]{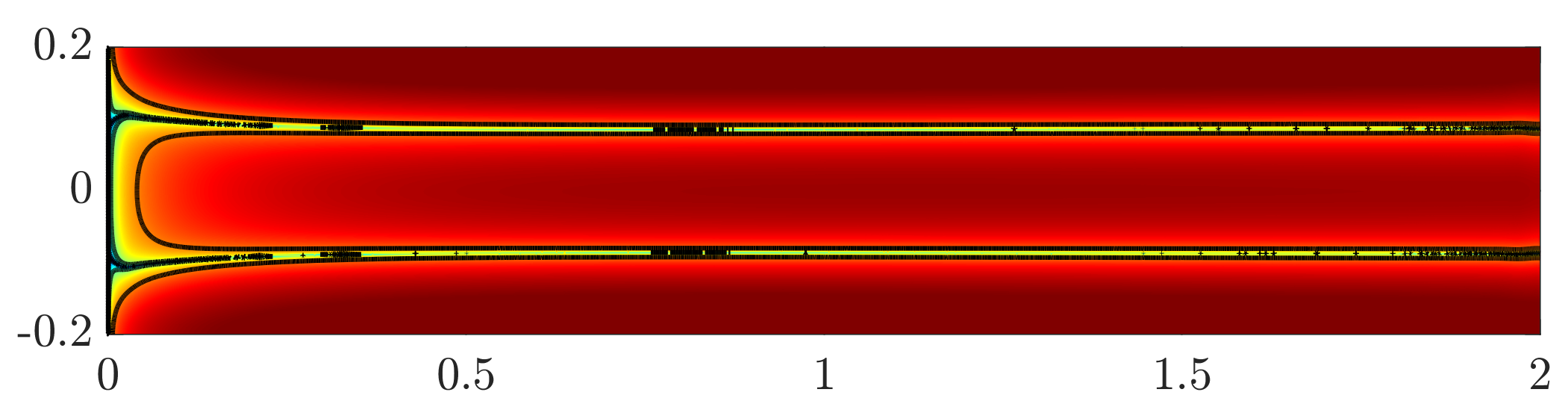}
    \hspace*{-.1cm}
    \includegraphics[height=1.55cm]{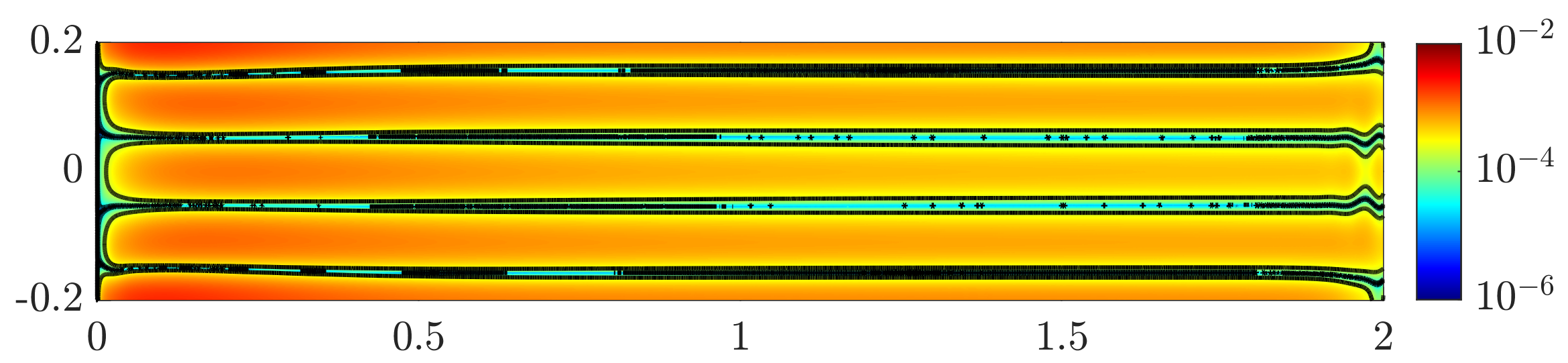}
    \hspace*{-.3cm}
    \includegraphics[height=1.55cm]{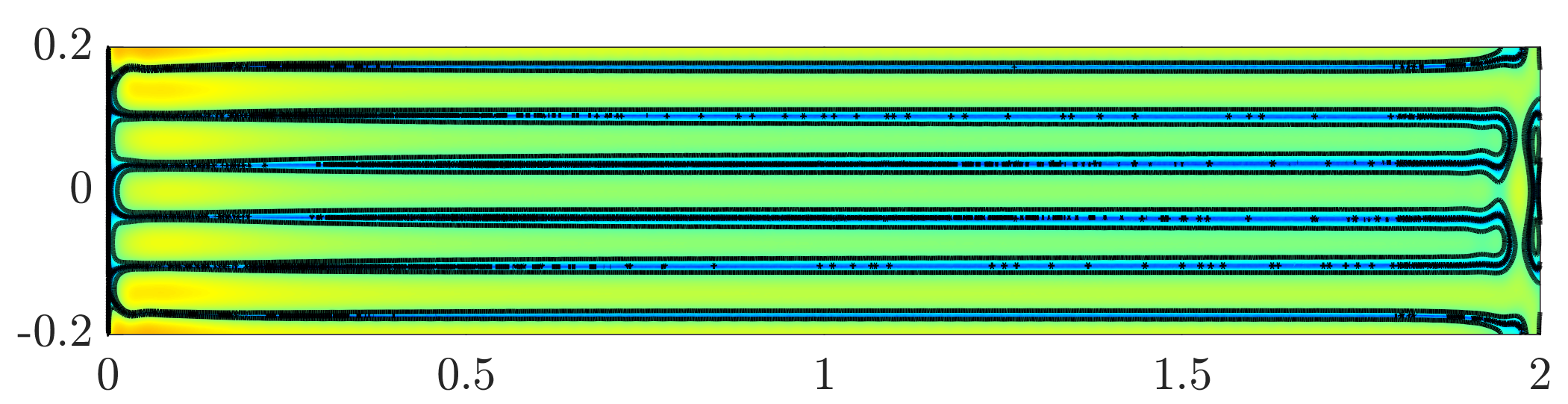}
    \hspace*{-.1cm}
    \includegraphics[height=1.55cm]{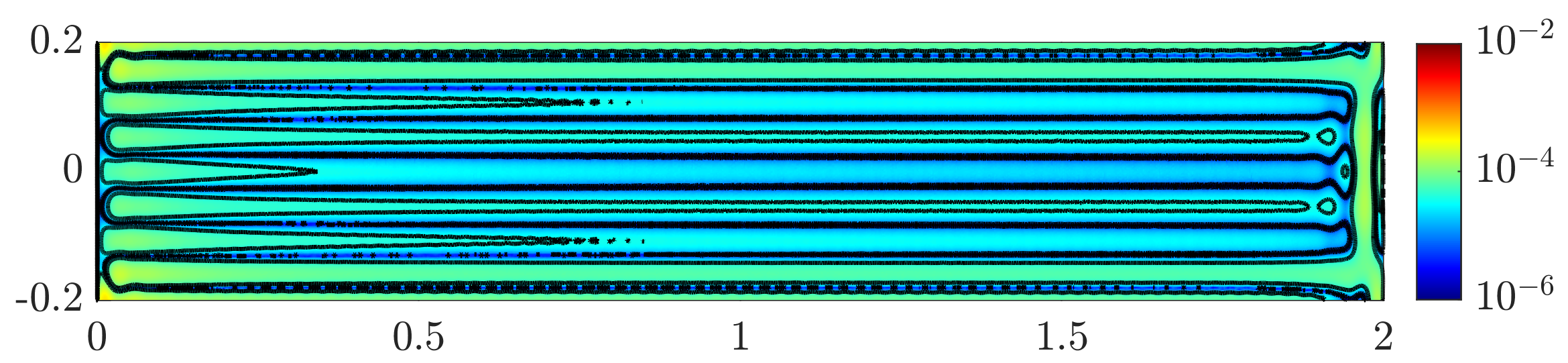}
     \caption{Poiseuille flow test case: spatial distribution of the absolute modelling error associated with the HiMod approximation for $m= 1$, $3$ (top, left and right) and $m=5$, $7$ (bottom, left and right).}
     \label{fig:errhimod}
   \end{figure}

In the left panel of Figure~\ref{basis_fs}, we show the four HiPhom$\varepsilon$ basis functions $\chi_i$, for $i=1, \ldots, 4$. We observe that HiPhom$\varepsilon$ employs symmetric functions only, due to the symmetry of the velocity profile with respect to the $x$-axis.\\
In the right panel of the same figure, we compare the HiPhom$\varepsilon$ modes $\chi_1$ and $\chi_2$ with the 
HiMod educated modal basis functions 
$\varphi_i$, for $i=1, 2, 4$, 
which are derived by setting ${\mathcal L}_S=D \partial_{yy}$ in \eqref{eq:SLE1}. It is evident 
the considerable similarity between $\chi_1$ and $\varphi_2$ as well as the good matching between $\chi_2$ and $\varphi_4$.  
The HiMod modal basis 
includes also non-symmetric functions  (see, e.g., mode $\varphi_1$ in the figure), which 
turn out to be redundant. This justifies the discrepancy in the numbering of the matching HiPhom$\varepsilon$ and HiMod basis functions. 
The expectation is that HiMod involves a higher number of modal functions with respect to the HiPhom$\varepsilon$ approach in order to guarantee a certain accuracy to the reduced solution (compare Figures~\ref{fig:errpoinewnew} and~\ref{fig:errhimod}). 
%

\begin{figure}[!htbp]
    \centering
    \hspace*{-.3cm}
    \includegraphics[height=1.55cm]{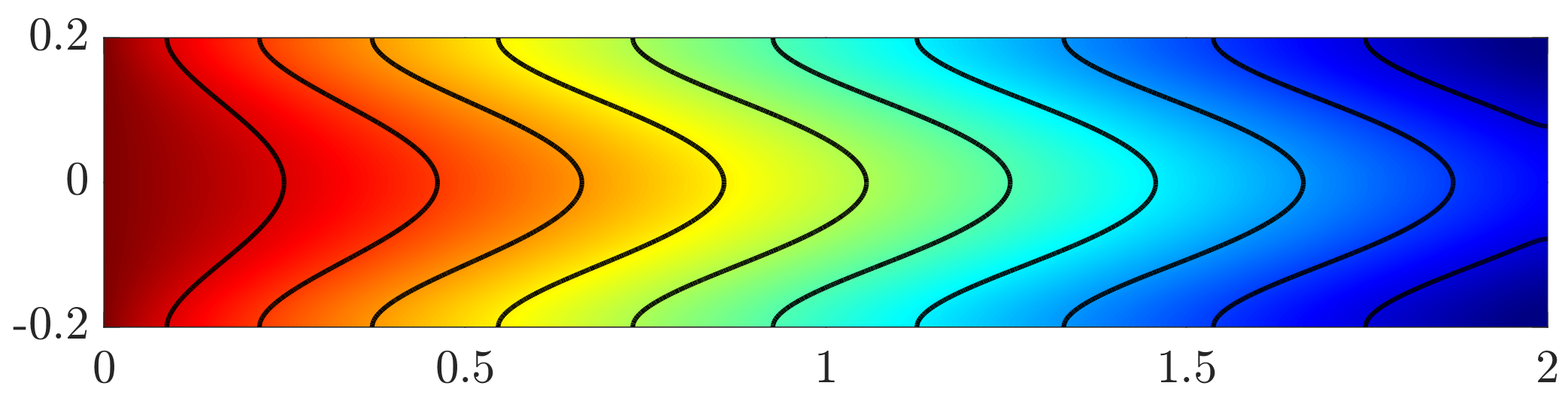}
    \hspace*{-.1cm}
    \includegraphics[height=1.55cm]{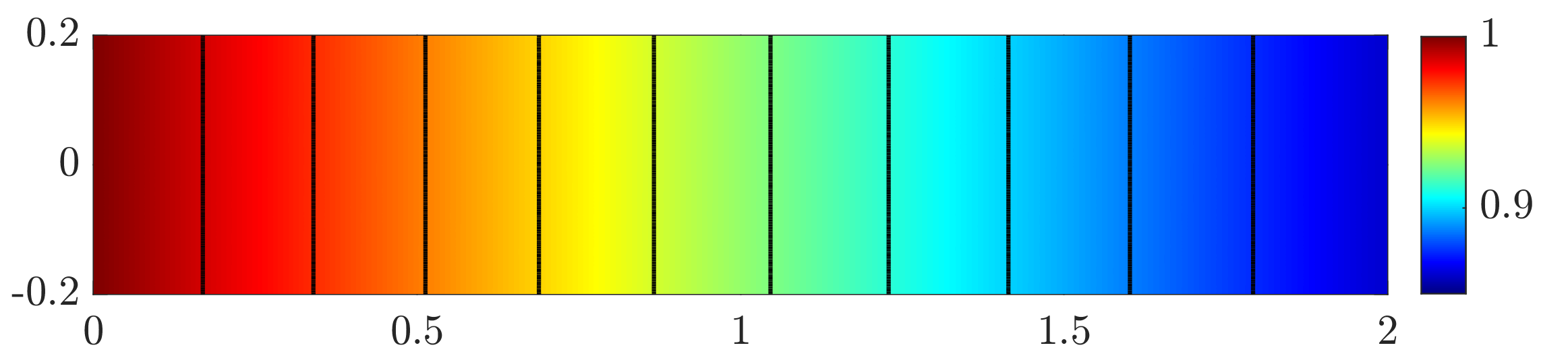}
    \hspace*{-.3cm}
    \includegraphics[height=1.55cm]{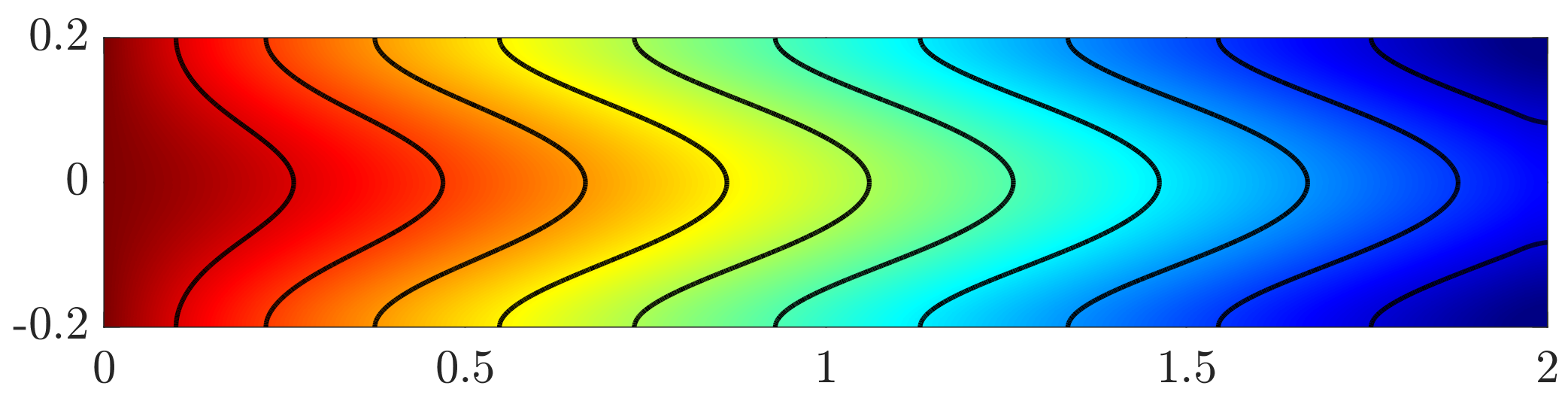}
    \hspace*{-.1cm}
    \includegraphics[height=1.55cm]{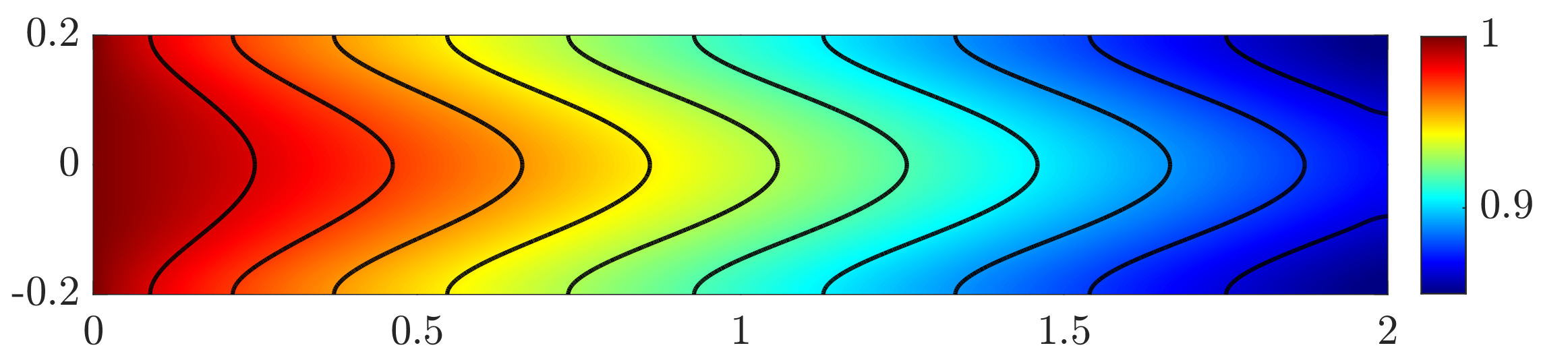}
     \caption{Poiseuille flow test case: contour plot of the reference solution (top left) and of the HiPhom$\varepsilon$ approximation for $m = 1$ (top-right), and $m=2$, $3$ (bottom, left and right).}
     \label{fig:respoinewnew_HiPhome}
   \end{figure}
\begin{figure}[!htbp]
    \centering
    \hspace*{-.3cm}
    \includegraphics[height=1.55cm]{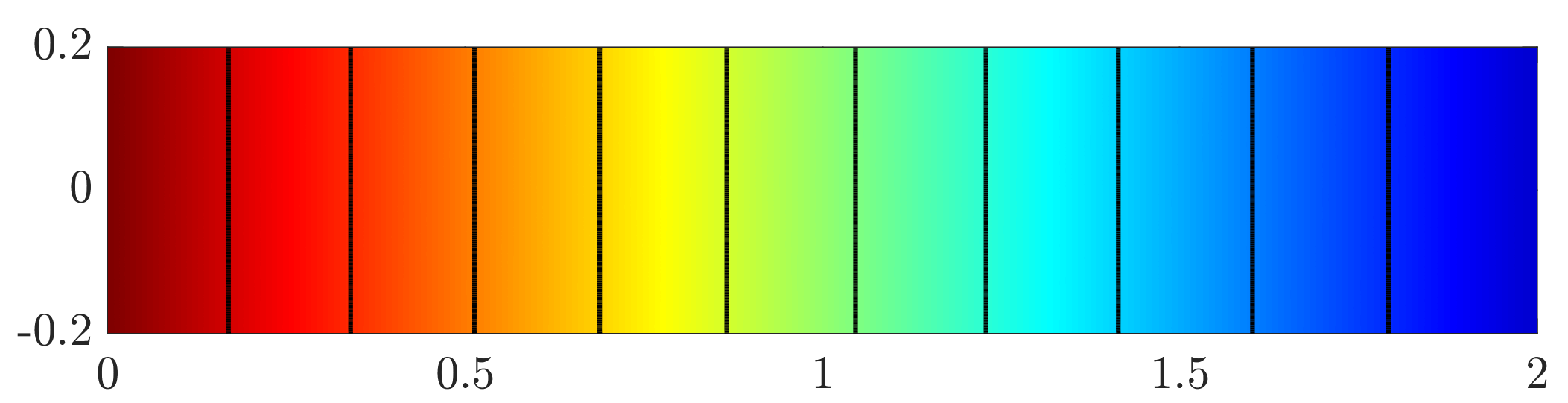}
    \hspace*{-.1cm}
    \includegraphics[height=1.55cm]{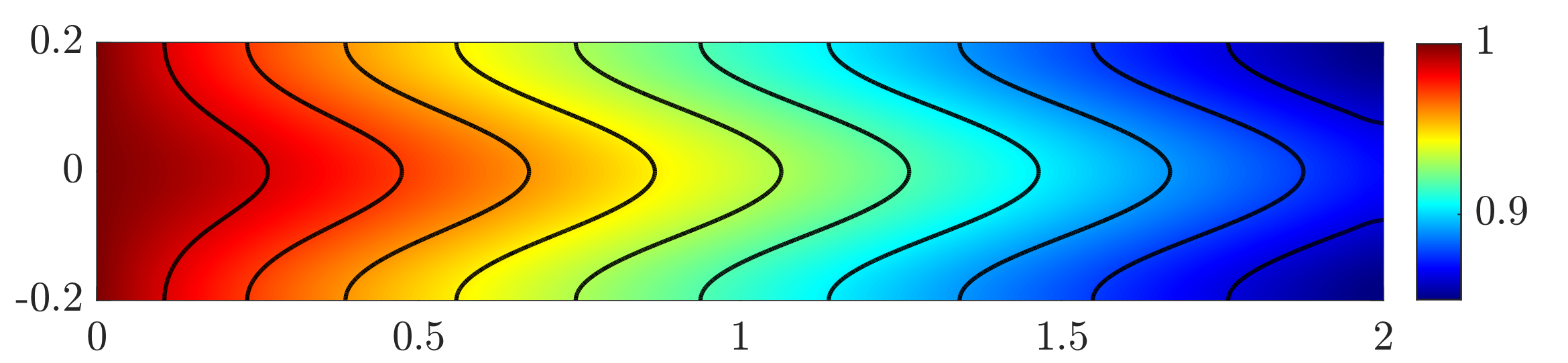}
    \hspace*{-.3cm}
    \includegraphics[height=1.55cm]{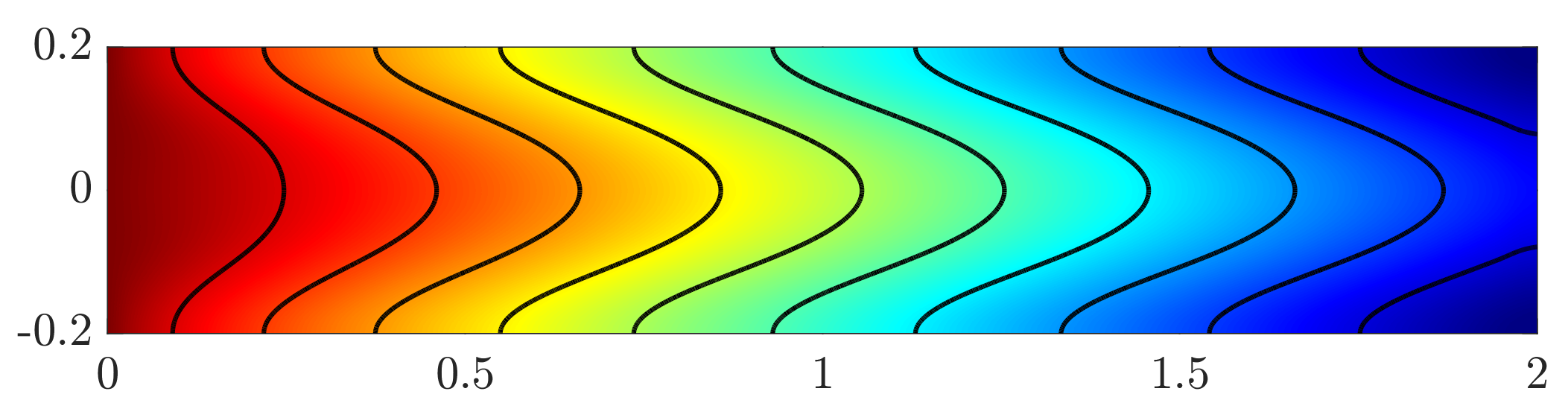}
    \hspace*{-.1cm}
    \includegraphics[height=1.55cm]{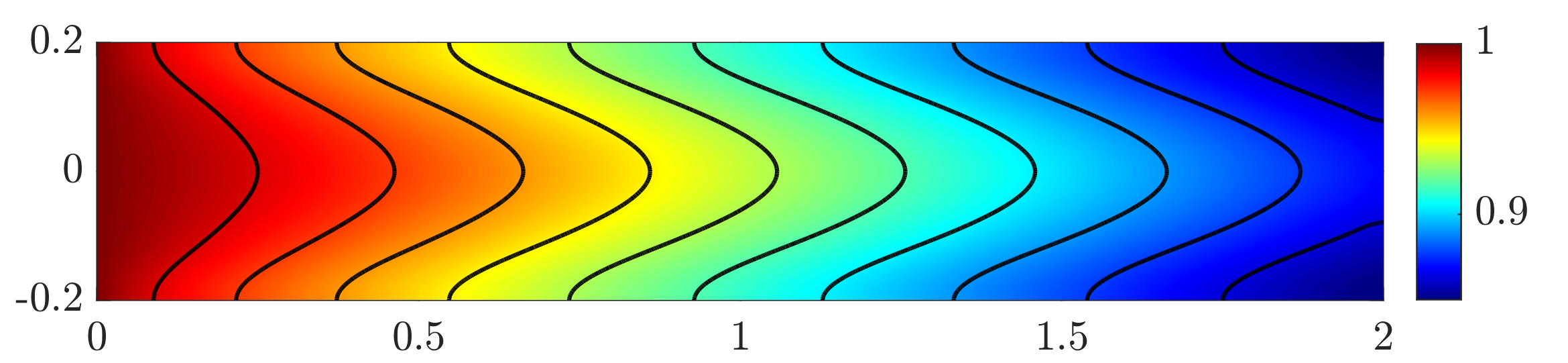}
     \caption{Poiseuille flow test case: contour plot of the HiMod approximation for $m = 1$, $2$ (top, left and right), and $m=3$, $4$ (bottom, left and right).}
     \label{fig:respoinewnew_HiMod}
   \end{figure}
   
For comparison purposes, we adopt as reference solution 
the (full) approximation to problem \eqref{eq:6} yielded by a linear finite element scheme associated with 
an unstructured uniform tessellation of $\Omega_\varepsilon$ consisting of $338202$ triangles (see the top-left panel in Figure~\ref{fig:respoinewnew_HiPhome}). 
 It is evident that the concentration field deforms according to the imposed parabolic velocity profile.
 
For the HiPhom$\varepsilon$ and HiMod reduction procedures we approximate the main dynamics by employing linear finite elements based on a uniform partition of the supporting fibre $\Omega_\varepsilon^{1D}$, with $h=0.0125$, while resorting to an increasing number of modes to describe the transverse dynamics. 
Figures~\ref{fig:respoinewnew_HiPhome} and~\ref{fig:respoinewnew_HiMod} show the HiPhom$\varepsilon$ and the HiMod approximation, respectively, for different choices of the modal index $m$. From a qualitative viewpoint, the solutions yielded by the two reduction procedures are very similar, for a fixed number of modes. Figures~\ref{fig:errpoinewnew} and \ref{fig:errhimod} offer a more quantitative comparison between the HiPhom$\varepsilon$ and the HiMod approximations, by showing
the spatial distribution of the 
absolute error between the reference finite element approximation and the HiPhom$\varepsilon$ and HiMod solution, respectively. Notably, for HiMod we report results obtained with odd values of $m$ to consider solutions embedding an increasing number of symmetric modes. 
The HiPhom$\varepsilon$ error is mainly localized close to the inlet and the outlet boundaries of the domain (see Figure \ref{fig:errpoinewnew}). 
As expected, the discrepancy between the full and the reduced solution decreases when increasing the modal index. For instance, 
when $m=2$ some visible difference in the HiPhom$\varepsilon$ solution curvature as compared to the reference solution (compare the top-left and the bottom-left panels in Figure~\ref{fig:respoinewnew_HiPhome}). This is reflected by values of the local error of the order of $10^{-3}$ close to inlet and outlet boundaries (see Figure \ref{fig:errpoinewnew} top right). When increasing the number of modes, the error progressively reduces to the order $10^{-5}$ in the majority of the spatial domain (see Figure \ref{fig:errpoinewnew} bottom right).\\
Figure~\ref{fig:errhimod} shows that the HiMod reduction requires additional modal functions with respect to the HiPhom$\varepsilon$ approach to reach a certain accuracy. Notably, employing $7$ educated modal functions yield a modelling error which on average is $\mathcal{O}(10^{-4})$. 
Moreover, the error associated with the HiMod reduction is
essentially independent of $x$. In contrast, the HiPhom$\varepsilon$ approximation displays an uneven distribution of errors along the $x$-axis with maximum values localized close to the inlet and outlet boundaries.

\subsubsection{Convergence analysis} \label{sec:erranapoi}
In this section, we numerically investigate the convergence of the error associated with the HiPhom$\varepsilon$ approximation, a theoretical analysis being beyond the goal of the paper. To this aim,
we evaluate the modelling error, $\|c - c_{m, \varepsilon} \|_{L^2{(\Omega_\varepsilon)}}$, with respect to the $L^2(\Omega_\varepsilon)$-norm as well as the Quantity of Interest (QoI) $\| c \|_{L^2(\Omega_\varepsilon)}$, by defining the error $J(c, c_{m, \varepsilon}) :=  | \|c\|_{L^2{(\Omega_\varepsilon)}} - \|c_{m, \varepsilon}\|_{L^2{(\Omega_\varepsilon)}} |$. The selected QoI is
directly related to the scalar dissipation rate, which measures the solute mixing within the domain~\cite{bolster2011mixing}.
%
\begin{figure}[!htb]
     \hspace*{-0.35cm}\includegraphics[width=0.52\linewidth]{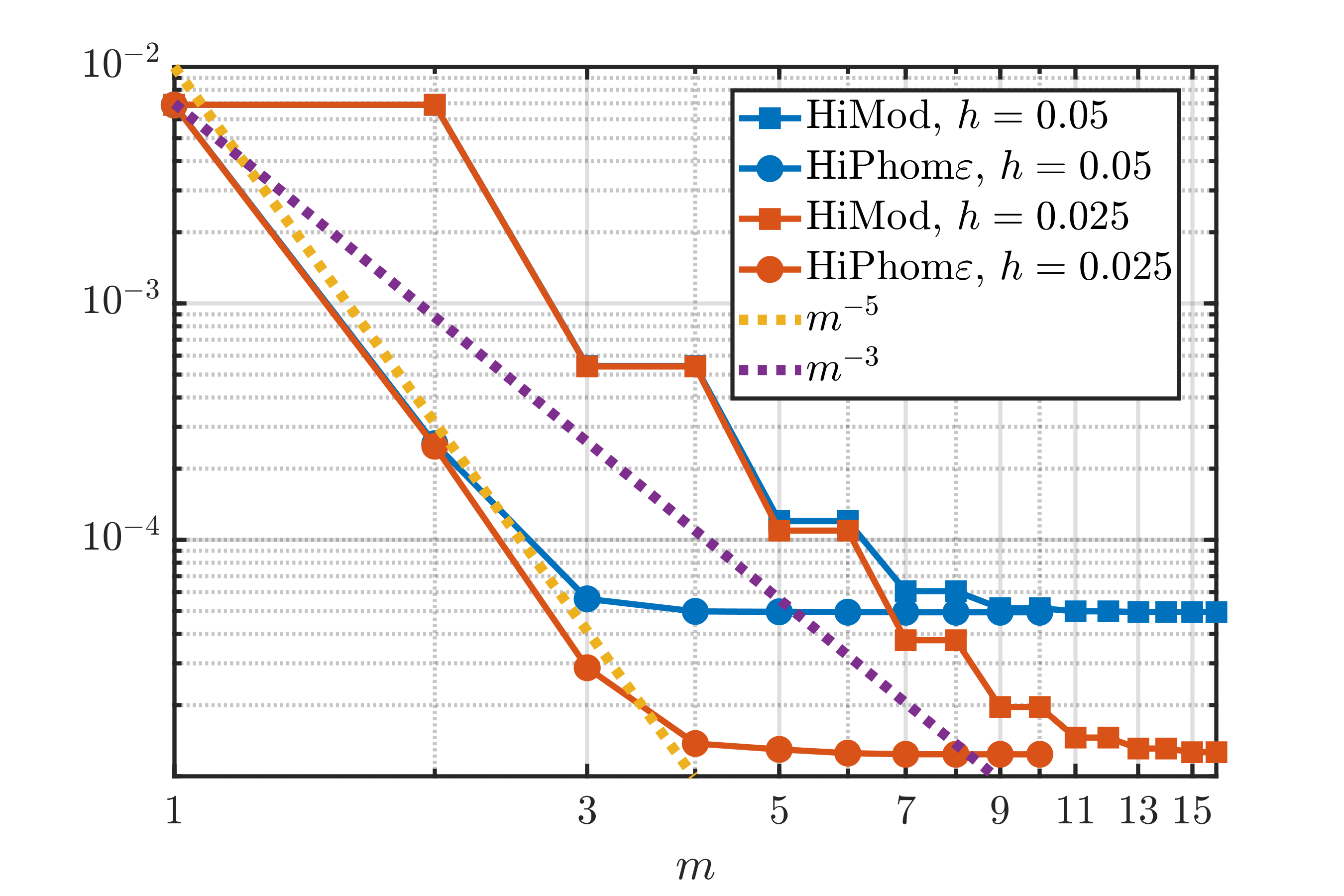}
    \hspace*{-0.35cm}\includegraphics[width=0.52\linewidth]{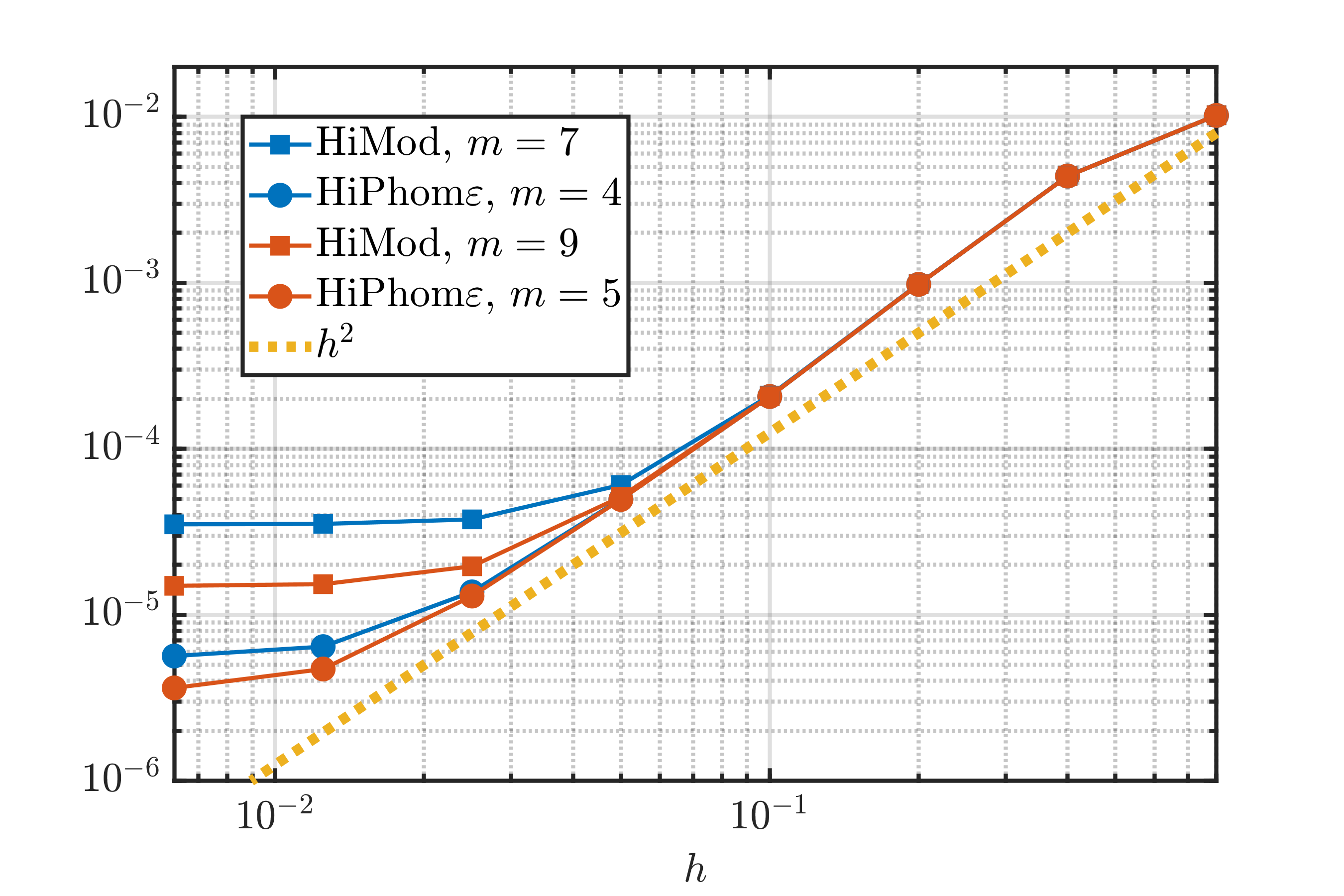}
    \caption{Poiseuille flow test case: 
    $L^2(\Omega_\varepsilon)$-norm of the modelling error associated with the HiPhom$\varepsilon$ and the HiMod approximations as a function of $m$ for different choices of $h$ (left) and as a function of $h$ for different choices of $m$ (right). The dot lines provide a reference trend for the convergence rate.}
    \label{fig:errbym}
\end{figure}

Figure~\ref{fig:errbym}, left shows the $L^2(\Omega_\varepsilon)$-norm of the modelling error associated both with the HiPhom$\varepsilon$ and the HiMod discretisation as a function of the number of modes and for two different choices of the step size used to discretise the supporting fibre $\Omega_\varepsilon^{1D}$. 
As expected the accuracy of both the reduced solutions decreases with $h$ and increases with $m$, as detected by the error trends. When increasing the number of modes $m$ the error values are independent of the grid size $h$ for small $m$, but tend to display an $h$-dependent asymptotic error value attained for large number of modes (large $m$). This indicates a stagnation of the approximation error taking place when the longitudinal approximation error is dominant over the modal one. When considering the same value of $h$, the asymptotic error value is the same for both approaches. However, HiPhom$\varepsilon$ converges faster to the asymptotic value, as the method exhibits a faster convergence as compared to HiMod. 
Finally, the step-wise pattern trend characterising the HiMod error highlights the redundancy of the non-symmetric HiMod basis functions, i.e. error drops are associated only with odd values of $m$, these latter corresponding to the occurrence of symmetric basis functions (see Figure~\ref{basis_fs}).
\\
The right panel in Figure~\ref{fig:errbym} displays the
$L^2(\Omega_\varepsilon)$-norm of the HiPhom$\varepsilon$ and of the HiMod error as a function of the discretisation size $h$, when considering a different number of modal functions. This time we observe a stagnation of the error value for small $h$, where the error is dominated by inaccuracies in capturing the transverse behaviour. The HiPhom$\varepsilon$ allows attaining smaller errors than HiMod when $h \to 0$ with smaller values of $m$, thus displays an optimal performance. For sufficiently large values of $h$, the order of convergence turns out to be quadratic.
\begin{figure}[h!tbp]
    \hspace*{-0.35cm}\includegraphics[width=0.52\linewidth]{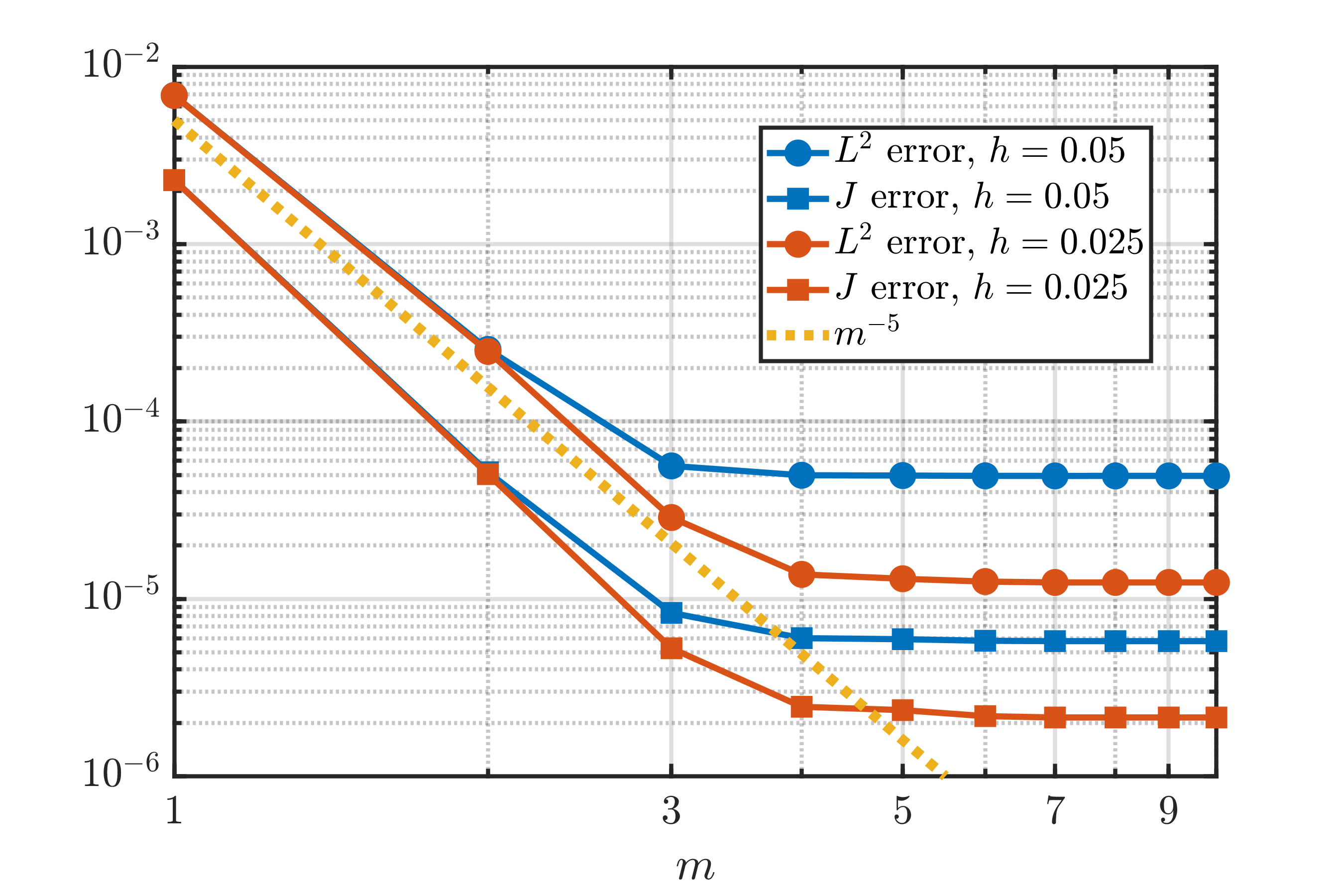}
    \hspace*{-0.35cm}\includegraphics[width=0.52\linewidth]{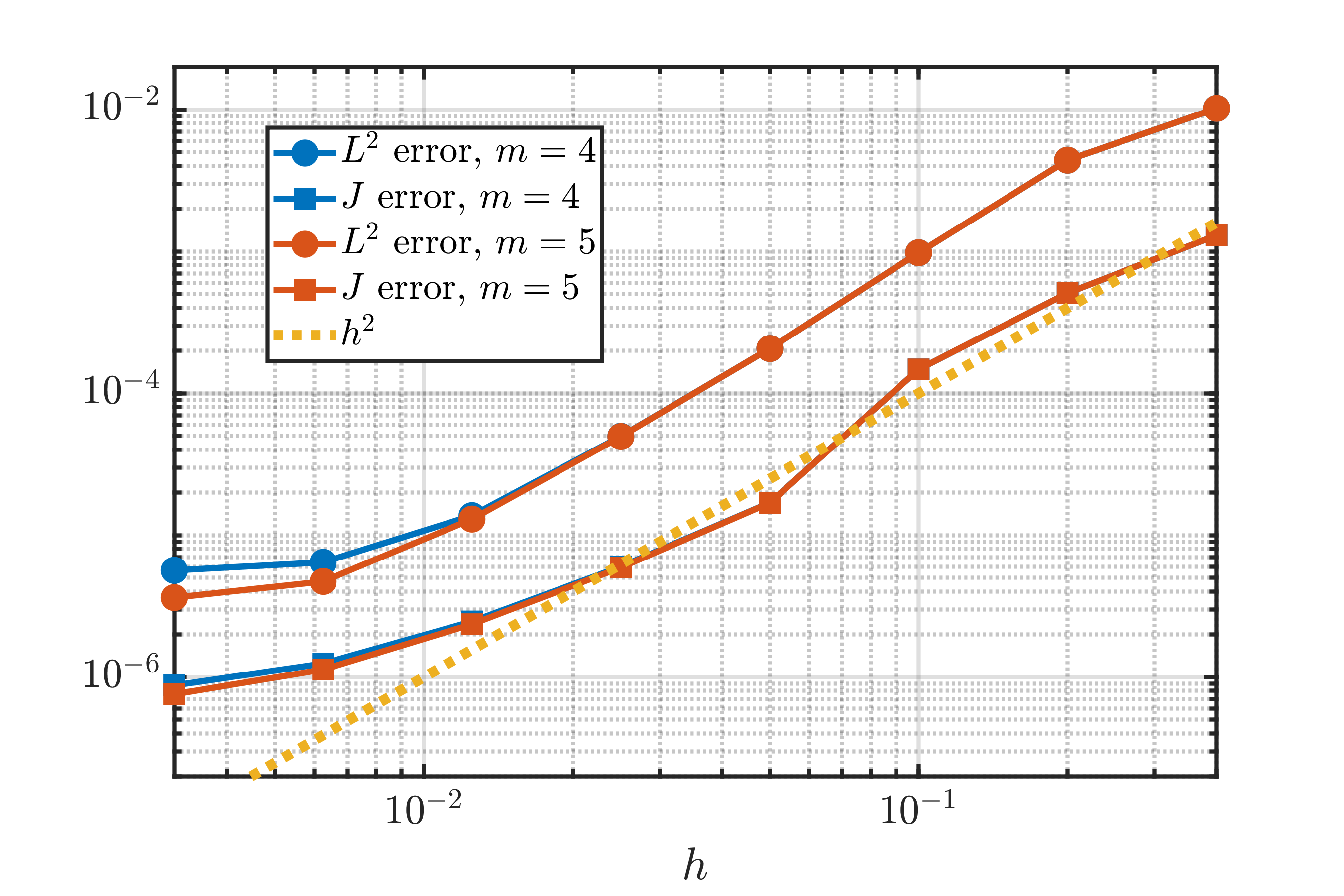}
    \caption{Poiseuille flow test case: 
    error associated with the HiPhom$\varepsilon$ approximation with respect to the $L^2(\Omega_\varepsilon)$-norm and the QoI as a function of $m$ and for different choices of $h$ (left) and as a function of $h$ for different choices of $m$ (right). The dot lines provide a reference trend for the convergence rate.}
    \label{fig:j2errfirst}
\end{figure}

The left panel of Figure~\ref{fig:j2errfirst} shows the trend of $J$, i.e. the error associated with the QoI, as a function of $m$ and for two choices of $h$. The error quickly decreases when increasing $m$ and reducing $h$, before stagnating towards a constant value. This confirms the reliability 
of the HiPhom$\varepsilon$ approach also in reproducing physically meaningful quantities for the application at hand.
For comparison purposes, we also show the trend of the $L^2(\Omega_\varepsilon)$-norm of the HiPhom$\varepsilon$ error, for the same choices of $m$ and $h$. The rate of convergence characterising the two different error norms is very similar (about $\mathcal{O}(m^{-5})$). In particular, 
the HiPhom$\varepsilon$ approximation turns out to be particularly effective in surrogating the QoI 
(with a gain of about one order of accuracy with respect to the evaluation of the $L^2(\Omega_\varepsilon)$-norm of the modelling error).\\
The right panel in Figure~\ref{fig:j2errfirst} shows the reduction of the HiPhom$\varepsilon$ error evaluated 
as a function of $h$ and for $m=4$ and $m=5$. The choices for $m$ are associated with a stagnation trend in the left panel, for $h=0.05$ and $h=0.025$. The plot on the right highlights that an improvement in terms of accuracy can be achieved by further reducing $h$, although such a gain remains limited. Finally, the rate of convergence is quadratic for both the considered error norms.

\subsection{A boundary layer flow}\label{sec:experiments2}
\begin{figure}[!htb]
\centering
    \includegraphics[width=0.5\textwidth]{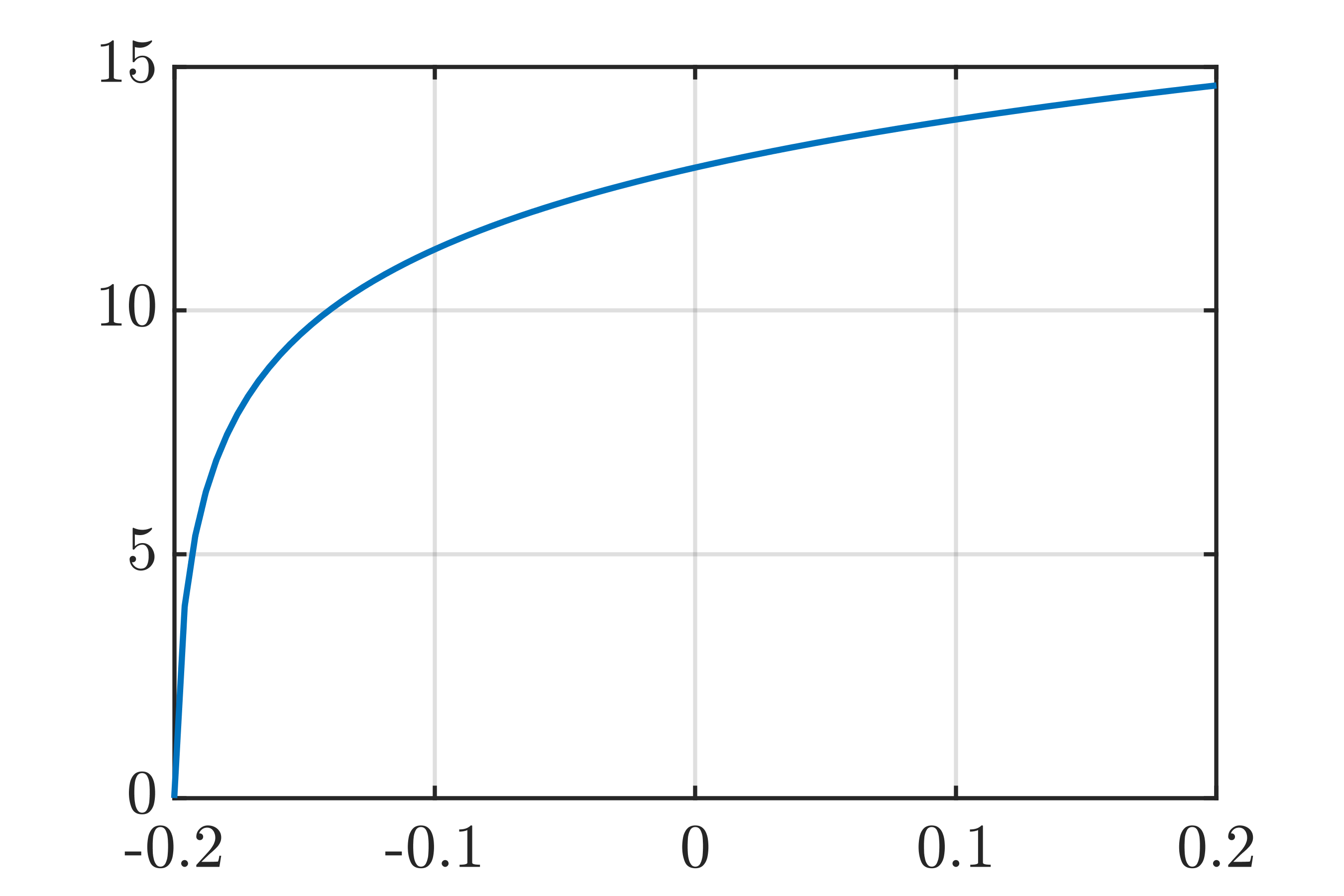}
    \caption{The boundary layer flow test case: log-law velocity profile.}
    \label{fig:low_law}
\end{figure}

We now test the HiPhom$\varepsilon$ method 
by modelling an ADR problem in a stationary configuration characterized by a log-law velocity profile. This advective field is usually employed to approximate the vertical profile of a stream-wise velocity in an open channel flow~\cite{TimPr:07,Stroo:02}.\\
Thus, we select the non-symmetric velocity profile
\begin{equation}\label{eq:boundary1}
    \mathbf{u}=[u(z), 0]^T = \left[\frac{1}{k}\ln(z+\varepsilon+d) + C, 0\right]^T,
\end{equation}
with $k=0.41$ the Von Kármán constant, $\varepsilon=0.2$, $d = 0.001$ and $C = -\log( d )/k$ (see Figure~\ref{fig:low_law}). Problem \eqref{eq:6} is approximated in the 
domain $\Omega_{\varepsilon} = (0,2) \times (-\varepsilon,\varepsilon)$, by setting $\sigma=1$, $f=0$, $c_B=1$, and $D=1$ so that $D_\varepsilon=0.2$. Finally, the boundary portions coincide with $\Gamma_{\varepsilon}^i = \{x=0\}\times (-\varepsilon/2,\varepsilon/2)$, $\Gamma_{\varepsilon}^o = \{x=2\}\times (-\varepsilon/2,\varepsilon/2)$, and  $\Gamma_{\varepsilon}^w = [0,2]\times \{\varepsilon/2,\varepsilon/2\}$.

\begin{figure}[tb]
\hspace*{-.8cm}
    \includegraphics[width=0.55\textwidth]{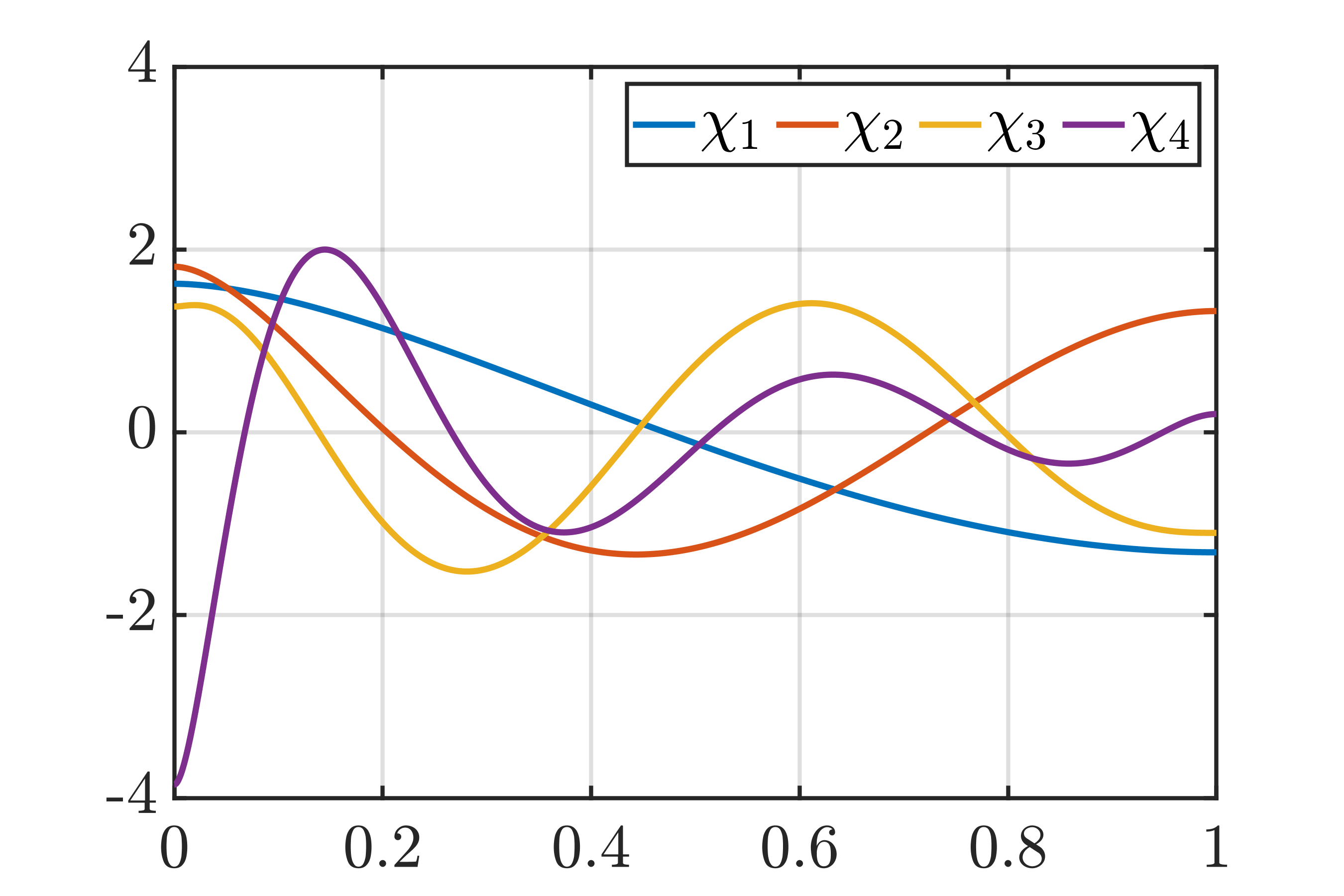}
   \hspace*{-.8cm}
   \includegraphics[width=0.55\textwidth]{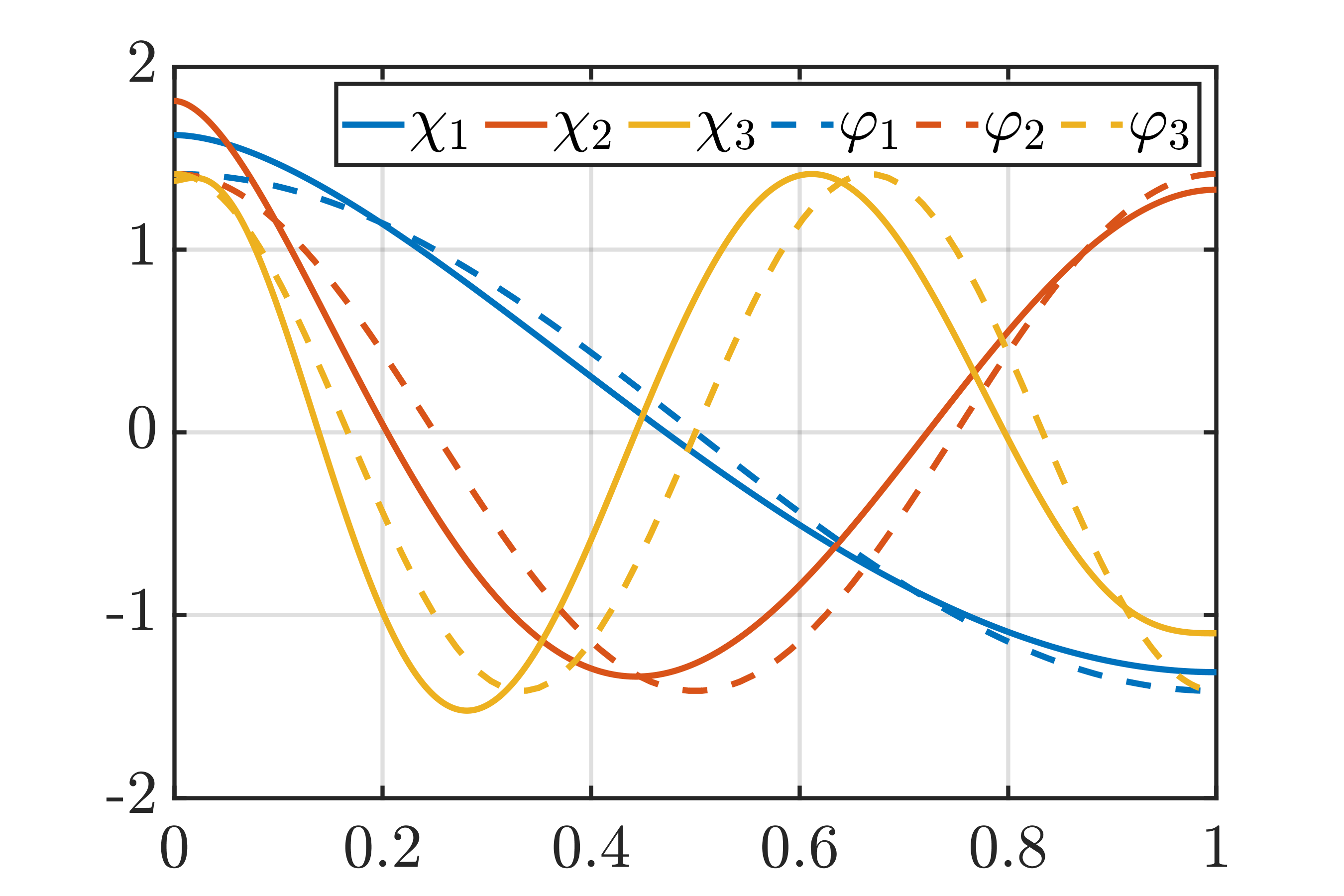}\label{fig:modalbasis}
    \caption{The boundary layer flow test case: HiPhom$\varepsilon$ modal basis functions $\chi_i$, for $i=1, \ldots, 4$ (left); comparison among the first $3$ HiPhom$\varepsilon$ and HiMod modes (right).}
    \label{fig:Tc2_data}
\end{figure}

Figure~\ref{fig:Tc2_data}, left  gathers the HiPhom$\varepsilon$ basis functions $\chi_i$, with $i=1\ldots, 4$, which properly take into account the non-symmetric shape of the velocity distribution.
In the panel on the right, we compare functions $\chi_1$, $\chi_2$ and $\chi_3$, with the corresponding HiMod educated modes, $\varphi_1$,  $\varphi_2$ and $\varphi_3$, obtained by setting $\mathcal L_s=D\partial_{yy}$ in \eqref{eq:SLE1}. The mismatch between functions $\chi_i$ and $\varphi_i$ 
becomes more evident for increasing values of $i$. In particular, 
for this case study, both the HiPhom$\varepsilon$ and the HiMod modes are non-symmetric functions, consistently with the flow features. 

\begin{figure}[tb]
    \centering
    \hspace*{-.3cm}
    \includegraphics[height=1.55cm]{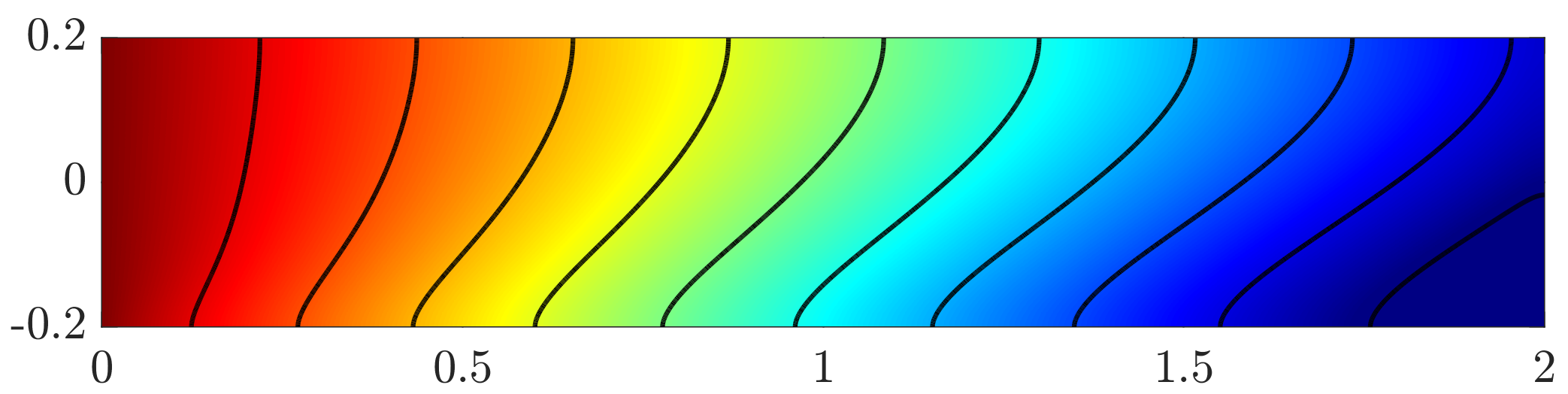}
    \hspace*{-.1cm}
    \includegraphics[height=1.55cm]{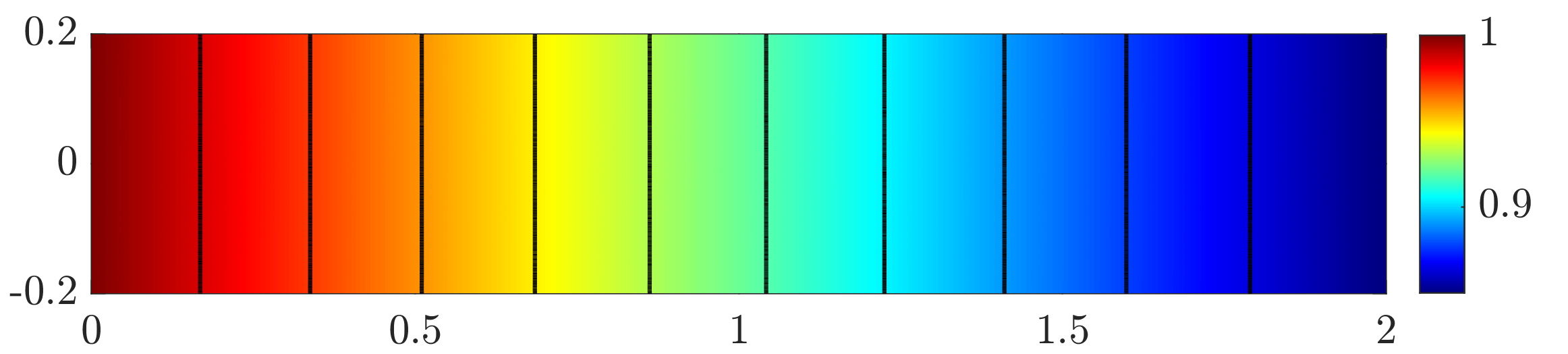}
    \hspace*{-.3cm}
    \includegraphics[height=1.55cm]{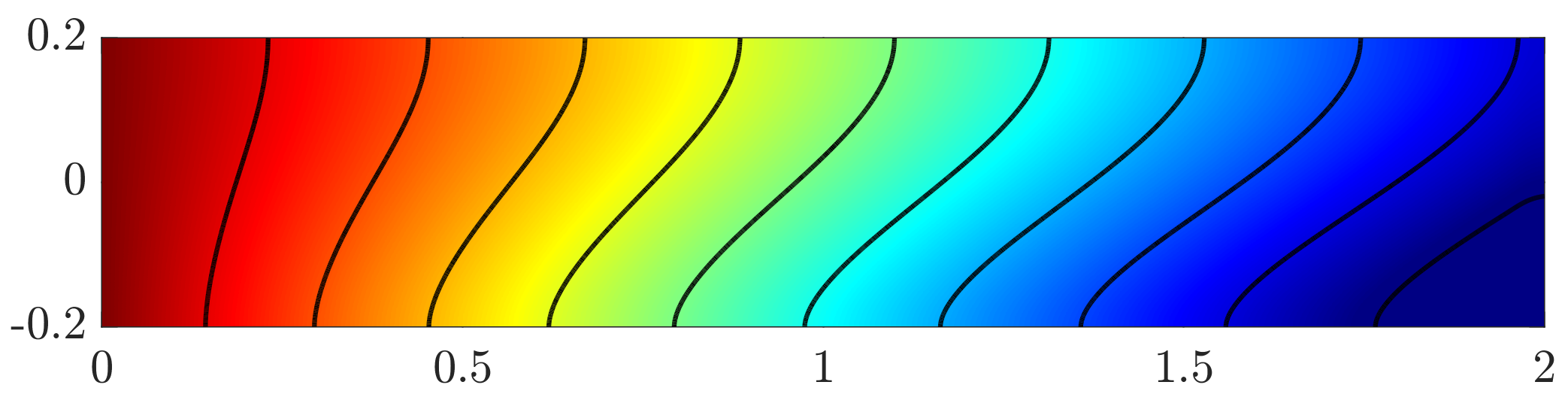}
    \hspace*{-.1cm}
    \includegraphics[height=1.55cm]{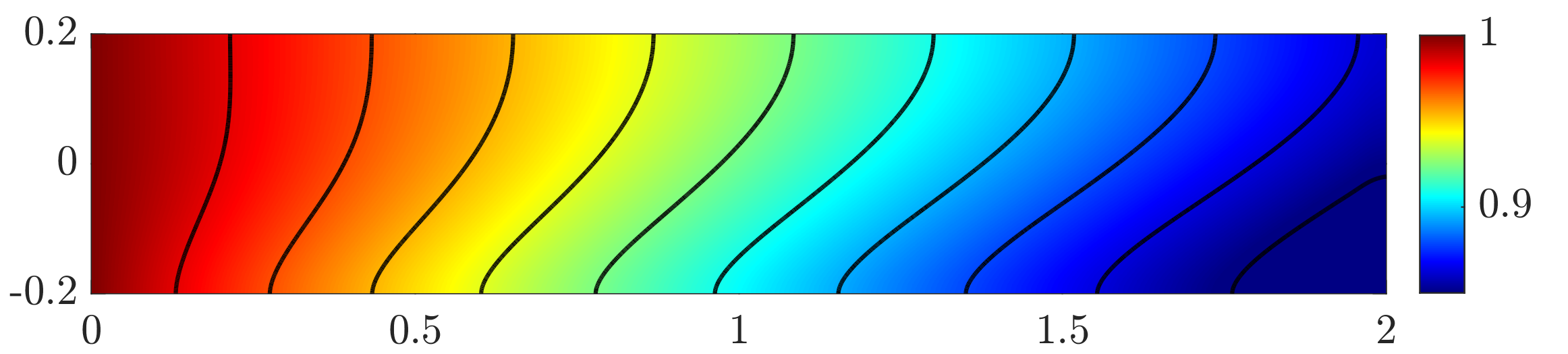}
     \caption{The boundary layer flow test case: 
     contour plot of the reference solution (top left) and of the HiPhom$\varepsilon$ approximation for $m = 1$ (top-right), and $m=2$, $3$ (bottom, left and right).}
\label{fig:logrespoinewnew_HiPhome}
   \end{figure}

To evaluate the reliability of the HiPhom$\varepsilon$ discretisation, we take as reference solution the (full) approximation to problem \eqref{eq:6} based on linear finite elements, associated with an unstructured mesh of $\Omega_\varepsilon$ constituted by 
$906640$ triangles (see Figure~\ref{fig:logrespoinewnew_HiPhome}, top-left).
The HiPhom$\varepsilon$ model reduction is carried out by discretising the leading dynamics with linear finite elements based on a uniform partition of the supporting fibre with step size $h=0.0125$, and by gradually increasing the cardinality of the HiPhom$\varepsilon$ modal basis. 
Figure~\ref{fig:logrespoinewnew_HiPhome}, top-right and bottom shows the HiPhom$\varepsilon$ discretisation for $m=1$, $2$ and $3$.
It turns out that $3$ HiPhom$\varepsilon$ modes suffice to have a qualitatively reliable reduced solution. 

\begin{figure}[tb]
    \centering
    \hspace*{-.3cm}
    \includegraphics[height=1.55cm]{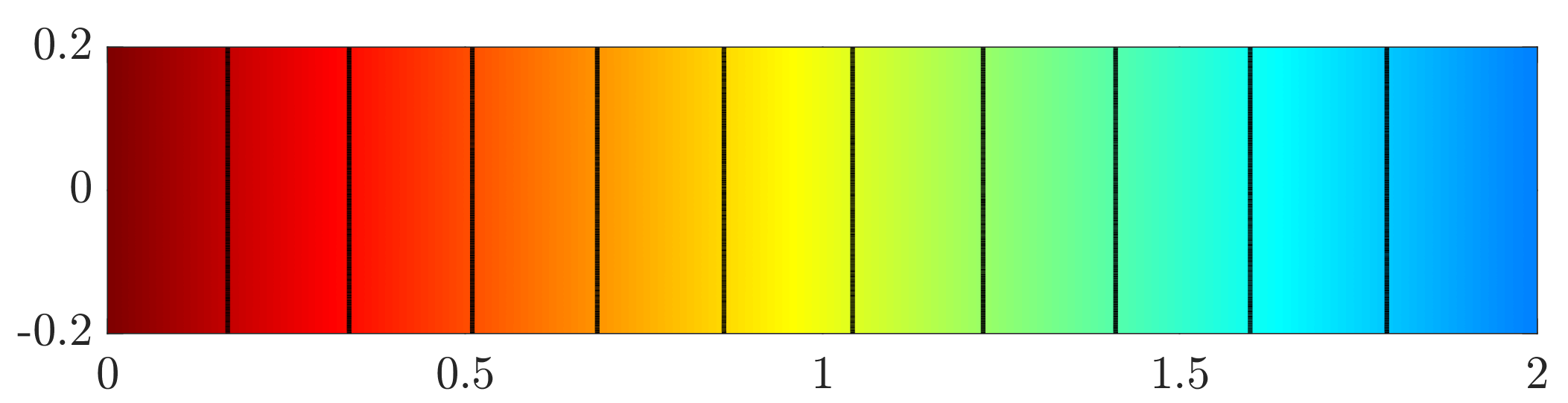}
    \hspace*{-.1cm}
    \includegraphics[height=1.55cm]{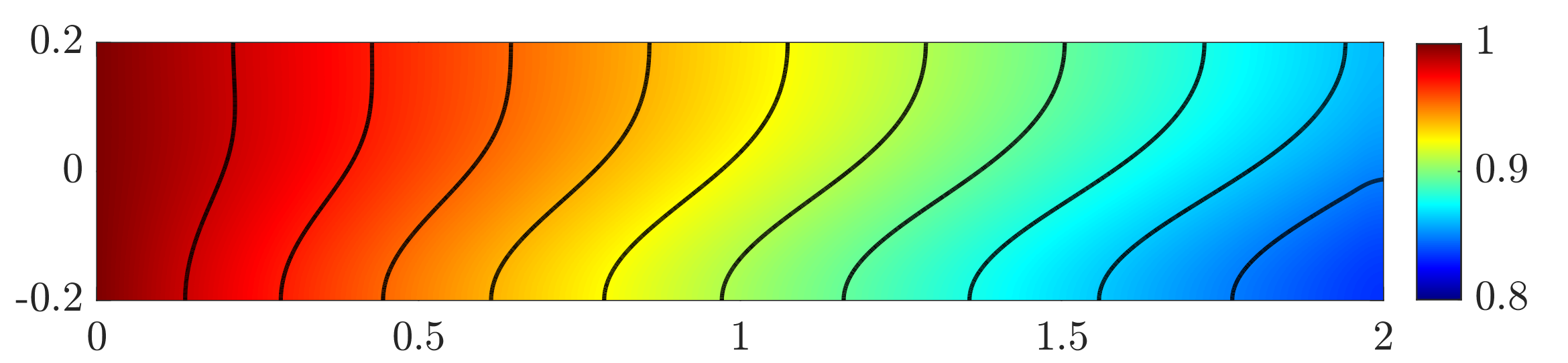}
    \hspace*{-.3cm}
    \includegraphics[height=1.55cm]{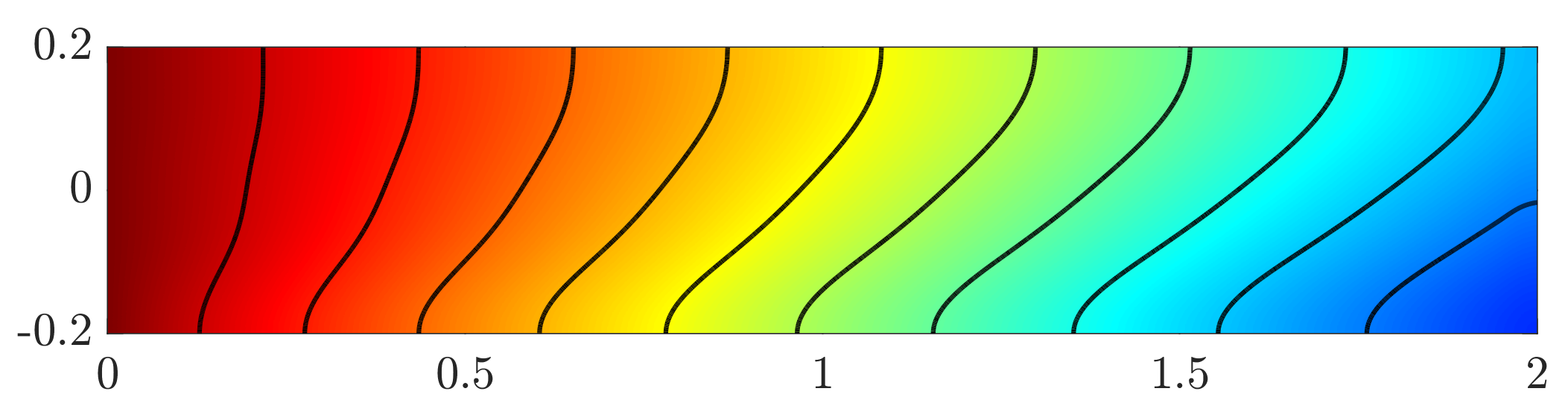}
    \hspace*{-.1cm}
    \includegraphics[height=1.55cm]{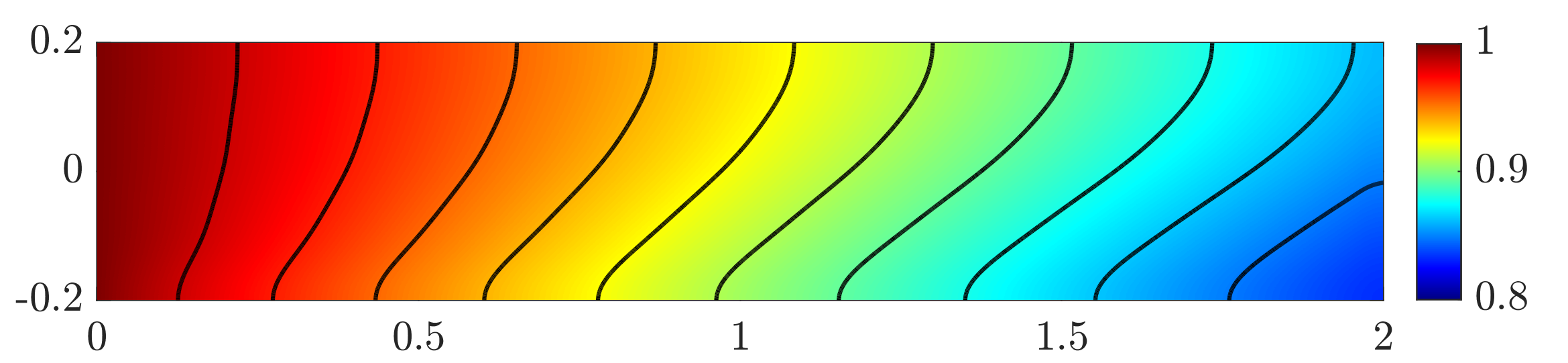}
     \caption{The boundary layer flow test case: contour plot of the HiMod approximation for $m = 1$, $2$ (top, left and right), and $m=3$, $4$ (bottom, left and right).}
     \label{fig:logrespoinewnew_HiMod}
   \end{figure}

\begin{figure}[!htbp]
    \centering
    \hspace*{-.3cm}
    \includegraphics[height=1.55cm]{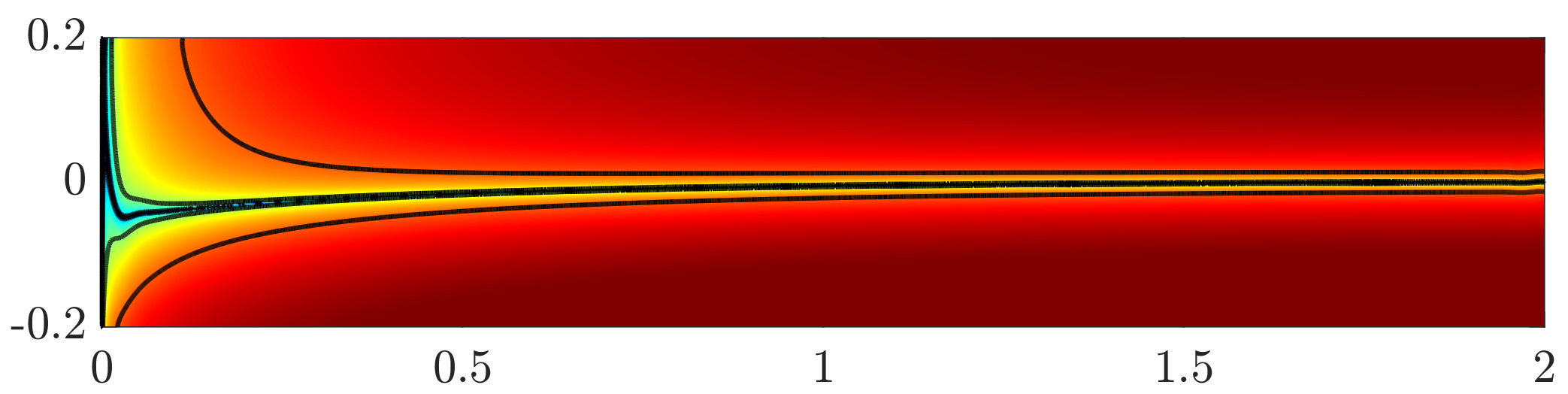}
    \hspace*{-.1cm}
    \includegraphics[height=1.55cm]{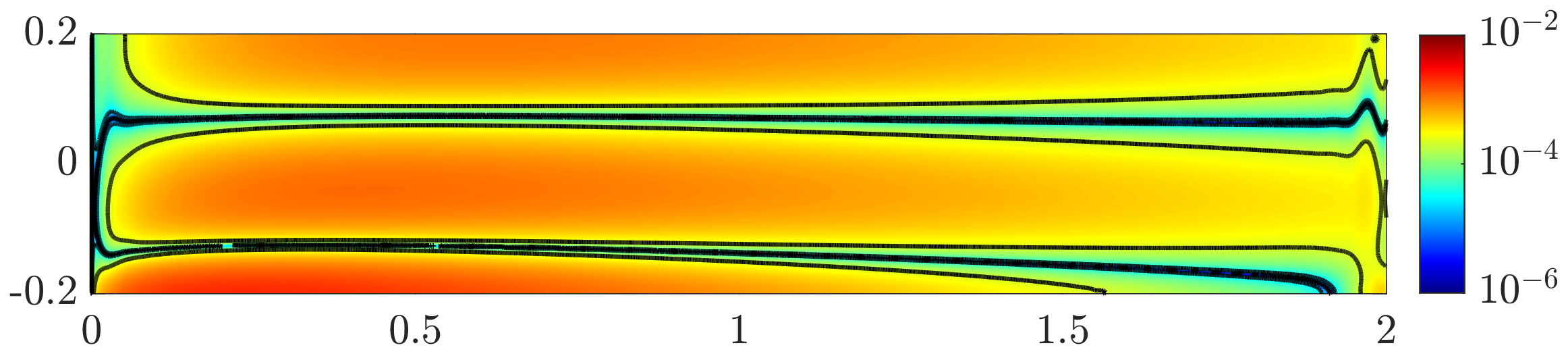}
    \hspace*{-.3cm}
    \includegraphics[height=1.55cm]{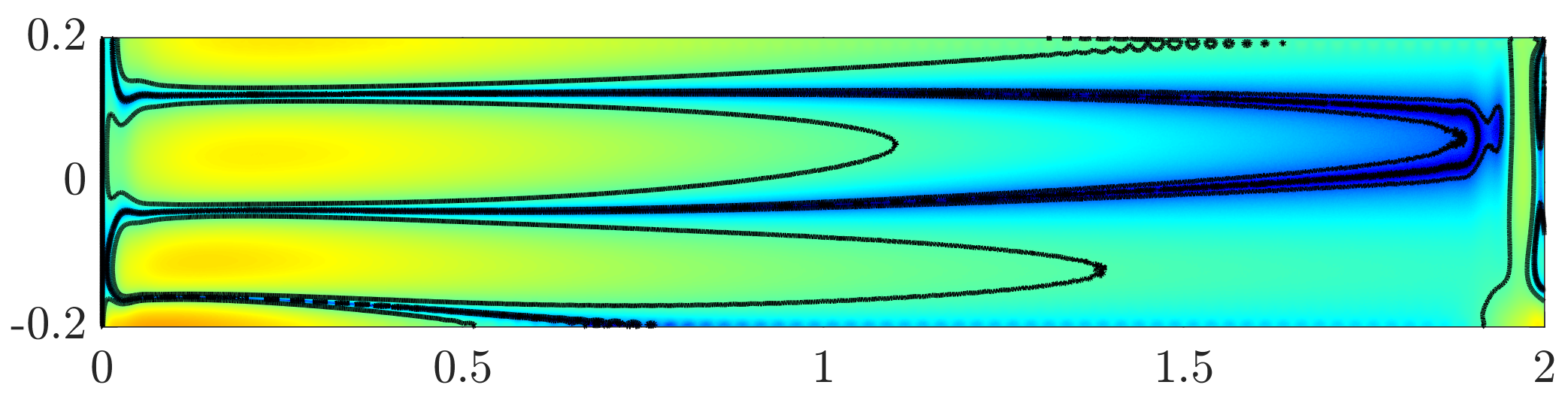}
    \hspace*{-.1cm}
    \includegraphics[height=1.55cm]{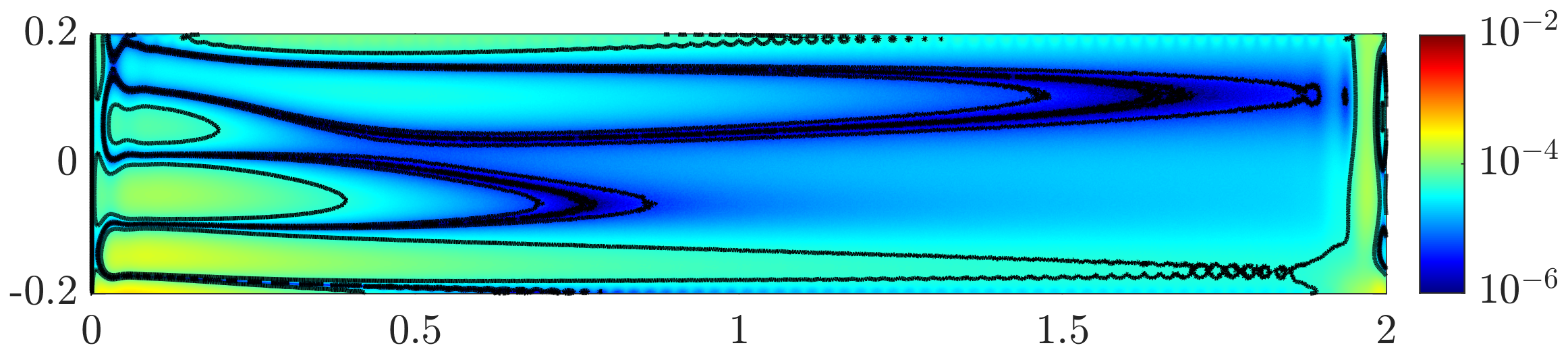}
     \caption{The boundary layer flow test case: spatial distribution of the absolute modelling error associated with the HiPhom$\varepsilon$ approximation for $m= 1$, $2$ (top, left and right) and $m=3$, $4$ (bottom, left and right).}
     \label{fig:logerrpoinewnew}
   \end{figure}

\begin{figure}[!htbp]
    \centering
    \hspace*{-.3cm}
    \includegraphics[height=1.55cm]{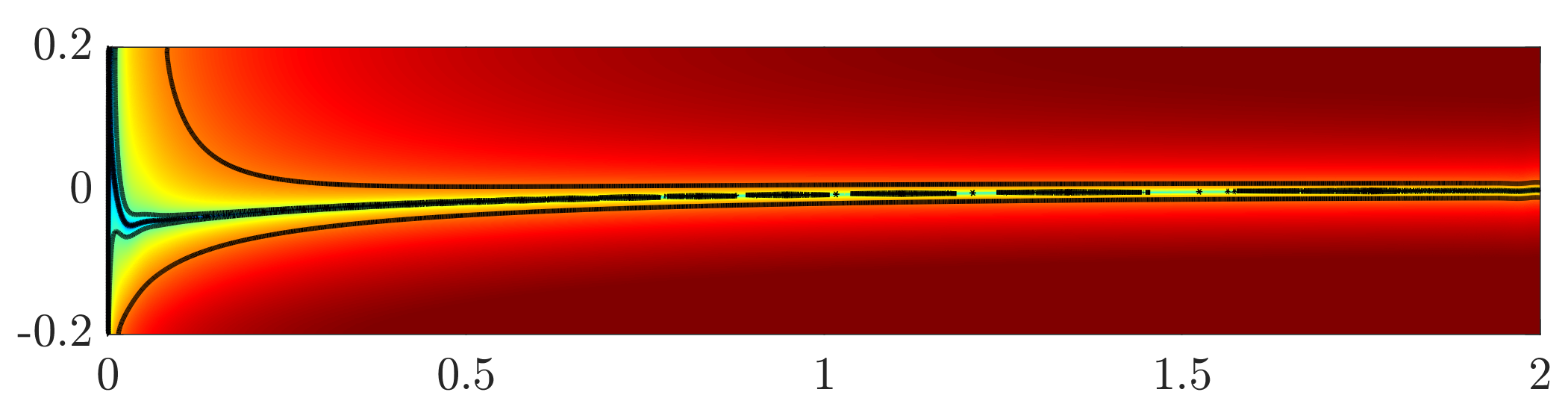}
    \hspace*{-.1cm}
    \includegraphics[height=1.55cm]{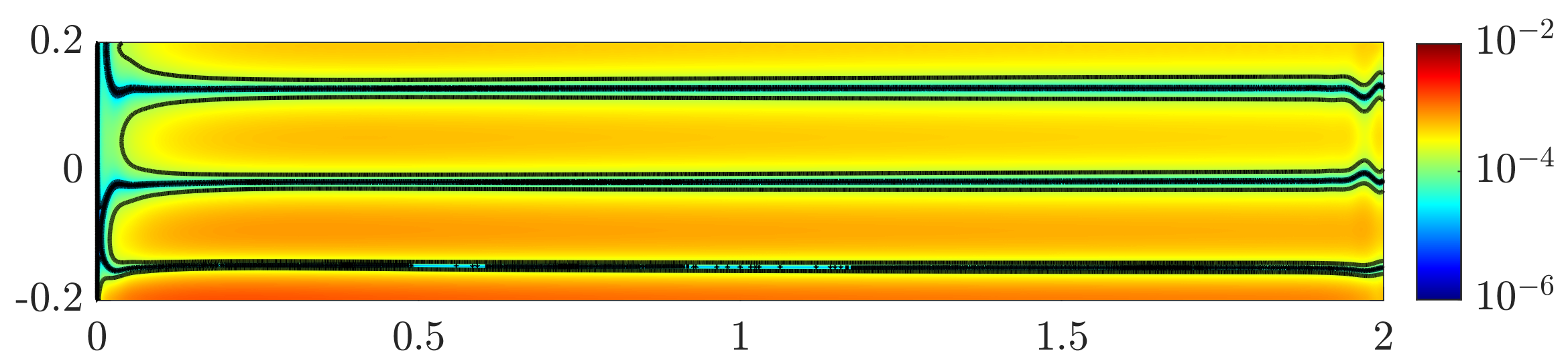}
    \hspace*{-.3cm}
    \includegraphics[height=1.55cm]{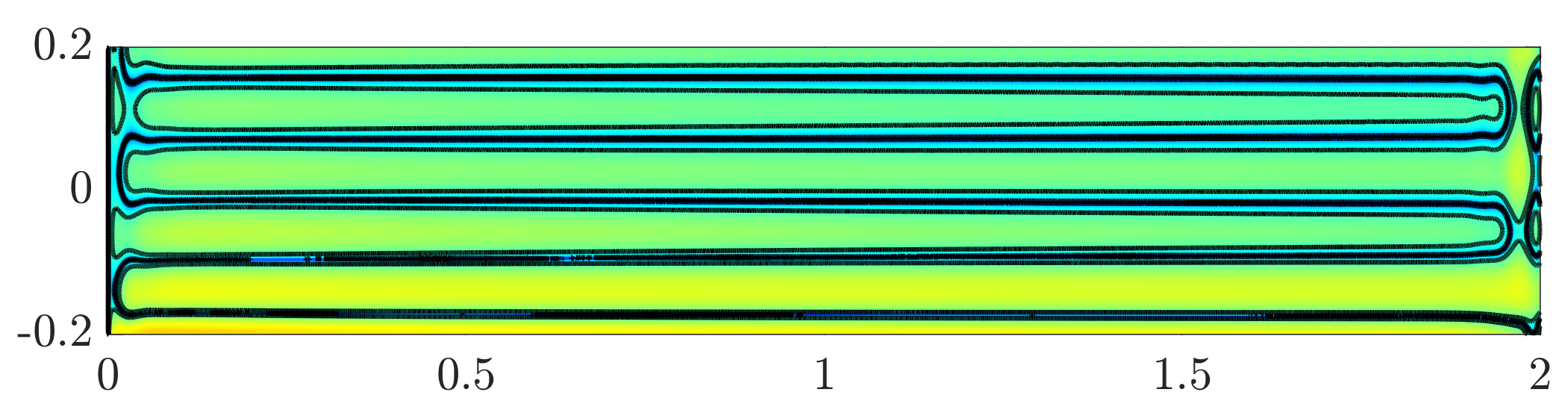}
    \hspace*{-.1cm}
    \includegraphics[height=1.55cm]{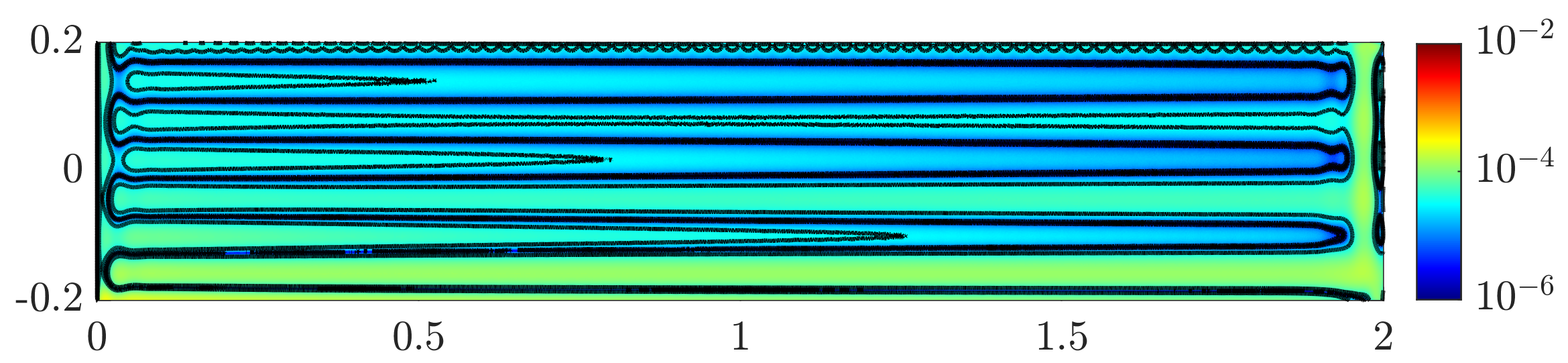}
     \caption{The boundary layer flow test case: spatial distribution of the absolute modelling error associated with the HiMod approximation for $m= 1$, $2$ (top, left and right) and $m=3$, $4$ (bottom, left and right).}
     \label{fig:logerrhimod}
   \end{figure}

Figure~\ref{fig:logrespoinewnew_HiMod} displays the HiMod approximation 
yielded by preserving the same
discretisation as for the HiPhom$\varepsilon$ reduction along $\Omega_\varepsilon^{1D}$, while using the educated modal functions in Figure~\ref{fig:Tc2_data}, right to model the transverse dynamics. The correct trend of functions $\varphi_i$ considerably promotes the reliability of the HiMod approximation, $4$ modes guaranteeing an accurate reduced solution, at least from a qualitative viewpoint. \\
To make the comparison between HiPhom$\varepsilon$ and HiMod more quantitative, we provide in Figures~\ref{fig:logerrpoinewnew} and~\ref{fig:logerrhimod} the spatial distribution of the absolute modelling error characterising the two reduction procedures. 
HiMod and HiPhom$\varepsilon$ yield errors of similar magnitude and distribution for $m=1$, $2$, whereas, for $m=3$, $4$, 
the spreading of the error across the whole domain exhibit a very different pattern for the two approaches, with a more regular distribution for the HiMod discretisation.

\subsubsection{Convergence analysis}
We replicate the analysis carried out in Section~\ref{sec:erranapoi}, with the goal of verifying whether a more challenging configuration compromises the decay of the error characterising the HiPhom$\varepsilon$ approach. 
\begin{figure}[tbp]
 \hspace*{-0.35cm}\includegraphics[width=0.52\linewidth]{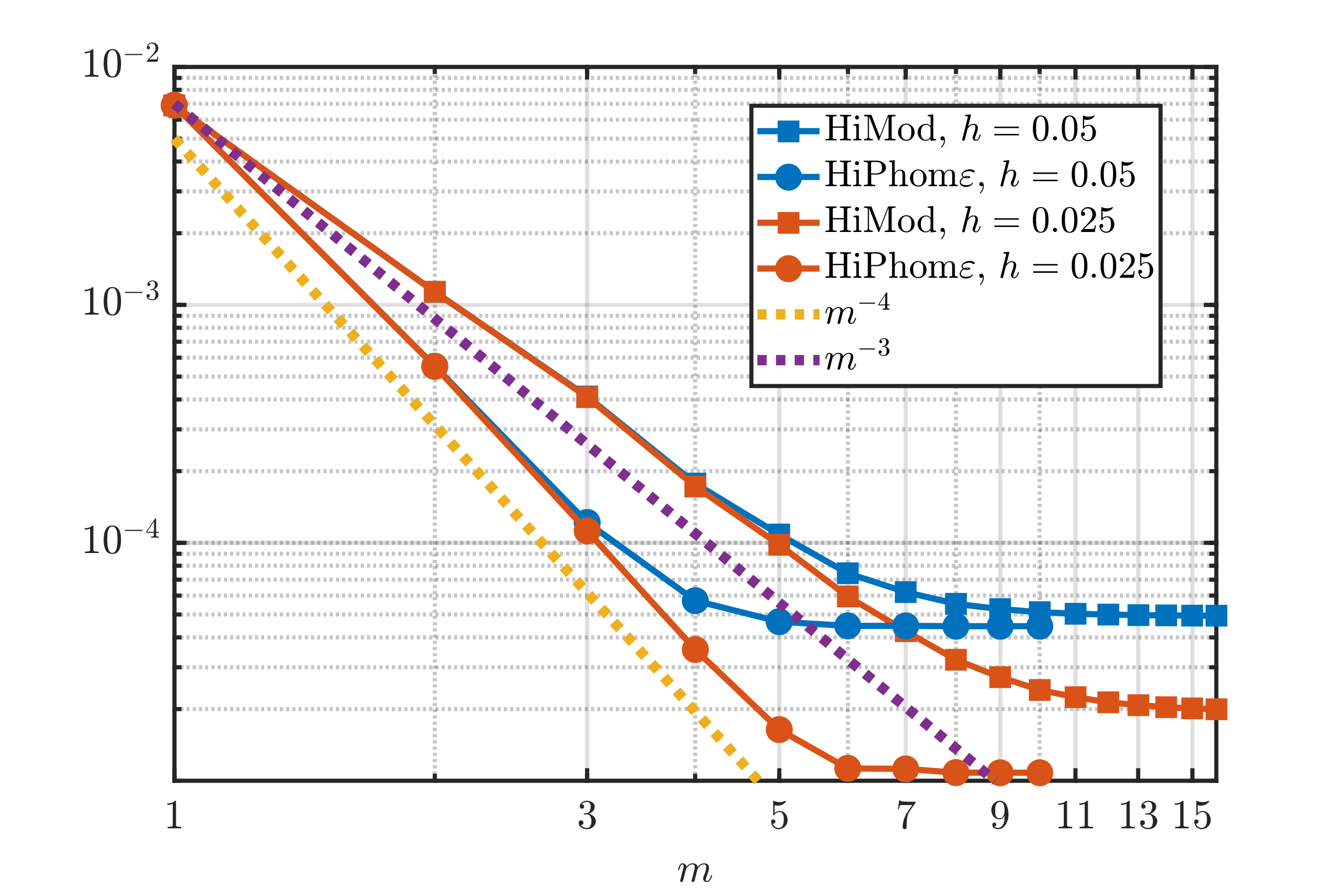}
    \hspace*{-0.35cm}   \includegraphics[width=0.52\linewidth]{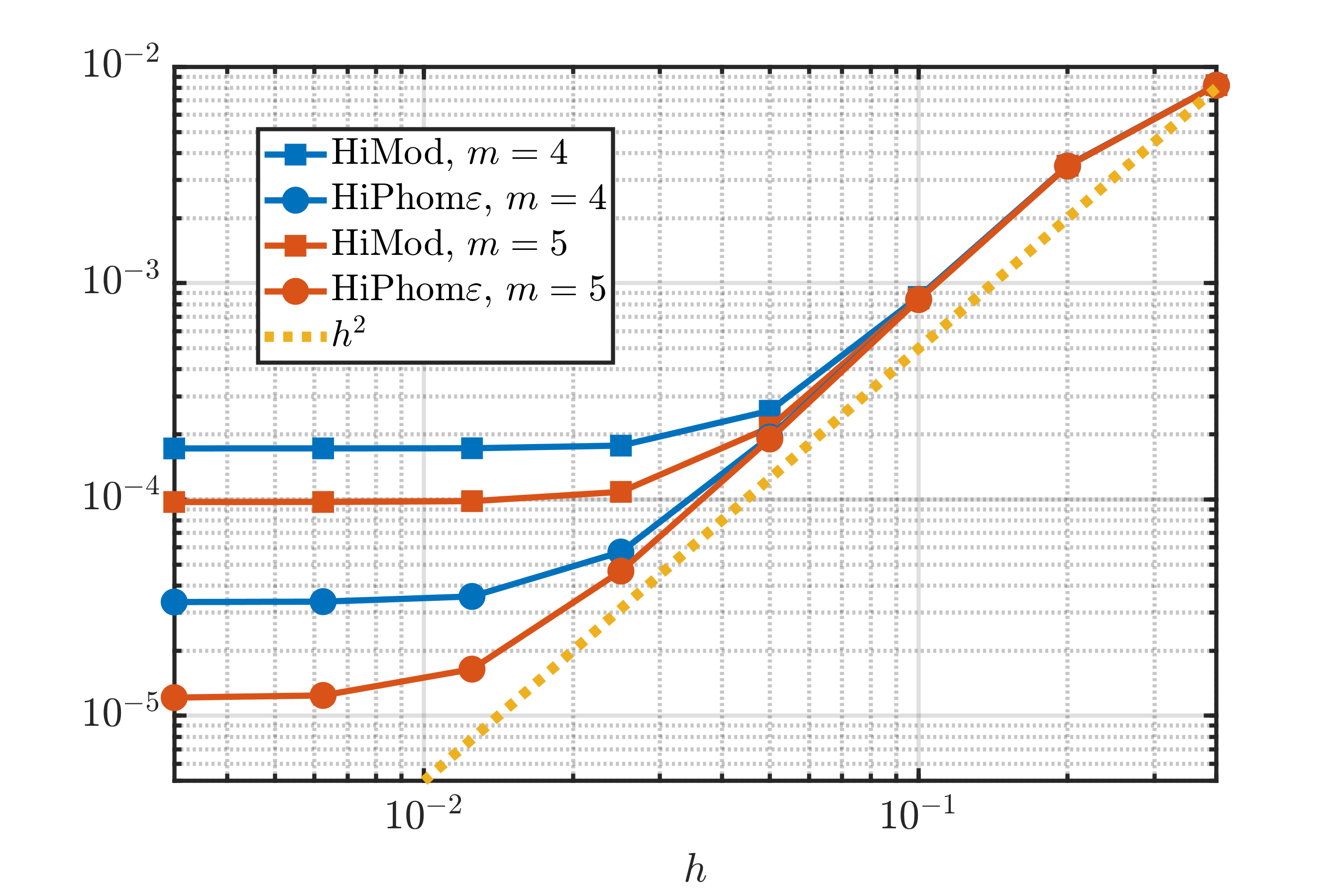}
    \caption{The boundary layer flow test case: $L^2(\Omega_\varepsilon)$-norm of the modelling error associated with the HiPhom$\varepsilon$ and the HiMod approximations as a function of $m$ for different choices of $h$ (left) and as a function of $h$ for different choices of $m$ (right). The dot lines provide a reference trend for the convergence rate.}
    \label{fig:errfin1}
\end{figure}

In Figure~\ref{fig:errfin1} we provide the converge history for the $L^2(\Omega_\varepsilon)$-norm of the HiPhom$\varepsilon$ and of the HiMod modelling errors as a function of the modal index (left), for the choices $h=0.05$ and $h=0.025$ of the discretisation step along the supporting fibre, and as a function of $h$ (right), for two different choices of $m$. 

For the two selected spatial discretisation steps, HiPhom$\varepsilon$ exhibits 
a higher rate of model convergence with respect to the HiMod reduction ($\mathcal{O}(m^{-4})$ versus $\mathcal{O}(m^{-3})$). As discussed in section \ref{sec:erranapoi}. this result implies that the HiPhom$\varepsilon$ modelling error stagnates earlier for lower values of $m$ if compared to the HiMod reduction.
In particular, the maximum accuracy guaranteed by the two reduction procedures is about the same for $h=0.05$. Vice versa, the HiPhom$\varepsilon$ approach outperforms the HiMod reduction for $h=0.025$, with an accuracy equal to $10^{-5}$ when using $m=6$ modes to be compared with an error equal to $9\cdot 10^{-4}$ when employing $11$ educated modal basis functions. In addition, a cross-comparison with the plot in Figure~\ref{fig:errbym}, left highlights a deterioration of the HiPhome convergence rate, together with a monotonic decay of the HiMod error, all the educated modes being now instrumental to the modelling of the transverse dynamics.

As far as the finite element convergence is concerned, we still detect a superior performance of HiPhom$\varepsilon$ when compared with HiMod, with a more striking mismatch with respect to the Poiseuille configuration (see the right panel in Figure~\ref{fig:errbym}). For instance, about one order of accuracy is gained by the HiPhom$\varepsilon$ reduction for $m = 5$ with respect to HiMod. Finally, the quadratic rate of convergence is preserved ($\mathcal{O}(h^2)$), in accordance with the standard finite element theory.
\begin{figure}[h!tbp]
    \hspace*{-0.35cm}\includegraphics[width=0.52\linewidth]{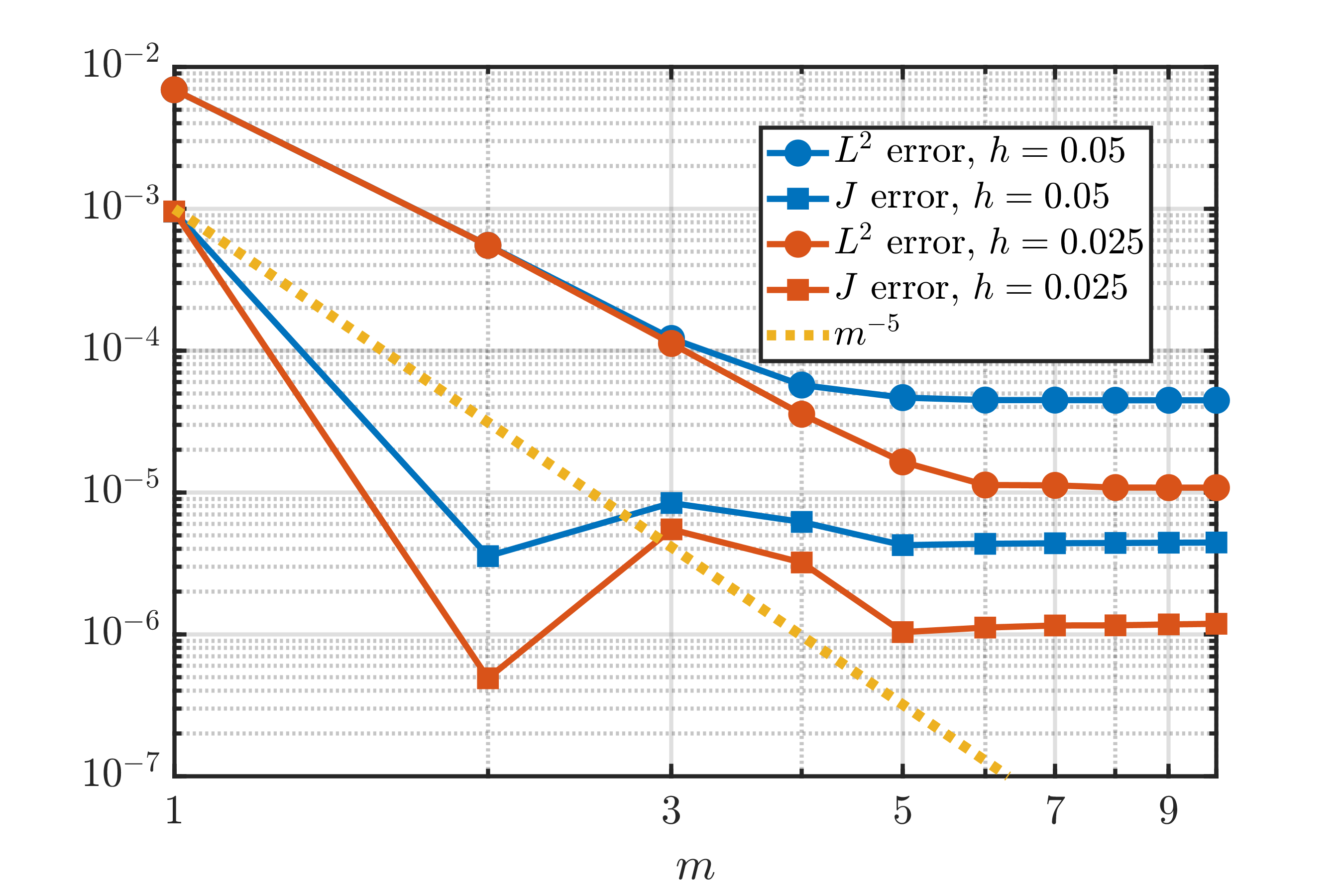}
    \hspace*{-0.35cm}\includegraphics[width=0.52\linewidth]{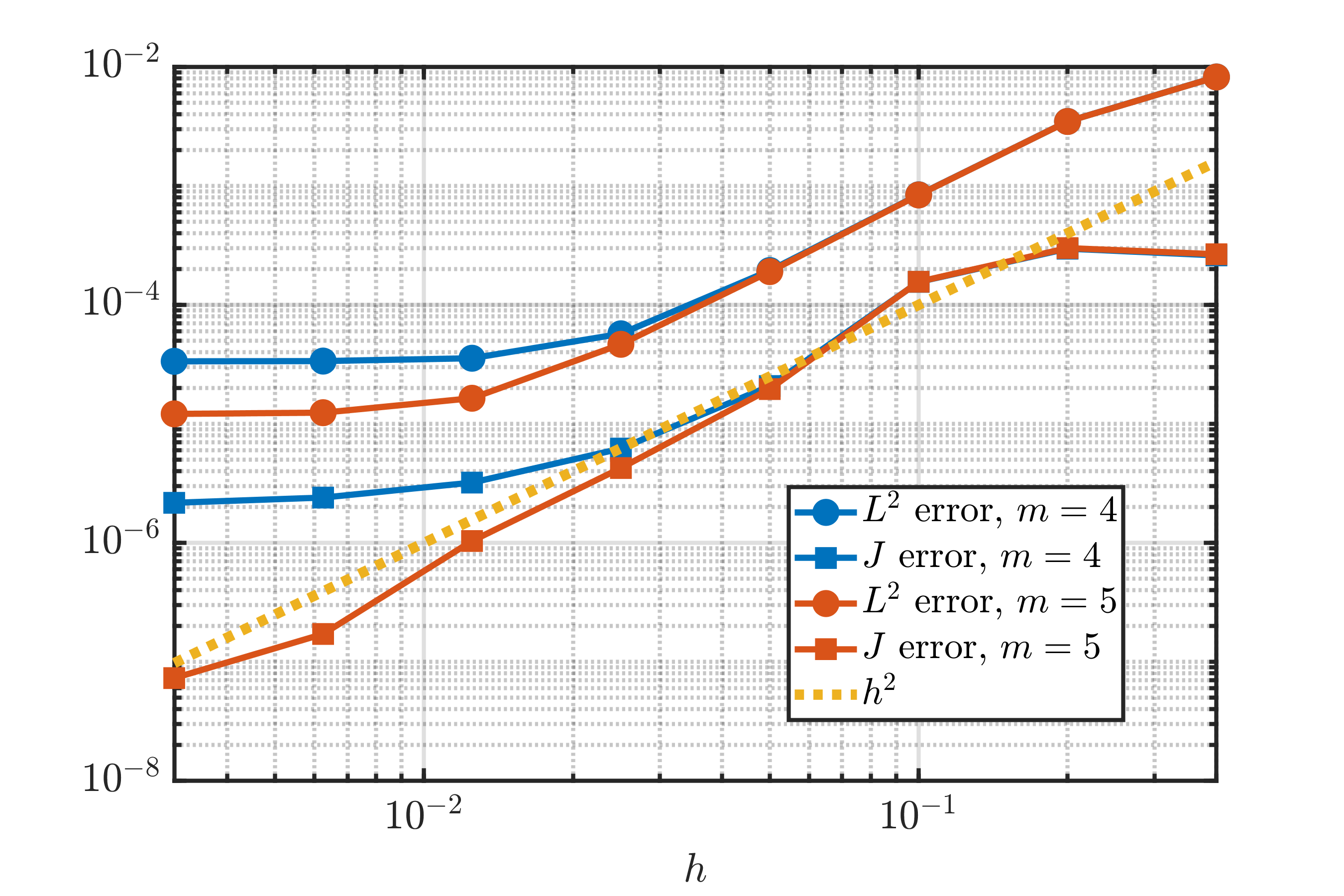}
    \caption{The boundary layer flow test case: 
    error associated with the HiPhom$\varepsilon$ approximation with respect to the $L^2(\Omega_\varepsilon)$-norm and the QoI as a function of $m$ and for different choices of $h$ (left) and as a function of $h$ for different choices of $m$ (right). The dot lines provide a reference trend for the convergence rate.}
    \label{fig:logj2errfirst}
\end{figure}

Finally, in Figure~\ref{fig:logj2errfirst} we investigate 
the capability of the HiPhom$\varepsilon$ reduced solution in reproducing the QoI $\| c \|_{L^2(\Omega_\varepsilon)}$,
as done for the Poiseuille test case. Analogously to Figure~\ref{fig:j2errfirst}, we compare the convergence history for the error
$J$ with the trend of the $L^2(\Omega_\varepsilon)$-norm of the HiPhom$\varepsilon$ modelling error, as a function of the modal index for two different choices of $h$ (panel on the left), and as a function of $h$ for two diverse values of the modal index (panel on the right).\\
Both the panels confirm that the HiPhom$\varepsilon$ model reduction is very sound in reproducing the QoI, provided that at least $3$ HiPhom$\varepsilon$ modes are used to approximate the transverse dynamics, the discrepancy with respect to the  $L^2(\Omega_\varepsilon)$-norm of the modelling error being preserved, in line with results obtained with Poiseuille flow. The rates of convergence are the same as in Figure~\ref{fig:j2errfirst}.

\subsubsection{The unsteady setting}
As the last check, we address a time-de\-pen\-dent variant of the boundary layer flow configuration. To this aim, we select the interval $(0, 0.4)$ as the 
time window of interest, and we use the explicit Euler scheme to discretise the time dependence of $c$ in \eqref{eq:6}, with a time step $\Delta t$ equal to $0.005$. Concerning the other problem data, we preserve the same values as for the steady case, except for the reaction coefficient that we set to zero.

The HiPhom$\varepsilon$ model reduction is performed 
by resorting to the first four HiPhom$\varepsilon$ modes, while exploiting the same discretisation along $\Omega_\varepsilon^{1D}$ as for the steady setting.
Relying on expansion \eqref{HP_exp}, we observe that the HiPhom$\varepsilon$ basis functions do not depend on time. Moreover, from Proposition~\ref{lem:theo1}, it follows that the reactive term is not involved in the differential problem \eqref{eq:chi_i_plus_cond_a}-\eqref{eq:chi_i_plus_cond_b} defining function $\chi^*_i$ (namely, functions $\chi_i$). Thus, the HiPhom$\varepsilon$ modes for the unsteady setting exactly coincide with the ones characterising the steady simulation (i.e., with functions in Figure~\ref{fig:Tc2_data}, left).

\begin{figure}[t]
    \centering
    \hspace*{-.3cm}
    \includegraphics[height=1.55cm]{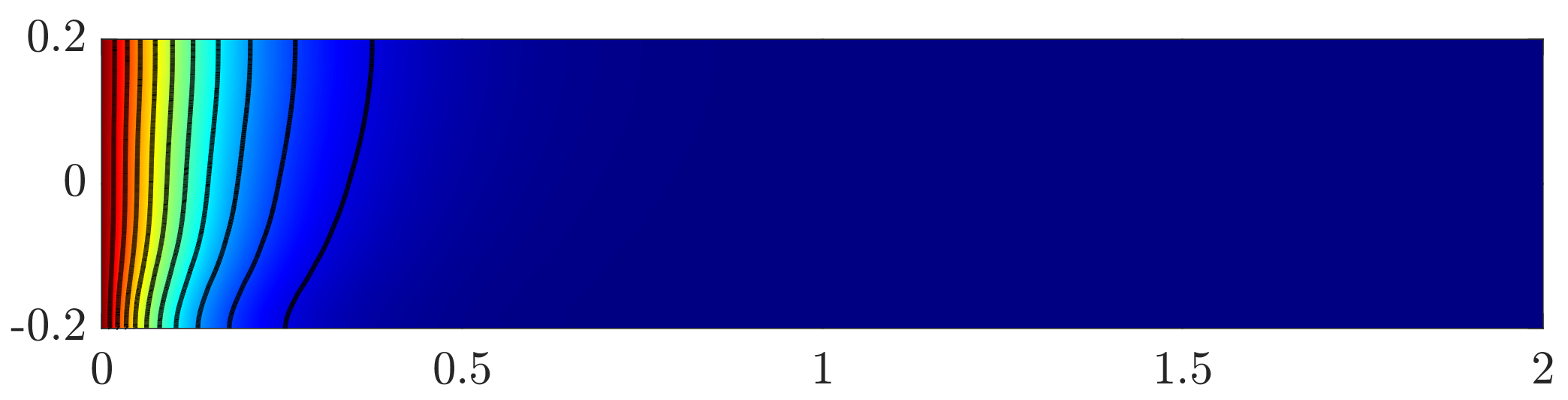}
    \hspace*{-.1cm}
    \includegraphics[height=1.55cm]{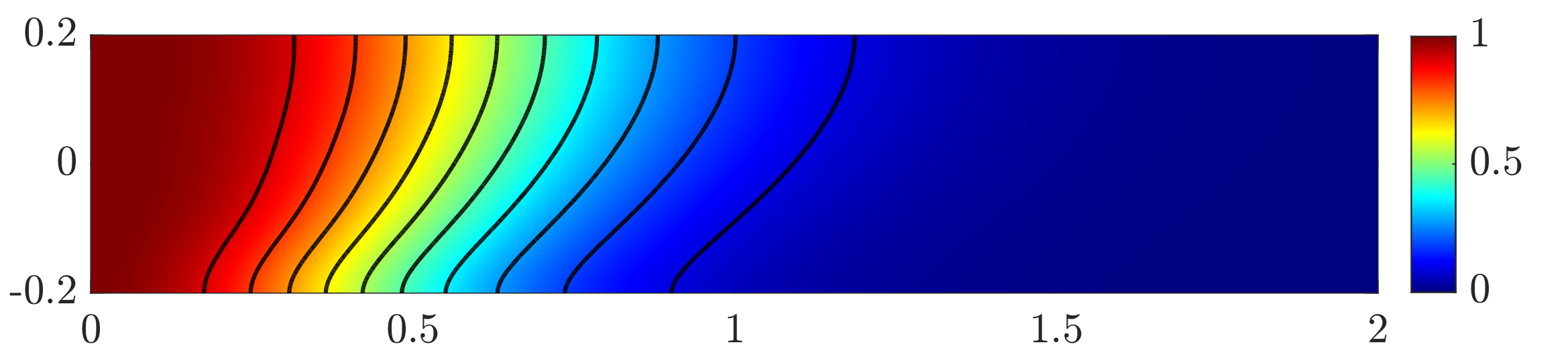}
    \hspace*{-.3cm}
    \includegraphics[height=1.55cm]{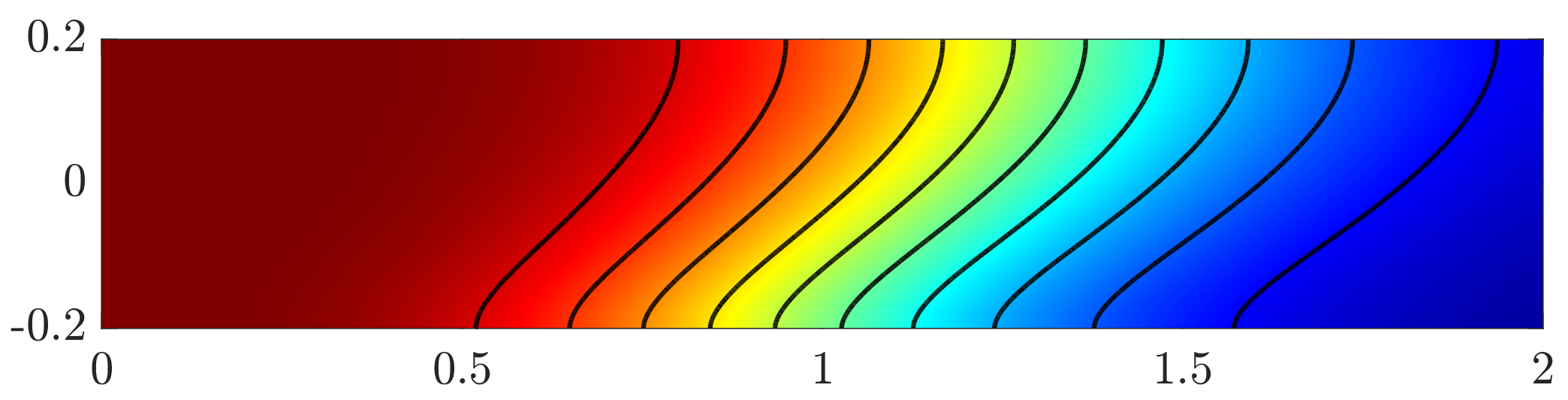}
    \hspace*{-.1cm}
    \includegraphics[height=1.55cm]{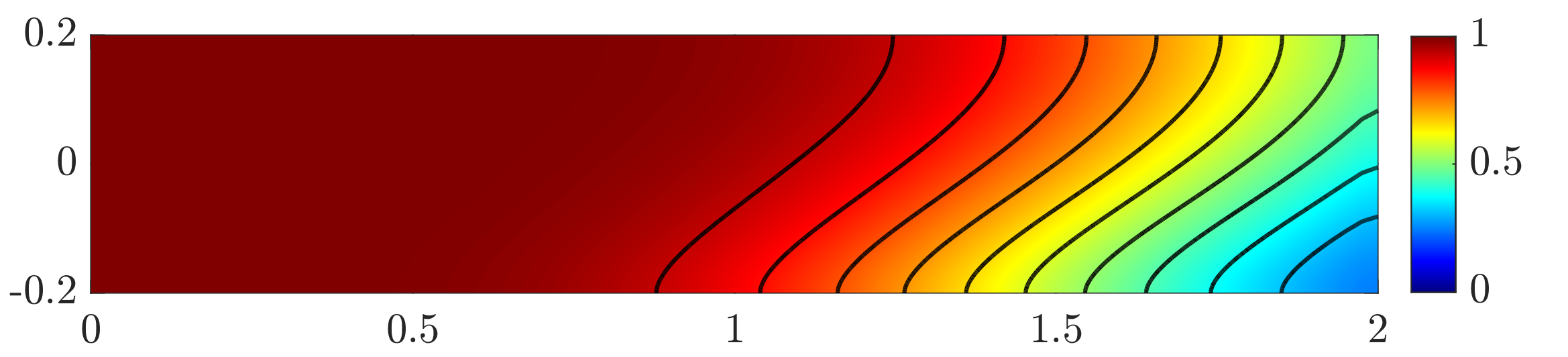}
    \hspace*{-.3cm}
    \includegraphics[height=1.55cm]{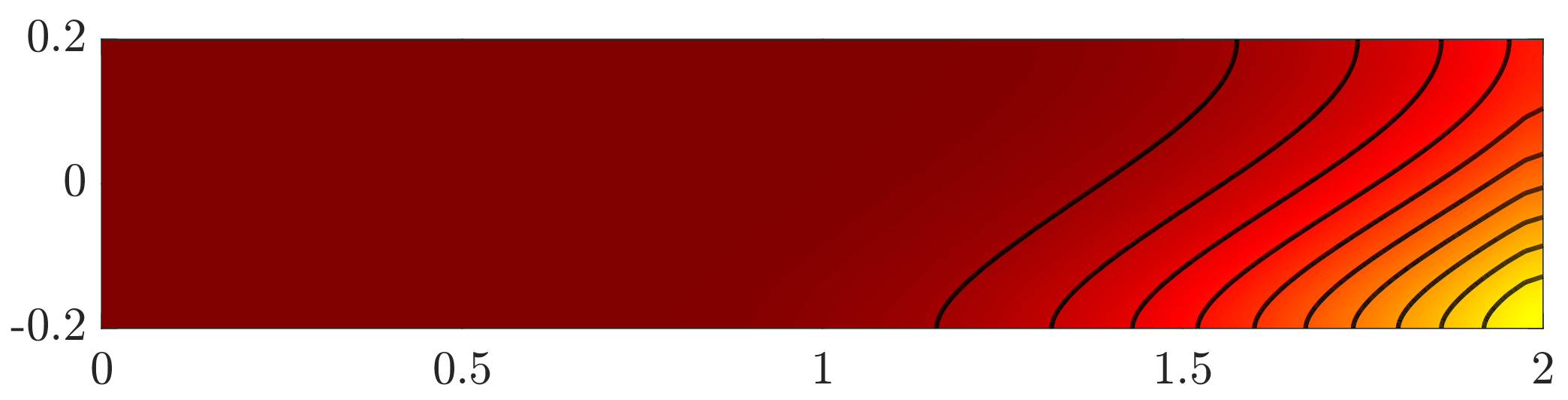}
    \hspace*{-.1cm}
    \includegraphics[height=1.55cm]{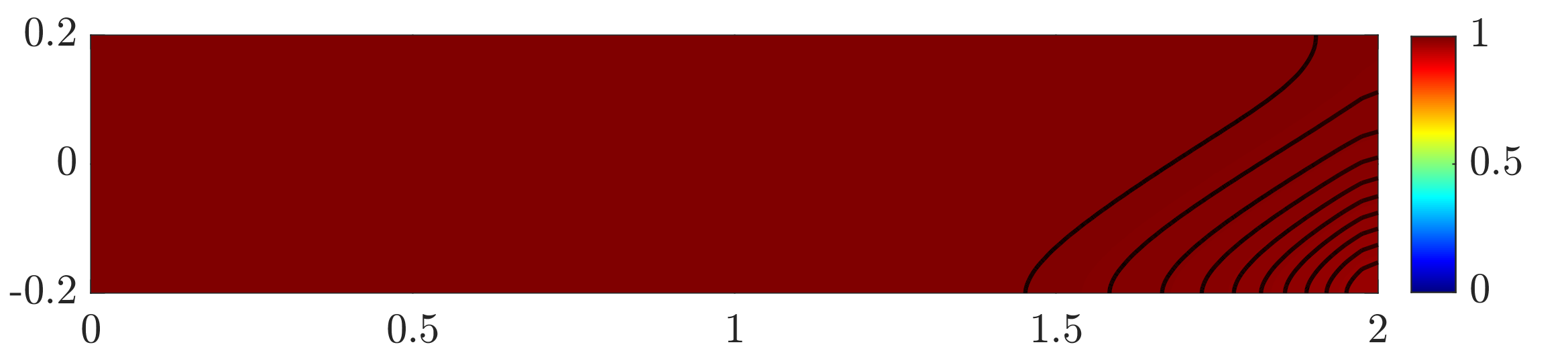}
     \caption{The unsteady boundary layer flow test case: contour plot of the HiPhom$\varepsilon$ approximation at the time instants $t = 0.01$, $0.05$, $0.1$, $0.15$, $0.2$, $0.3$ (top-bottom, left-right).}
  \label{fig:logVel_unsteadyHiPhome}
\end{figure}

\begin{figure}[t]
    \centering
    \hspace*{-.3cm}
    \includegraphics[height=1.55cm]{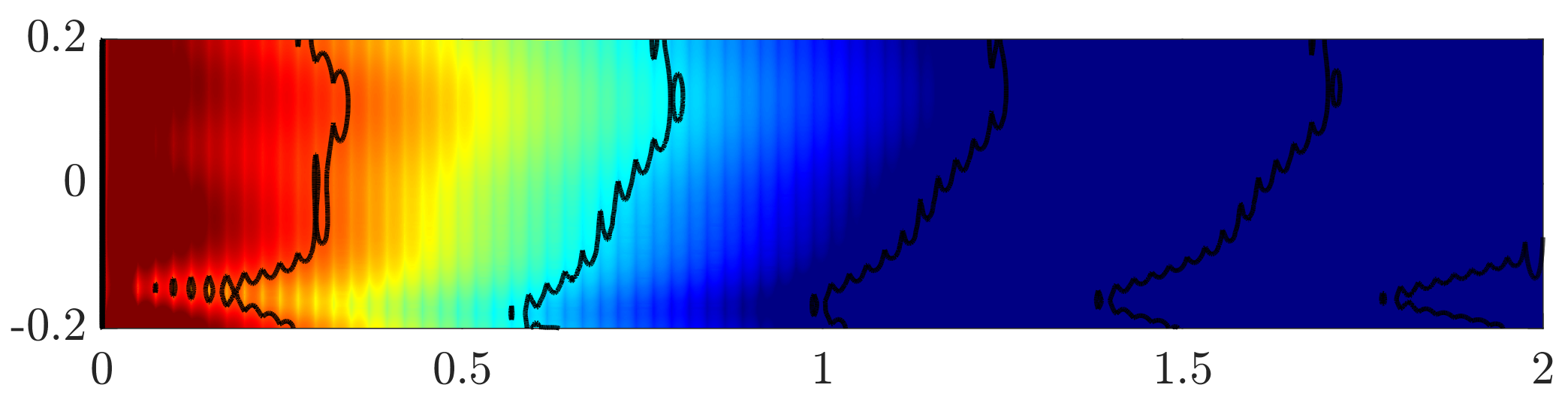}
    \hspace*{-.1cm}
    \includegraphics[height=1.55cm]{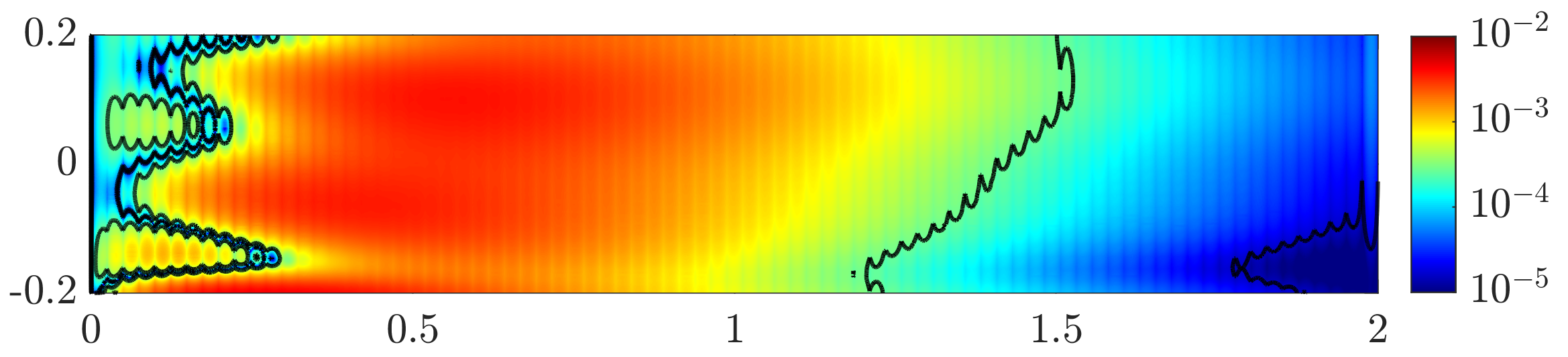}
    \hspace*{-.3cm}
    \includegraphics[height=1.55cm]{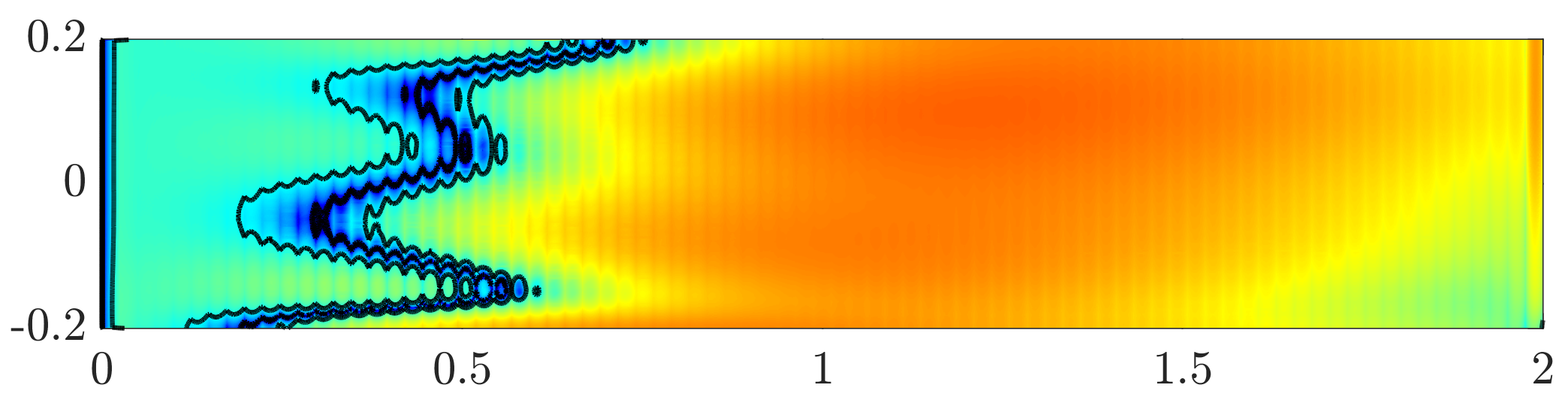}
    \hspace*{-.1cm}
    \includegraphics[height=1.55cm]{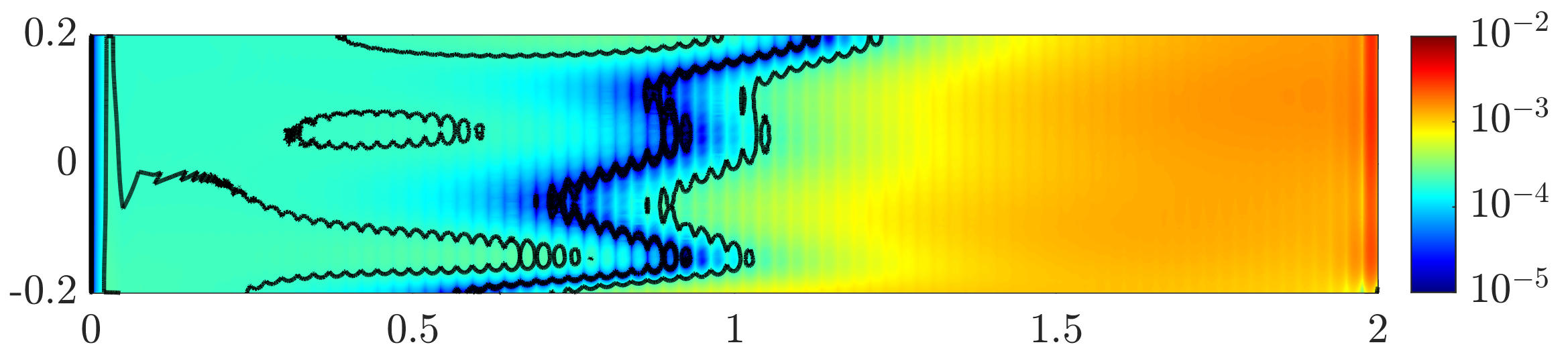}
    \hspace*{-.3cm}
    \includegraphics[height=1.55cm]{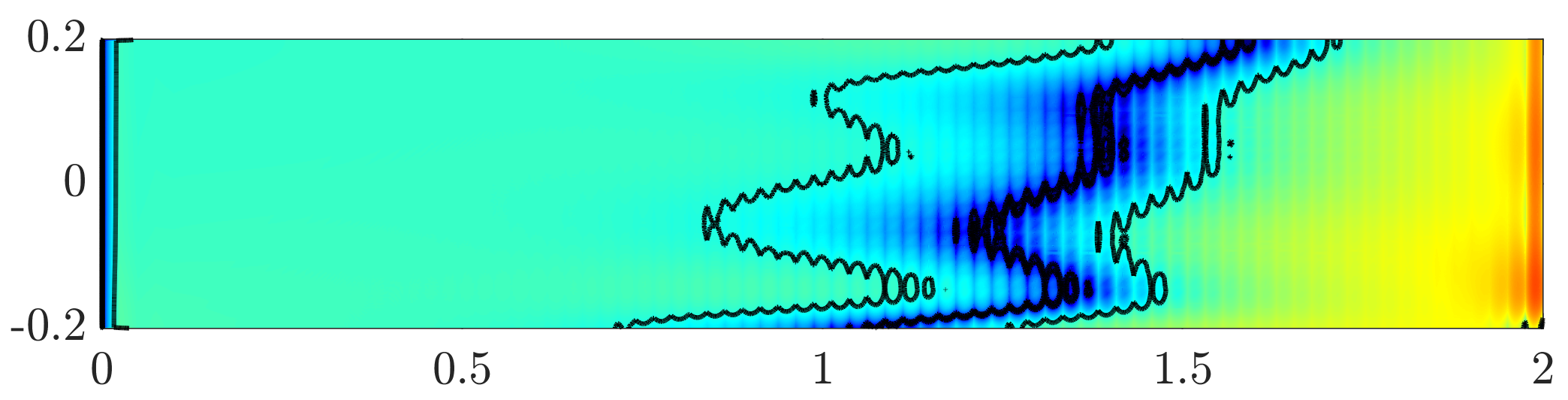}
    \hspace*{-.1cm}
    \includegraphics[height=1.55cm]{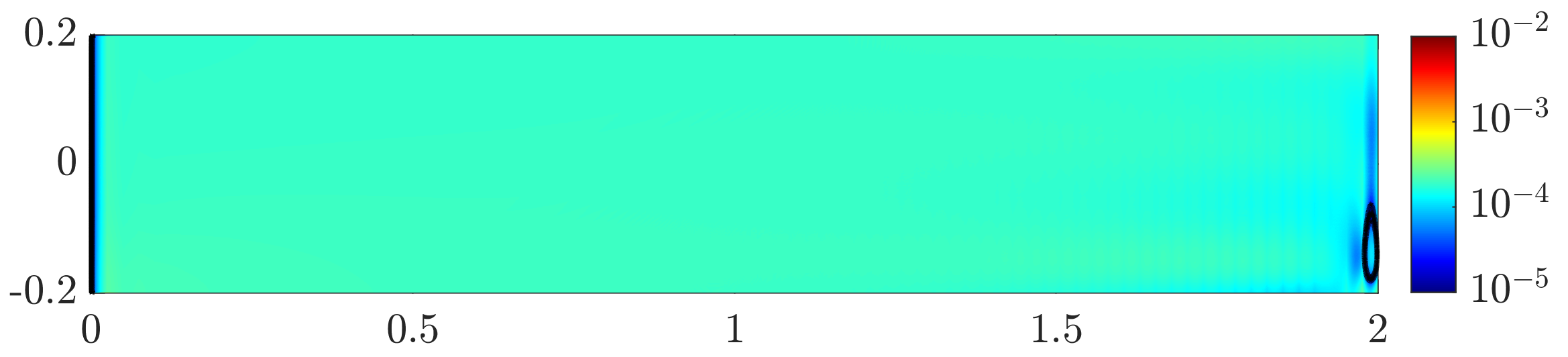}
     \caption{The unsteady boundary layer flow test case: spatial distribution of the absolute modelling error associated with the HiPhom$\varepsilon$ approximation at the time instants $t = 0.01$, $0.05$, $0.1$, $0.15$, $0.2$, $0.3$ (top-bottom, left-right).}
  \label{fig:logVel_unsteadyHiPhomeerr}
\end{figure}

Figure~\ref{fig:logVel_unsteadyHiPhome} shows the contour plot of the HiPhom$\varepsilon$ approximation at $6$ successive times. 
The solute concentration develops until it reaches the steady state around $t=0.3$. The accuracy of the HiPhom$\varepsilon$ discretisation is evaluated by taking as reference solution the linear finite element approximation of $c$ in \eqref{eq:6}, associated with a uniform tessellation of $\Omega_\varepsilon$ consisting of $906640$ triangles, while preserving the time discretisation adopted for the reduced model. As it is highlighted by Figure~\ref{fig:logVel_unsteadyHiPhomeerr}, the spatial distribution of the 
absolute modelling error
follows the propagation of the solute, by reaching a uniform distribution, characterized by an order of magnitude equal to $10^{-4}$, when the steady regime is achieved (bottom-right panel).
%
%
\begin{figure}[t]
    \centering
    \hspace*{-.3cm}
    \includegraphics[height=1.55cm]{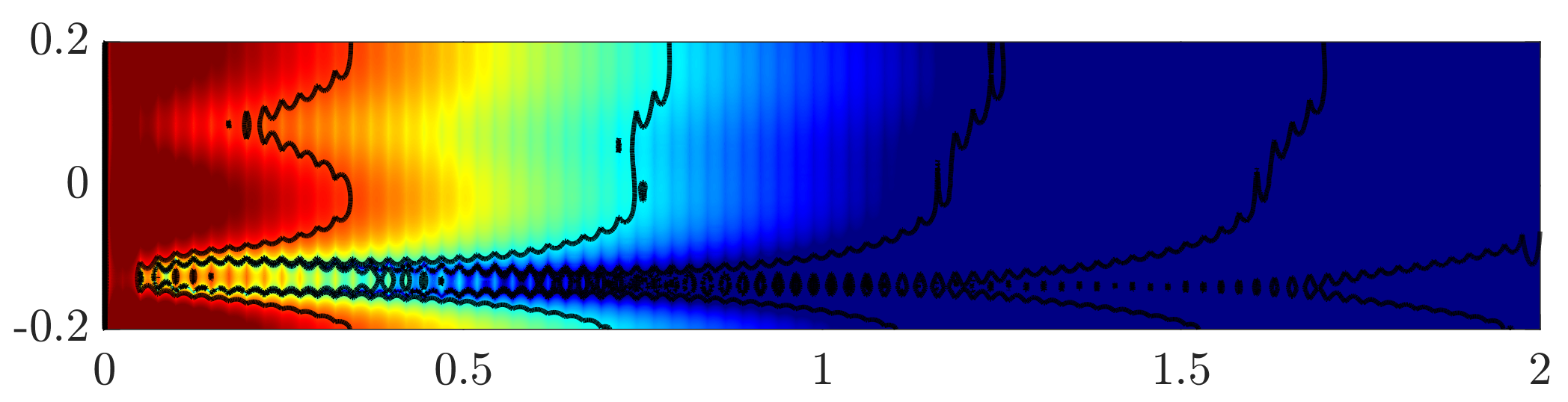}
    \hspace*{-.1cm}
    \includegraphics[height=1.55cm]{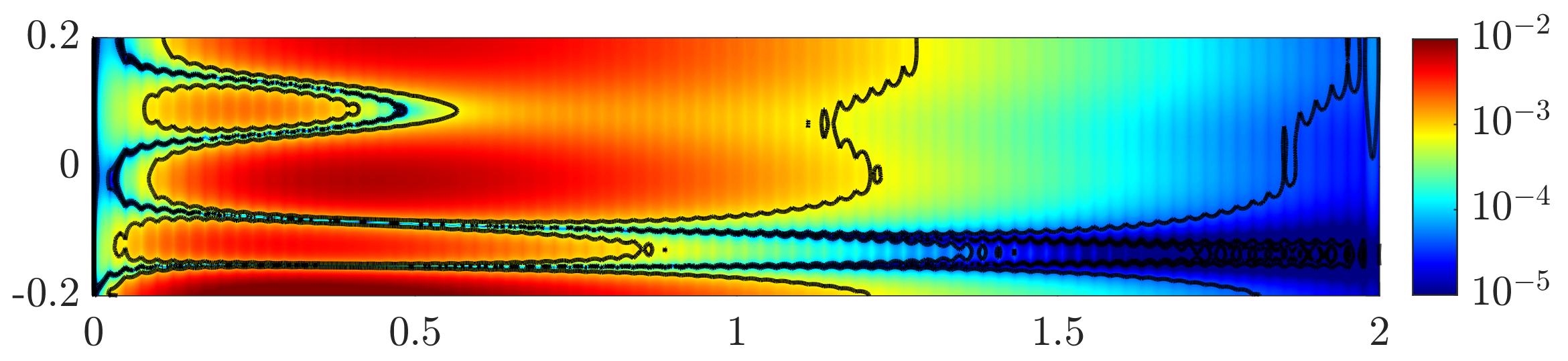}
    \hspace*{-.3cm}
    \includegraphics[height=1.55cm]{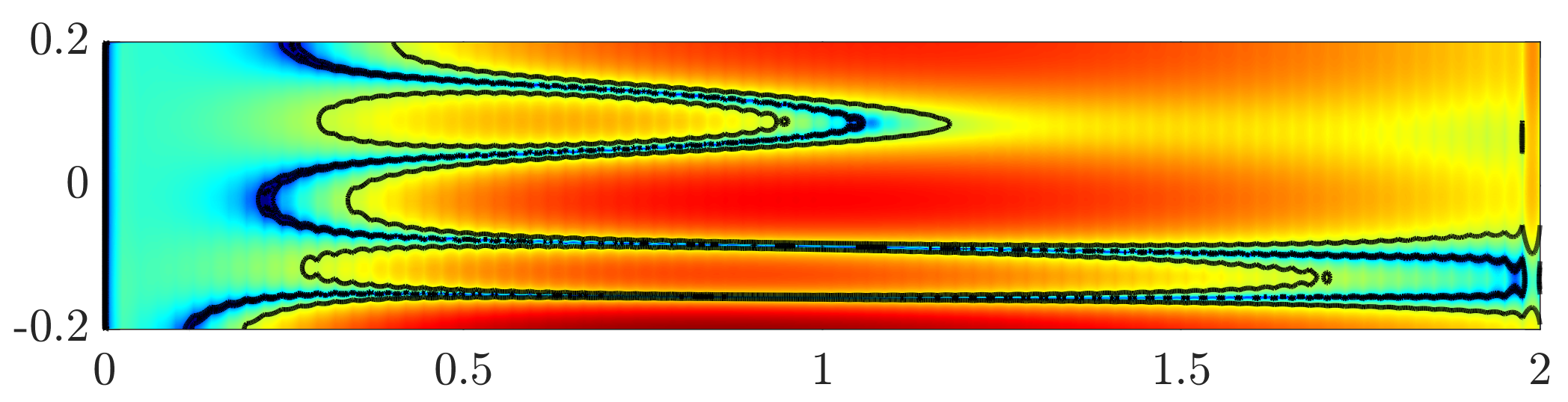}
    \hspace*{-.1cm}
    \includegraphics[height=1.55cm]{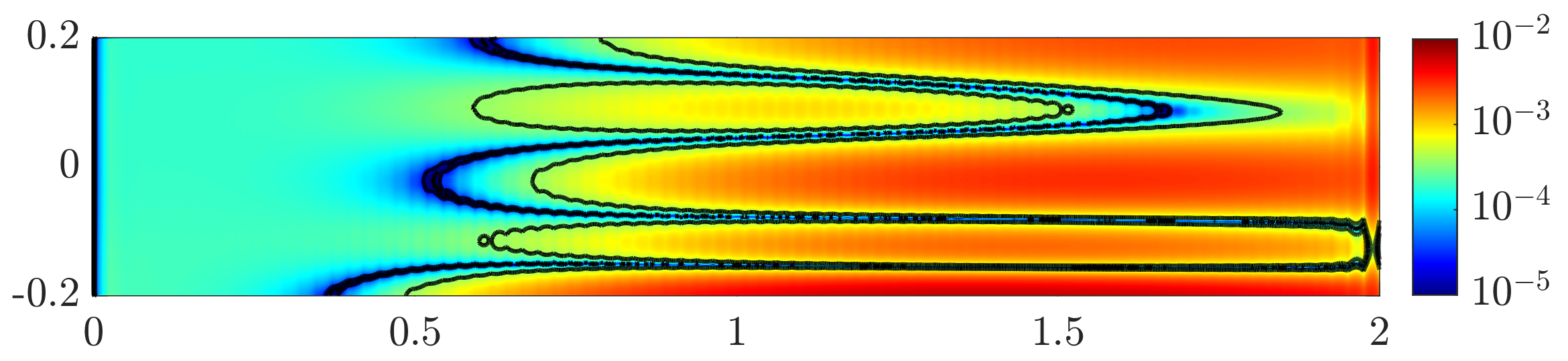}
    \hspace*{-.3cm}
    \includegraphics[height=1.55cm]{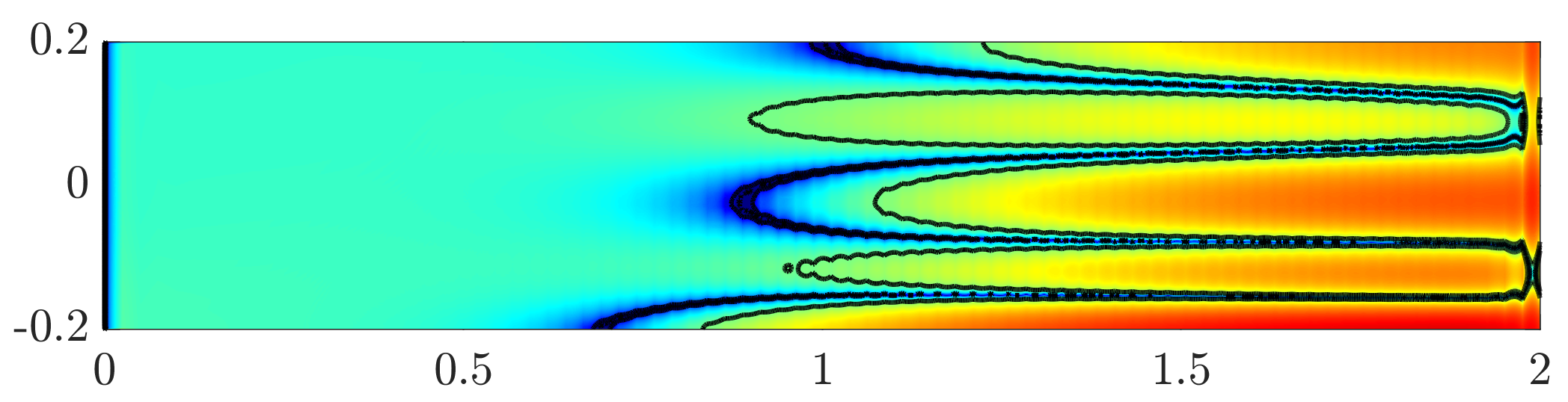}
    \hspace*{-.1cm}
    \includegraphics[height=1.55cm]{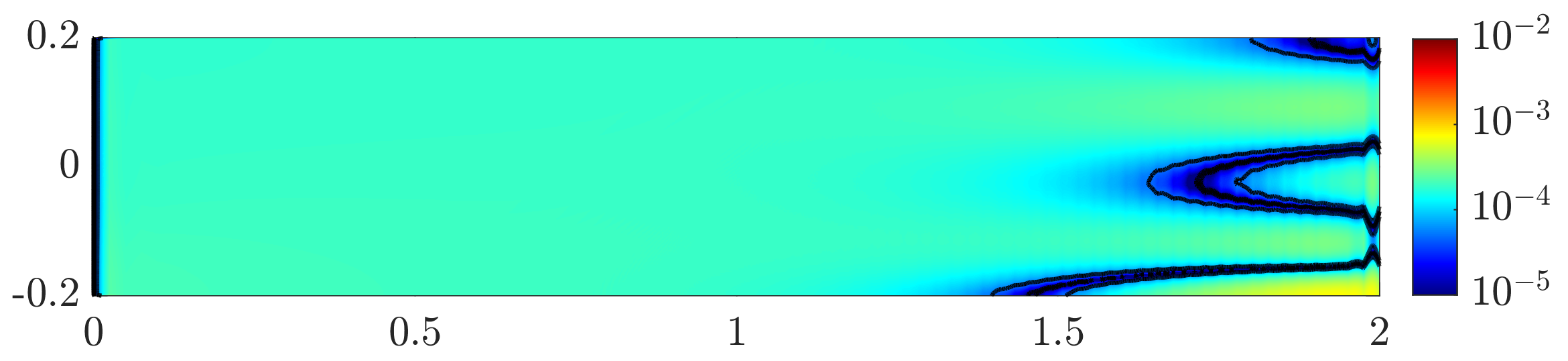}
     \caption{The unsteady boundary layer flow test case: 
     spatial distribution of the absolute modelling error associated with the HiMod approximation at the time instants $t = 0.01$, $0.05$, $0.1$, $0.15$, $0.2$, $0.3$ (top-bottom, left-right).}
  \label{fig:logVel_unsteadyHiModerr}
\end{figure}

\begin{figure}[!htb]
    \centering
    \includegraphics[width=0.6\linewidth]{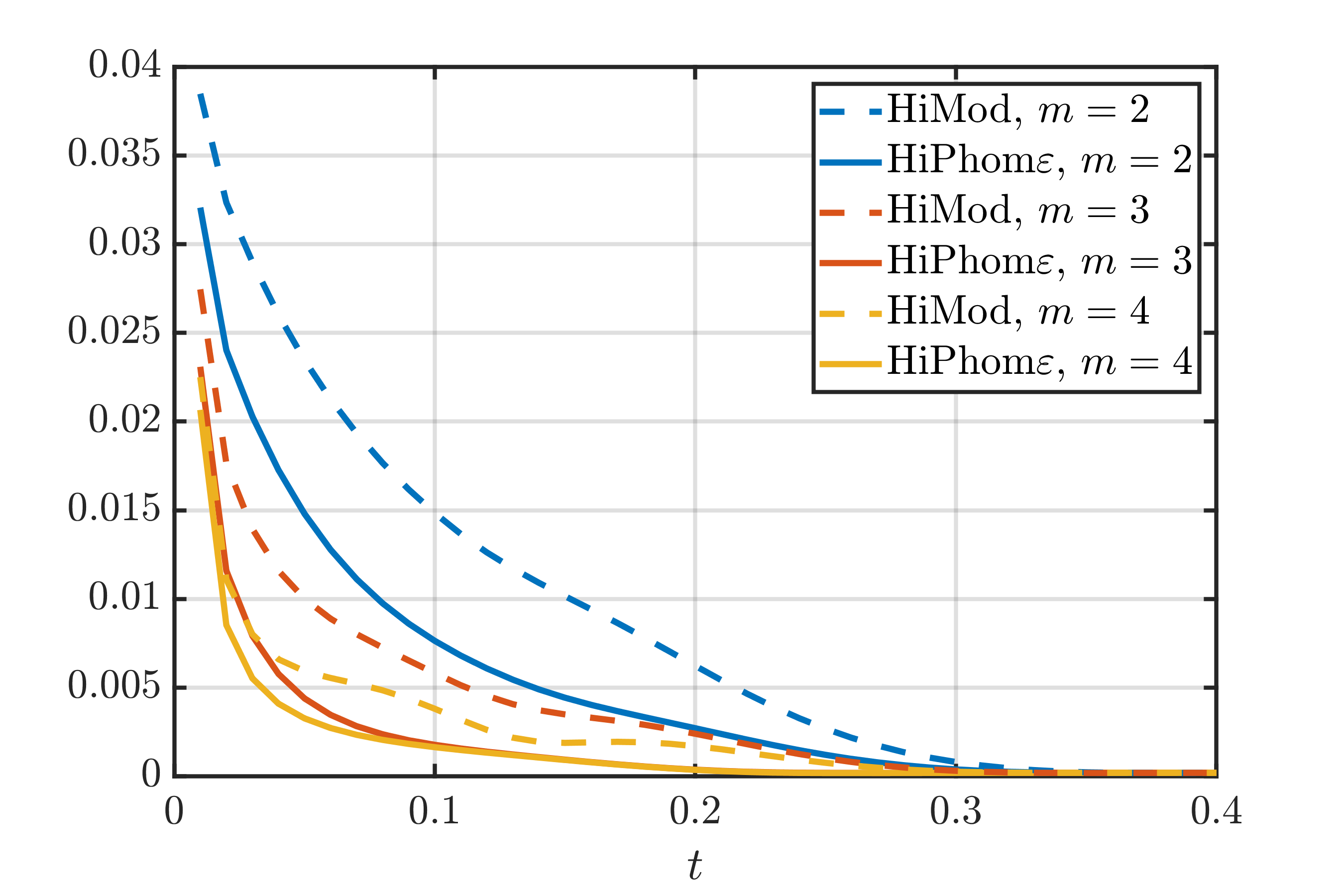}
    \caption{The unsteady boundary layer flow test case: time evolution of the $L^2(\Omega_\varepsilon)$-norm of the modelling error associated with the HiPhom$\varepsilon$ (solid lines) and the HiMod (dashed lines) approximations, for different choices of $m$.}
    \label{fig:errInTime}
\end{figure}

For comparison purposes, we approximate problem \eqref{eq:6} also through the HiMod reduction, by adopting the same spatial discretisation along $\Omega_\varepsilon^{1D}$ as well as the same time discretisation as for the HiPhom$\varepsilon$ approach, while modelling the transverse dynamics with the first $4$ educated modal basis functions derived by setting $\mathcal L_s=D\partial_{yy}$ in \eqref{eq:SLE1}. From a qualitative viewpoint, the HiPhom$\varepsilon$ and the HiMod solutions are practically indistinguishable. To provide a more quantitative cross-comparison between the two reduced solutions, we analyze the spatial of the absolute HiMod modelling error in Figure~\ref{fig:logVel_unsteadyHiModerr}. 
The corresponding panels in Figures~\ref{fig:logVel_unsteadyHiPhomeerr} and~\ref{fig:logVel_unsteadyHiModerr} show that 
the HiMod reduction guarantees a lower accuracy with respect to the HiPhom$\varepsilon$ approach, in particular at intermediate times, i.e. at time levels $t \in [0.1, 0.2]$ (see Figures \ref{fig:logVel_unsteadyHiPhomeerr} - \ref{fig:logVel_unsteadyHiModerr} middle row of plots and bottom left). The performance of the HiPhom$\varepsilon$ reduction is confirmed by Figure~\ref{fig:errInTime} that displays the time evolution 
of the $L^2(\Omega_\varepsilon)$-norm of the modelling error associated with both the HiPhom$\varepsilon$ and the HiMod discretisations, with three different choices of the modal index $m$.
Independently of the selected reduction process, the modelling error diminishes when approaching the end of the time window, with an evident improvement of the HiPhom$\varepsilon$ over the HiMod approximation in the pre-asymptotic dynamics. Finally, we remark that the discrepancy between the HiPhom$\varepsilon$ and the HiMod procedures reduces when we increase the number of modal functions employed in the modelling of the solute transverse dynamics.
%

\section{Conclusions}\label{sec:conclusions}

In this work, we have extended the two-scale expansion for an advection-diffusion-reaction problem in long rectangular domains, to arbitrary orders and interpreted it in the spirit of a Hierarchical Model (HiMod) reduction, by exploiting the separable structure of the asymptotic solution. This lead us to define a new approach, that we named HiPhom$\varepsilon$ (HIgh-order Projection-based HOMogEnisation), which, indeed, combines HiMod with the asymptotic expansion underpinning formal homogenisation theory.

We have verified the potential of this idea as a versatile, efficient model order reduction method, with problem-specific basis functions that adapt to the underlying partial differential equation. The
numerical results for the transport in a rectilinear domain are very encouraging, capturing with high accuracy the evolution of the concentration profile in the whole domain, including the regions near the boundaries where standard two-scale homogenisation fails. This approach can therefore be considered both as an enriched one-dimensional numerical method, as well as an extended projection-based high-order homogenisation. In fact, by solving a system of coupled one-dimensional problems, it becomes possible to approximate the pre-asymptotic dynamic behaviour of the system and thus increase the accuracy of standard homogenised solutions. In particular, the desired accuracy can be achieved by properly balancing the discretisation along the mainstream and the number of modes.

Future works include a more extensive theoretical analysis of the approximation properties of the corrector functions and the resulting separable expansions as well as the extension of HiPhom$\varepsilon$ towards three-dimensional curvilinear domains, more complex boundary conditions and flow profiles that require non-standard two-scale asymptotics, as in \cite{municchi2020generalized}.


\section*{Acknowledgments}
SP acknowledges the European Union’s Horizon 2020 research and innovation programme under the Marie Sklodowska-Curie Actions, grant agreement 872442 (ARIA, Accurate Roms for Industrial Applications).
MI thanks 
IES-R3-170302, Royal Society International Exchanges 2017 Round 3.

\bibliographystyle{siamplain}
\bibliography{refnodoi}

\end{document}